\documentclass{amsart}
\usepackage{amssymb, amsmath}
\usepackage{bbm}
\usepackage{calrsfs}
\usepackage{mathrsfs}
\font\fndef=cmbxti12
\font\sps=cmti5

\def\largerightarrow{   -\negthinspace\negthinspace -\negthinspace\negthinspace
                                       \negthinspace\longrightarrow }
                                       
\begin{document}
\title{\Large Regular Algebraic K-Theory for groups -- Part I}\author{U. Haag}
\date{\today}
\maketitle
\hfil \Large{\bf{(1)20 years after -- instead of an introduction.}}\hfil
\bigskip\bigskip\bigskip
\par\noindent
\Large
The work on Regular Algebraic $K$-theory for groups presented below was begun in 2004 and completed around 2010. 
The principal reasons to start this investigation were twofold: for one thing the goal was to find a convenient intermediate tool between ordinary group homology on one side and topological (operator algebra) $K$-theory in order to find a new approach to the (rational) Novikov Conjecture on higher signatures. The second goal was to find a natural "algebraic" extension of algebraic $K$-theory to higher dimensions. Both goals were achieved to some extent by this exposition, c.f. \cite{H3}. Of course a workable proof for the Novikov Conjecture in a really general setting is still at large, while on the other hand the precise relation of Regular Algebraic $K$-theory with Quillen's algebraic $K$-theory and related homology theories like group homology, topological $K$-theory etc. is still very much obscure, so the question remains where to place this theory among the bouquet of existing $K$-theories. Algebraic $K$-groups made their first appearance in the mathematical literature in 1904 in the form of central group extensions studied by I. Schur in his exposition "Über die Darstellung der endlichen Gruppen durch gebrochene lineare Substitutionen" which he termed {\it multipliers}. The theory below may be seen as a natural extension of this notion to higher dimensions. The term {\it higher Schur multipliers} often is used for the higher group homology of a group $\, G\, $. This is because in dimension $2$ one has an isomorphism of the algebraic $K$-group and the corresponding homology group of $\, G\, $. In dimensions $\, \geq 3\, $ there still is a natural map from the algebraic $K$-group into the corresponding homology group which however in some cases is far from being an isomorphism (injective or surjective). This map is given an studied in Appendix B of part III of this treatise (still unpublished).    
The material presented here is still in a very raw form and needs to be refined and completed in many ways. Nevertheless the author wants to give the interested reader an outlook on the established results since it is uncertain when or whether the theory arrives at a more accomplished state. 
Some questions that remain to be tackled include the formulation of a dual "Co-$K$"-theory for which there are some ideas but no more than that, the computation of some nontrivial examples for the higher $K$-groups, and as mentioned the clarification of its relation with (rational) group homology and the operator $K$-theory of the group $C^*$-algebra. Part III of this article series contains partial answers for the last problem and gives a full computation of the higher $K^J$-groups in case of a finite group (they are trivial, which is a both striking and disturbing fact since it raises the question whether the higher $K$-groups might be trivial in general). Of course one may argue that three algebraic $K$-groups are just about enough to distinguish two finite groups and the author knows no example where two different finite groups have the same algebraic $K$-groups. This exemplifies the need for a computation of some nontrivial examples. 
A different branch of evolution would be to try to extend the definition in a meaningful way to cover the case of topological groups or monoids.
\bigskip\bigskip
\par\noindent    
\hfil \Large{\bf{Section I -- Basic constructions}} \hfil 
\par\bigskip\bigskip\noindent
{\bf Definition 1.} A subgroup $\, N\subseteq F\, $ will be called { \it full } iff $\, N\, =\, [ N , F ]\, $.
\par\bigskip\noindent
{\bf Definition 2. (a)} Given a group $\, N\, $ and a group $\, F\, $ with a specified map 
$\, \alpha : F \rightarrow {\it Aut} ( N )\, $ consider the category of such pairs where a morphism is given by a pair of homomorphisms $\, \varphi : N \rightarrow N'\, $ and 
$\, \widetilde\varphi : F \rightarrow F'\, $, such that $\, \varphi\, $ is compatible with 
$\, \widetilde\varphi\, $ in the sense that 
$\, \varphi ( {\alpha}_z ( x ))\, =\, {\alpha}'_{\widetilde\varphi ( z )} ( \varphi ( x ))\, $. We will also require that there is a specified map $\, \iota : N \rightarrow F\, $ such that the composition with $\, \alpha\, $ is equal to the natural map $\, N \rightarrow {\it Inn} ( N )\, $. In particular if $\,\iota\, $ is an inclusion it identifies $\, N\, $ with a normal subgroup of $\, F\, $ which acts on $\, N\, $ in the obvious way (by the adjoint representation). Then the pair $\, ( N , F )\, $ will be called {\it normal}.
\par\medskip\noindent
{\bf (b)} Given a pair $\, ( \widetilde N , \widetilde F )\, $ together with a surjective map onto 
$\, ( N , F )\, $ the group $\, \widetilde N\, $ will be called an $F$-central extension of $\, N\, $ if the action of $\, \widetilde F\, $ on $\, \widetilde N\, $ reduces to an action of $\, F\, $ and this action is trivial on the kernel of the map $\, \widetilde N \rightarrow N\, $. (In particular $\,\widetilde N\, $ is a central extension of $\, N\, $.)
\par\medskip\noindent
{\bf (c)} With assumptions as above the $F$-central extension $\, ( \widetilde N , F )\, $ is called 
{\it universal } iff the following conditions are satisfied :
\par\smallskip\noindent
(i) For any morphism $\, ( N , F ) \rightarrow ( N' , F' )\, $ and any $F'$-central extension $\,\widetilde N'\, $ of $\, N'\, $ there is a map $\, \widetilde N \rightarrow \widetilde N'\, $ over $\, N \rightarrow N'\, $ compatible with $\, F \rightarrow F'\, $ and uniquely determined on the $F$-central kernel of 
$\, \widetilde N \rightarrow N\, $, which is supposed to lie in the subgroup 
$\, [ \widetilde N , F ]\, $.
\par\smallskip\noindent
(ii) If $\, (\widetilde {\widetilde N} , F )\, $ is an $F$-central extension of $\, ( \widetilde N , F )\, $, then there is a splitting $\, ( \widetilde N , F ) \rightarrow ( \widetilde {\widetilde N} , F )\, $.
\par\medskip\noindent
An $F$-central extension which only satisfies part (i) of the Definition above will be called 
{\it semiuniversal }.
\par\bigskip\noindent
{\bf Definition 3. (a)} For any group $\, F\, $ let $\, U_F\, $ denote the free group on the generators 
$\, \left\{ u_x\, \vert\, x\neq 1 , x \in F\right\}\, $ and let $\, J_F\, $ denote the kernel of the natural evaluation map $\, U_F \rightarrow F\, $.
\par\medskip\noindent
{\bf (b)} For $\, N\subseteq F\, $ a normal subgroup denote $\, U_{N , F }\, $ the kernel of the map 
$\, U_F \rightarrow U_{F / N}\, $ induced by $\, F \rightarrow F / N\, $ and let $\, J_{ N , F }\, $ denote the kernel of the map $\, U_{ N , F } \rightarrow N\, $ which is the restriction of the map 
$\, U_F \rightarrow F\, $.
For $\, n \geq 2\, $ define inductively $\, U^n_F := U_{U^{n-1}_F}\> ,\> 
U^n_{N , F} := U_{J^{n-1}_{N , F}\, ,\, U^{n-1}_F}\> $, 
$\> J^n_{N , F} := J_{J^{n-1}_{N , F}\, ,\, U^{n-1}_F}\, $.
\par\medskip\noindent
{\bf (c)} Also put $\, {\mathcal J}_{N , F} = J_{N , F} \cap [\, U_{N , F}\, ,\, U_F\, ]\, ,\, 
{\mathcal J}^n_{N , F} := J^n_{N , F} \cap [\, U^n_{N , F}\, ,\, U^n_F\, ]\, $ and 
$\, {\mathcal R}_{N , F} := [\, J_{N , F}\, ,\, U_F\, ]\, +\, [\, U_{N , F}\, ,\, J_F\, ]\, $. 
\par\bigskip\noindent
{\bf Lemma 1.} The group $\, U_{N , F }\, $ is free on the generators
\smallskip
$$ \left\{ u^{k_1}_{{\overline x}_1}\cdots u^{k_n}_{{\overline x}_n}\, u_x\, u^{-1}_{\overline x }\, 
u^{- k_n}_{{\overline x}_n} \cdots u^{- k_1}_{{\overline x}_1}\, \vert\, x , x_1 , \cdots , x_n \in F\,\right\}  $$
\smallskip
where $\, \left\{ \overline x \right\}\, $ is some fixed set of representatives of the classes $\, x N\, $ in 
$\, F\, $ and $\, k_1 , \cdots , k_n \in \mathbb Z\, $. One makes the convention that $\, \overline x = 1\, $ if $\, x \in N\, $ and puts $\, u_{\overline x} = 1\, $ in this case.
\par\bigskip\noindent
{\it Proof. }\quad The proof is left to the reader\qed
\par\bigskip\noindent
{\bf Theorem 1.}\quad Let $\, N\, $ be a full subgroup of the perfect group $\, F\, $. Then the pair 
$\, ( N , F )\, $ admits a universal $F$-central extension.
\par\bigskip\noindent
{\it Proof.}\quad We claim that the extension
\smallskip
$$ {\overline N}^{ _{\sps F}}\, =\, {\left[\, U_{N , F}\, ,\, U_F\, \right] \over 
             \left[\, J_{N , F}\, ,\, U_F\,\right]\, +\, \left[\, U_{N , F}\, ,\, J_F\,\right] }     $$
\smallskip
is a universal $F$-central extension of $\, N\, $. Note that it is $F$-central. If
\smallskip 
$$ ( N , F )\buildrel ( {\varphi }_N , \varphi )\over\largerightarrow ( N' , F' )  $$
\smallskip
is a morphism and $\, ( \widetilde N' , F' )\, $ a $F'$-central extension of $\, N'\, $ define a map 
$\, \widetilde {\varphi }_N : U_{N , F}\, \longrightarrow \widetilde N'\, $ compatible with the map 
$\, U_F\rightarrow F\buildrel\varphi\over\longrightarrow F'\, $ by sending $\, u_x u^{-1}_{\overline x}\, $ to an arbitrary preimage of $\, {\varphi }_N ( x {\overline x}^{-1} )\, $ in $\, \widetilde N'\, $ and inductively
\smallskip 
$$ \widetilde {\varphi }_N ( u^{\pm }_{\overline x}\, e\, u^{\mp }_{\overline x}\, )\, =\, 
\varphi ( \overline x )^{\pm }\, \widetilde {\varphi  }_N ( e )\, \varphi ( \overline x )^{\mp}  $$
\smallskip
if $\, \widetilde {\varphi }_N ( e ) \, $ is already defined. Since $\, \widetilde N'\, $ is $F'$-central 
$\, \widetilde {\varphi }_N\, $ maps $\, [\, J_{N,F} ,\, U_F\, ]\, $ and $\, [\, U_{N,F} ,\, J_F\, ]\, $ into 1, so that restricting to 
$\, [\, U_{N , F} ,\, U_F\, ]\, $ one gets a well defined map 
$\, ( \overline N , F )\longrightarrow (\widetilde N' , F' )\, $ over $\, ( N , F )\longrightarrow (N' , F' )\, $. We claim that $\, \overline N\, $ is a minimal $F$-normal subextension of 
\smallskip
$$ U_{N , F } \over \left[\,\, J_{ N , F }\, ,\, U_F\,\right]\, +\,\left[\, U_{N , F }\, ,\, J_F\,\right]  $$
\smallskip
over $\, N\, $. Indeed, given a smaller subextension $\, N_1\subseteq \overline N\, $ one gets 
$\, [ N_1 , F ] = [ \overline N , F ] = \overline N\, $, since $\, \overline N\, $ is full. Now let 
$\, ( \widetilde {\overline N} , F )\, $ be some $F$-central extension of $\, ( \overline N , F )\, $. Since $\, F\, $ is perfect this is actually an $F$-central extension of $\, ( N , F )\, $. Then one gets a splitting 
$\, ( \overline N , F ) \longrightarrow ( \widetilde {\overline N} , F )\, $ and since $\, \overline N\, $ is full any map $\, ( \overline N , F ) \longrightarrow ( N' , F' )\, $ is necessarily unique. Thus conditions (i) and (ii) of Definition 2 (c) are satisfied and $\, ( \overline N , F )\, $ is universal \qed
\par\bigskip\noindent
{\bf Definition 4.} Let $\, N\, $ be a normal subgroup of a perfect group $\, F\, $. Define the second relative  
$K$-group $\, K_2 ( N , F )\, $ of the pair $\, ( N , F )\, $ to be the kernel of the universal $F$-central extension of the pair $\, ( [ N , F ] , F )\, $ ( note that $\, ( [ N , F ]\, $ is full in $\, F\, $). Clearly, the assignment $\, ( N , F ) \rightsquigarrow K_2 ( N , F )\, $ gives a covariant functor from the category of normal subgroups of perfect groups to abelian groups.
\par\bigskip\noindent
Our goal is to inductively define higher relative $K$-groups of the pair $\, ( N , F )\, $. The rigorous definition will have to wait until section 6, but we will give a short sketch of the steps to be taken here. First assume that $\, N\, $ is full. Then one inductively constructs pairs 
$\, ( R^n_{N , F}\, ,\, U^n_F )\, $, where each pair admits a (relative) universal 
$U^n_F$-central extension of the form 
$\, ( {\overline R}^n_{N , F}\, ,\, U^n_F )\, $ (see the next section for the precise definition) and one defines $\, K_{n+2} ( N , F )\, $ to be the kernel of 
$\, {\overline R}^n_{N , F}\longrightarrow R^n_{N , F}\, $. Then, for arbitrary normal $\, N\subseteq F\, $, one puts $\, K_n ( N , F ) := K_n ( [ N , F ] , F )\, $. We start the induction putting 
$\, ( R^0_{N , F}\, ,\, U^0_F )\, :=\, ( N , F )\, $. Also assume that at each step there is an almost canonical map $\, U_{R^n_{N , F }\, ,\, U^n_F}\buildrel {\alpha }_n\over\longrightarrow {\overline R}^n_{N , F}\, $ which is surjective and lifts the canonical map 
$\, U_{R^n_{N , F}\, ,\, U^n_F}\longrightarrow R^n_{N , F}\, $ (almost canonical means that it is canonical on commutators 
$\, [\, U_{R^n_{N , F}\, ,\, U^n_F}\, ,\, U^{n+1}_F\, ]\, $). Under these assumptions consider the commutative diagram
\bigskip
$$ \vbox{\halign{ #&#&#&#&#\cr
$J^{n+2}_{N , F}$ & $\largerightarrow$ & \hfil $U^{n+2}_{N , F}$\hfil & $\largerightarrow$ & 
$R^{n+1}_{N , F}$ \cr
\hfil $\Biggm\downarrow$\hfil && \hfil $\Biggm\downarrow$\hfil && $\,\,\Biggm\downarrow$ \cr
$J^{n+2}_F$ & \hfil $\largerightarrow$\hfil & \hfil $U^{n+2}_F$\hfil & \hfil $\largerightarrow$\hfil &
$U^{n+1}_F\, =\, U_{R^n_{N , F}\, ,\, U^n_F} \rtimes U_{U^n_F / R^n_{N , F}}$\hfil \cr
\hfil $\Biggm\downarrow$\hfil && \hfil $\Biggm\downarrow$\hfil && $\,\,\Biggm\downarrow$ \cr
$P^{n+1}_{N , F}$ & \hfil $\largerightarrow$\hfil & \hfil $U_{U^{n+1}_F / R^{n+1}_{N , F}}$\hfil & 
\hfil $\largerightarrow$\hfil &
${\overline R}^n_{N , F}\rtimes U_{U^n_F / R^n_{N , F}}\> .$\hfil \cr }}   \leqno\raise54pt\hbox{\, (1)}   $$
\par\bigskip\noindent
This is to be understood in the sense that one has already constructed the universal $U^n_F$-central extension $\, {\overline R}^n_{N , F}\, $ getting an almost canonical map as indicated. Then  define 
$\, R^{n+1}_{N , F}\, $ to be the kernel of this map and the rest of the diagram is obtained by completion. One needs to show that 
$\, R^{n+1}_{N , F}\, $ admits a (relative) universal $U^{n+1}_F$-central extension. This is the content of Theorem 2. For its proof we need a short series of Lemmas examining the structure of the subgroups 
$\, J_{N , F}\, $ and $\, J_F\, $. We begin with
\par\bigskip\noindent
{\bf Lemma 2.}\quad Let $\, J_F\, $ be the kernel of the map $\, U_F \rightarrow F\, $. Then $\, J_F\, $ is generated freely by the set of elements
\smallskip
$$  \left\{\, u_x u_y u^{-1}_{xy}\, ,\, u_x u_{x^{-1}}\,\vert\, x , y \in F\,\right\}\>   .   $$
\par\bigskip\noindent
{\it Proof.}\quad Suppose given a relation in the generators of $\, J_F\, $ as above. One can assume it to be cyclically reduced and to contain no proper subrelations. Also one can assume without loss of generality that it begins with some element $\, u_x u_y u^{-1}_{xy}\, $ (the other cases are treated similarly). As a relation in $\, U_F\, $ it is elementary, so one can write it as
\smallskip
$$ u_x u_x u^{-1}_{xy} (\cdots ) u_{xy} (\cdots ) u^{-1}_y (\cdots ) u^{-1}_x (\cdots )    $$
\smallskip
where the brackets denote elementary subrelations in $\, U_F\, $. If the first bracket is not to be a word in the generators then either $\, ( \cdots ) u_{xy}\, $ or $\, (\cdots ) u_{xy} u^{\pm }_z\, $ must be, where 
$\, u^{\pm }_z\, $ is the element to the right of $\, u_{xy}\, $. But then, since $\, u_{xy} \notin J_F\, $ the first case is excluded and the second implies $\, u^{\pm }_z = u_{(xy)^{-1}}\, $ which in turn implies that the bracket is a word in the generators, hence empty by the assumption that the relation has no proper subrelations. Similarly, one gets that the second and the third bracket are empty, thus giving a contradiction to the fact that our relation was reduced \qed
\par\bigskip\noindent
Before going on, we generalize a concept implicitely used above that gives $\, U_{N, F}\, $ its universal structure with respect to the inclusion 
$\, N \subseteq F\, $.
\par\bigskip\noindent
{\bf Definition 5. (a)} For any normal pair $\, ( N , F )\, $ define its {\it canonical $F$-central extension} by 
\smallskip
$$ {\overline N}^{ _{\sps F, can}}\, :=\, { U_{N , F} \over 
\bigl[\, J_{N , F}\, ,\, U_F\, \bigr] + \bigl[\, U_{N , F}\, ,\, J_F\,\bigr]}\>  $$
\par\bigskip\noindent
and put
\smallskip
$$  K^J_2 ( N , F )\> :=\> {J_{N , F}\>\cap\> \bigl[\, U_{N , F}\, ,\, U_F\,\bigr] \over
\bigl[\, J_{N , F}\, ,\, U_F\, \bigr] + \bigl[\, U_{N , F}\, ,\, J_F\,\bigr]} \> $$
\par\bigskip\noindent
Clearly $\, K^J_2\, $ is a covariant functor from normal pairs to abelian groups (also called 
relative Schur multipliers).
\par\bigskip\noindent
{\bf (b)}\quad Assume that $\, U\, $ is a free group with normal subgroup $\, K\, $, such that 
$\, G = U / K\, $ is free and given a splitting $\, s : G \rightarrow U\, $ of the quotient map we will say that a set of free generators $\, \{ x_k \}\, $ of a normal subgroup $\, J\, $ of $\, U\, $ has a {\,\it core } for 
$\, s\, $ iff there exists a subset $\, \{ e_i \} \subseteq \{ x_k \}\, $ with 
$\, \{ x_k \} = \{ s(x) e_i s(x)^{-1} \vert x \in G \}\, $. We will say that  $\, J\, $ {\it has a core for s} 
if there exists a set of free generators $\, \{ x_k \}\, $ with this property and that  $\, J\, $ {\it has a core for  s  modulo} $\, \bigl[ J , J \bigr]\, $  if there exists a basis $\, \{ x_k \}\, $ and a subset $\, \{ e_i \}\, $ such that each $\, x_k\, $ is congruent to an element $\, s(x) e_i s(x)^{-1}\, $ modulo $\, [\, J\, , J\, ]\, $. This implies that $\, \{\, \overline { s(x)\, e_i\, s(x)^{-1} }\, \}\, $ freely generates $\, J / [\, J\, , J\, ] \, $ as an abelian group.
\par\medskip\noindent
Eventually, we will also consider some weaker forms of this notion (compare with section 3).
\par\medskip\noindent
{\it Example.} Assume that for a given set of free generators $\, \{ z_i \}\, $ of $\, U\, $ there exists a set of free generators $\, \{ f_j \}\, $ of $\, G\, $ contained in the image of the set $\, \{ z_i \}\, $ and choose a splitting $\, s : G \rightarrow U\, $ mapping this basis to a subset of the $\, \{ z_i \}\, $. Put 
$\, y_j = s( f_j )\, $ and for the remaining $\, z_i\, $ define $\, e_i = z_i s( {\overline z}_i )^{-1}\, $ where 
$\, \overline z_i\, $ is the image of $\, z_i\, $ in $\, G\, $. Then 
$\, \{ e_i \}\, $ is a core of $\, K\, $ for $\, s\, $. If $\, \{ x_k \}\, ,\, \{ y_j \}\, $ as above then it is easy to see that 
$\, \{\, e_i\, ,\, [\, x_k , y^{\pm }_j\, ]\, \}\, $ is a full set of free generators of $\, K\, $ and that 
$\, \{\, \overline{ [\, x_k , y^{\pm }_j\, ]}\, \}\, $ is a set of free abelian generators of $\, [ U , K ] / [ K , K ]\, $ where $\, \overline {[\, x_k , y^{\pm }_j\, ]}\, $ is the image in this group.
\par\bigskip\noindent
{\bf Lemma 3.}\quad Let $\, J_{N , F}\, $ be the kernel of the map $\, U_{N , F} \rightarrow N\, $ and let
$\, s : F / N \nearrow F\, $ be a section to the quotient map $\, F \rightarrow F / N\, $. Denote elements in $\, N\, $ by $\, c , d\, $ etc., elements in $\, s ( F / N)\, $ by letters $\, x, y, z\, $ etc. and define a splitting 
$\, \lambda : J_{F / N} \rightarrow J_F\, $ by the rule 
$\, u_{\xi }\, u_{\eta }\, u^{-1}_{\xi\eta }\,\longmapsto\, 
u_{x}\, u_{y}\, u^{-1}_{x y}\, =\, 
u_{s (\xi )}\, u_{ s (\eta )}\, u^{-1}_{ s ( \xi ) s ( \eta )}\, $. Then $\, J_{N , F}\, $ admits a core for the splitting 
$\,\lambda\, $ which is given by the set $\, \bigl\{\, u_z\, B\, u^{-1}_z\, \bigr\}\, $ with $\, z\, $ in the image of $\, s\, $ and $\, B\, $ is the union of the following sets  
\smallskip
$$ B_1\, =\, \bigl\{\, u_c\, u_d\, u^{-1}_{cd}\, \bigr\}\quad ,\quad B_2\, =\, \bigl\{\, u_x\, u_d\, u^{-1}_{xd}\, \bigr\} \quad , \quad B_3\, =\, \bigl\{\, u_c\, u_{dy}\, u^{-1}_{cdy}\, \bigr\}\, .   $$
\par\smallskip\noindent 
{\it Proof.}\quad  Certainly, a core of $\, J_{N , F}\, $ for the lift $\,\lambda : J_{F / N}\rightarrow J_F\, $ (with respect to the standard bases) induced by $\, s\, $ is given by the set
\smallskip
$$ A\, =\, \bigl\{\, u_{cx}\, u_{dy}\, u^{-1}_{cxdy}\, u_{xy}\, u^{-1}_y\, u^{-1}_x\, \bigr\}     $$
\par\smallskip\noindent
with $\, c, d \in N\, $ and $\, x, y\, $ in the image of $\, s\, $. We claim that taking the subset $\, B\, $ of 
$\, A\, $ as above the set $\, \bigl\{\, u_z\, B\, u^{-1}_z\, \bigr\}\, $ with $\, z\, $ in the image of $\, s\, $ gives an equivalent core of $\, J_{N , F}\, $ for $\,\lambda ( J_{F / N} )\, $. One gets
\smallskip 
$$ u_z\, ( u_c\, u_d\, u^{-1}_{cd} )\, u^{-1}_z\, =\, ( u_z\, u_c\, u^{-1}_{zc} )\, 
( u_{c_z z}\, u_d\, u^{-1}_{c_z zd} )\, ( u_{zcd}\, u^{-1}_{cd}\, u^{-1}_z )         $$
\par\smallskip\noindent
with the first and the third bracket in $\, B^{\pm}_2\,  ,\, c_z = z c z^{-1}\, $, the second bracket in 
$\, A\backslash B\, $. Next we compute 
\smallskip
$$ u_z\, ( u_x\, u_d\, u^{-1}_{xd} )\, u^{-1}_z\, =\, u_z\, u_x\, u^{-1}_{zx}\, ( u_{zx}\, u_d\, u^{-1}_{zxd} )\,
u_{zd_x x}\, u^{-1}_{d_x x}\, u^{-1}_z   .    $$
\par\smallskip\noindent
Put $\, zx = (zx) c_{z,x}\, $ with $\, (zx)\, $ the element of $\, s ( F / N )\, $ corresponding to $\, zx\, $ and 
$\, c_{z,x} \in N\, $ the difference term. Then  $\,( u_{zx}\, u_d\, u^{-1}_{zxd} )\, =\, $
\smallskip
$$ ( u_{(zx) c_{z,x}}\, u^{-1}_{c_{z,x}}\, u^{-1}_{(zx)} )
\bigl[\, u_{(zx)}\, ( u_{c_{z,x}}\, u_d\, u^{-1}_{c_{z,x} d} )\, u^{-1}_{(zx)}\,\bigr]\, 
( u_{(zx)}\, u_{c_{z,x} d}\, u^{-1}_{(zx) c_{z,x} d} ) .   $$
\par\smallskip\noindent 
The middle bracket is of the form already considered and the other two are in $\, B^{\pm}_2\, $, so that modulo a commutator of basis elements in $\, A\, $ already considered with 
$\, u_z\, u_x\, u^{-1}_{zx} \in \lambda ( J_{F / N} )\, $ our term is of the form
\smallskip
$$ ( u_{zx}\,  u_d\, u^{-1}_{zxd} )\, ( u_z\, u_x\, u^{-1}_{zx}\, u_{zd_x x}\, u^{-1}_{d_x x}\, u^{-1}_z ) .   $$
\par\smallskip\noindent 
Changing the core on basis elements with $\, \lambda ( J_{F / N} )\, $ again gives a core for 
$\,\lambda\, $, so that the elements $\,\{ u_z\, B_2\, u^{-1}_z \}\, $ will add the subset
\smallskip 
$$  \bigl\{ u_x\, u_{dy}\, u^{-1}_{xdy}\, u_{xy}\, u^{-1}_y\, u^{-1}_x \bigr\} $$ 
\par\smallskip\noindent
with arbitrary $\, d , x , y\, $ to our list in an independent way. Finally consider 
$\, u_z\, ( u_c\, u_{dy}\, u^{-1}_{cdy} )\, u^{-1}_z\, =\, $
\smallskip
$$ ( u_z\, u_c\, u^{-1}_{zc} )\, ( u_{c_z z}\, u_{dy}\, u^{-1}_{c_z zdy}\, u_{zy}\, u^{-1}_y\, u^{-1}_z )\,
( u_z\, u_y\, u^{-1}_{zy}\, u_{zcdy}\, u^{-1}_{cdy}\, u^{-1}_z ) .  $$
\par\smallskip\noindent
The first bracket is in $\, B_2\, $ and the third is in the range of 
$\,\{ u_z\, B_1\, u^{-1}_z \}\, +\, \{ u_z\, B_2\, u^{-1}_z \}\, $, so that we can add the general elements of the form $\, \{ u_{cx}\, u_{dy}\, u^{-1}_{cxdy}\, u_{xy}\, u^{-1}_y\, u^{-1}_x \}\, $ with $\, c , d , y \neq 1\, $ to our list in an independent way. One checks that this list now exhausts all of $\, A\, $ , and 
$\,\{ u_z\, B\, u^{-1}_z \}\, $ is indeed a core for $\, \lambda ( J_{F/N} )\, $\qed 
\par\bigskip\noindent
{\bf Lemma 4}\quad With notation as in Lemma 3 let $\, \sigma : U_{F / N} \longrightarrow U_F\, $ be the lift induced by the section $\, s\, $. Then a set of generators for the free abelian group 
$\, J_{N , F} / [\, J_{N , F}\, ,\, J_{N , F}\, ]\, $ is given by arbitrary conjugates with elements of 
$\, \sigma ( U_{F / N} )\, $ of elements in the set $\, B_1\cup B_2 \cup \{ u_c\, B_3'\, u^{-1}_c \}\, $ where 
$\, B_3'\, $ is the subset $\,\{\, u_c\, u_y\, u^{-1}_{cy}\,\}\, $ of $\, B_3\, $ and $\, B_2\, $ denotes the set 
$\, \{ [\, u_x\, ,\, u_d\, ]\, u_d\, u^{-1}_{x d x^{-1}}\}\, $ which is congruent to the set denoted $\, B_2\, $ in Lemma 3 modulo the normalization of $\, B_3'\, $.
\par\bigskip\noindent
{\it Proof.}\quad Modulo $\, [\, J_{N , F}\, ,\, J_{N , F}\, ]\, $ and the normalization of $\, B_3'\, $ any conjugate of $\, B_3\, $ is congruent to a corresponding conjugate of $\, B_1\, $ and any element in 
$\, B_2\, $ can be written in the form given by Lemma 4. The normalization of the elements $\, B_3'\, $ is generated modulo $\, [\, J_{N , F}\, ,\, J_{N , F}\, ]\, $ by the set given by the Lemma. Then writing 
\smallskip
$$ \bigl[\, u_x\, ,\, u_c\, u_d\, u^{-1}_{cd}\, \bigr]\> =\> 
\bigl(\,\bigl[\, u_x\, ,\, u_c\, \bigr]\, u_c\, u^{-1}_{x c x^{-1}}\bigr)\> $$
$$ \bigl[\, u_{x c x^{-1}}\, ,\, \bigl(\, \bigl[\, u_x\, ,\, u_d\,\bigr]\, u_d\, u^{-1}_{x d x^{-1}}\,\bigr)\,\bigr]\> 
\bigl(\, \bigl[\, u_x\, ,\, u_d\,\bigr]\, u_d\, u^{-1}_{x d x^{-1}}\,\bigr)\> $$
$$ \bigl(\, u_{x c x^{-1}}\, u_{x d x^{-1}}\, u^{-1}_{x cd x^{-1}}\,\bigr)\>
\bigl(\, u_{x cd x^{-1}}\, u^{-1}_{cd}\,\bigl[\, u_{cd}\, ,\, u_x\,\bigr]\,\bigr)\>
\bigl(\, u_{cd}\, u^{-1}_d\, u^{-1}_c\,\bigr)  $$
\par\bigskip\noindent
one sees that modulo$\, [\, J_{N , F}\, ,\, J_{N , F}\, ]\, $ it is possible to replace the basis given by Lemma 3 by the generating set of Lemma 4. The normalization of 
$B_3'$-elements is obviously independent of the other two sets, and the other two sets are linear independent as elements of $\, J_{N , F}\, $, so they are independent in the abelianization of the corresponding enveloping subgroup of $\, J_{N , F}\, $. However it is not clear a priori that they are independent modulo $\, [\, J_{N , F} , J_{N , F}\, ]\, $ unless one can prove that the corresponding subgroup exhausts all of $\, J_{N , F}\, $  (which is certainly the case when all elements $\, c_{x , y}\, $ are trivial \qed
\par\bigskip\noindent
The following example is very important since it arises naturally in various contexts and since it provides a counterexample to exactness (compare with section 5).
\par\bigskip\noindent
{\bf Lemma 5.}\quad Assume that $\, C\, $ is a free abelian subgroup which is a direct summand of 
$\, B = C\times B / C\, $. Pick a basis 
$\,\{\, f_k\,\} \, $ for $\, C\, $ and an order on the generators $\,\{ f_k \}\, $. Then with notation as above the kernel $\, J_{C , B}\, $ of 
$\, U_{C , B} \rightarrow C\, $ is freely generated by the following list of elements
\smallskip
$$ \bigl\{ \> u^{\pm }_{x_1}\cdots u^{\pm }_{x_k}\,  u_e\, \bigl( u_c\, u_x\, u^{-1}_{cx} \bigr) \, 
u^{-1}_e\, u^{\mp }_{x_k}\cdots u^{\mp }_{x_1}\> \bigr\}\> ,$$
$$ \bigl\{\> u^{\pm }_{x_1}\cdots u^{\pm }_{x_k}\, 
u^{l_1}_{f_{k_1}}\cdots u^{l_s}_{f_{k_s}}\, \bigl[ u_x\, ,\,  u_{f_k} \bigr]\, 
u^{-l_s}_{f_{k_s}}\cdots u^{-l_1}_{f_{k_1}}\, u^{\mp }_{x_k}\cdots u^{\mp }_{x_1}\> \bigr\}\> , $$
$$ \bigl\{\> u^{\pm }_{x_1}\cdots u^{\pm }_{x_k}\, 
u_e\, \bigl( u_d\, u^{- m_r}_{f_{j_r}}\cdots u^{- m_1}_{f_{j_1}} \bigr)\, 
u^{-1}_e\, u^{\mp }_{x_k}\cdots u^{\mp }_{x_1}\> \bigr\}\> ,  $$
$$ \bigl\{\> u^{\pm }_{x_1}\cdots u^{\pm }_{x_k}\, 
u^{n_1}_{f_{i_1}}\cdots u^{n_p}_{f_{i_p}}\, \bigl[\, u_{f_i}\, ,\, u_{f_j}\,\bigr]\, 
u^{-n_p}_{f_{i_p}}\cdots u^{-n_1}_{f_{i_1}}\, u^{\mp }_{x_k}\cdots u^{\mp }_{x_1}\>\bigr\} $$
\par\bigskip\noindent
with 
$\, d = f^{m_1}_{j_1}\cdots f^{m_r}_{j_r}\, ,\, f_{j_1} < \cdots < f_{j_r}\, ,\, f_i < f_j\, ,\, 
f_{i_1} < \cdots < f_{i_p} \leq f_j\, ,\, f_{k_1} <\cdots < f_{k_s} \leq f_k\, ,\, e \in C\, ,\, 
x, x_1 ,\cdots , x_k \in s ( B / C )\, $ arbitrary.
\par\bigskip\noindent
{\it Proof.}\quad We first show that the set generates $\, J_{C , B}\, $. Since $\, C\, $ is a direct summand the section $\, s\, $ can be chosen to be a homomorphism so that $\,\lambda\, $ extends to a map 
$\, U_{B / C} \longrightarrow U_B\, $. Then the elements $\, c_{x , y}\, $ as in the proof of Lemma 3 are all equal to $\, 1\, $ and the conjugates of the set $\, B\, $ with arbitrary elements 
$\, u^{\pm }_{x_1}\cdots u^{\pm }_{x_k}\,\in\, \lambda ( U_{B / C} )\, $ constitute a basis for $\, J_{C , B}\, $.
A basis for $\, U_B\, $ is given by the union of the lifted basis of $\, U_{B / C}\, $, the standard basis of 
$\, U_C\, $ and the elements in the subset $\, B_3' = \{\, u_c\, u_y\, u^{-1}_{cy}\,\}\, $ of $\, B_3\, $. Modulo the normalization of $\, B_3'\, $ the set $\, B_2\,\, $ is congruent to the set $\,\{\, [\, u_x\, ,\, u_c\, ]\,\}\, $ and the set $\, B_3\, $ is congruent to $\, B_1\, $. 
Let $\, d \in C\, $ be an arbitrary element and write $\, d\, $ in the basis of $\, C\, $ as 
$\, f^{m_1}_{j_1}\cdots f^{m_r}_{j_r}\, $ with minimal coefficients $\, m_1,\cdots , m_r\, $ and define the degree of $\, d\, $ to be the sum 
$\, \vert m_1\vert +\cdots + \vert m_r\vert\, $. Let $\, d'\, $ be an element differing from $\, d\, $ by a basis element $\, f^{\pm }_k\,  $ but with smaller degree. Writing 
\smallskip
$$ \bigl[\, u_z\, ,\, u_d\, u_{f^{\pm }_k}\, u^{-1}_{d f_k}\,\bigr]\> =\> $$
$$\> \bigl[\, u_z\, ,\, u_d\,\bigr]\> 
\bigl(\, u_d\, \bigl[\, u_z\, ,\, u_{f^{\pm }_k}\,\bigr]\, u^{-1}_d\,\bigr)\> 
\bigl(\, u_d\, u_{f^{\pm }_k}\, u^{-1}_{d f^{\pm }_k}\,\bigr)\>
\bigl[\, u_{d f^{\pm }_k }\, ,\, u_z\,\bigr]\> \bigl(\, u_{d f^{\pm }_k}\, u^{-1}_{f^{\pm }_k}\, u^{-1}_d\,\bigr)\>  $$
\par\medskip\noindent
and
\smallskip
$$ \bigl[\, u_z\, ,\, u_{f^{-1}_k}\, u_{f_k}\,\bigr]\>  =\> 
\bigl[\, u_z\, ,\, u_{f^{-1}_k}\,\bigr]\> 
\bigl(\, u_{f^{-1}_k}\,\bigl[\, u_z\, ,\, u_{f_k}\,\bigr]\, u^{-1}_{f^{-1}_k}\, \bigr) $$
\par\medskip\noindent 
one finds that one can express (a certain conjugate of) the element $\, [\, u_x\, ,\, u_d\, ]\, $ by 
the same conjugate of $\, [\, u_x\, ,\, u_{d'}\, ]\, $ where 
$\, d'\, $ is of smaller degree modulo elements given by the second set of the Lemma and elements which are conjugates of certain $B_1$-elements by elements from $\,\lambda ( J_{F / N} )\, $. Iterating one can reduce to elements
\smallskip 
$$ u^{l_1}_{f_{k_1}}\cdots u^{l_s}_{f_{k_s}}\, [\, u_x\, ,\, u_{f_k}\, ]\, 
u^{-l_s}_{f_{k_s}}\cdots u^{-l_1}_{f_{k_1}} $$ 
\par\medskip\noindent
subject to the condition 
$\, f_{k_1} <\cdots < f_{k_s}\, $ modulo certain conjugates of elements of type $\, B_1\, $. It is easy to see  that the central parts (contained in $\, J_C\, $) of the third and the fourth set generate $\, J_C\, $. Performing a change of basis of 
$\, U_C\, $ by replacing 
$\, u_d\, $ with the corresponding central bracket of the third set one finds that $\, J_C\, $ is generated by conjugates of such elements and of elementary commutators $\, [ u_{f_i}\, ,\, u_{f_j}\, ]\, $ by arbitrary ordered coefficients $\, u^{l_1}_{f_{j_1}}\cdots u^{l_p}_{f_{j_p}}\, $ which again can be replaced by the corresponding elements $\, u_e\, $ in the third set. As for the fourth one notes on writing
\smallskip
$$ u^{\pm }_{f_k}\, \bigl[\, u_{f_i}\, ,\, u_{f_j}\,\bigr]\, u^{\mp }_{f_k}\> =\> 
\bigl[\, u^{\pm }_{f_k}\, ,\, u_{f_i}\,\bigr]\>\bigl(\, u_{f_i}\, \bigl[\, u^{\pm }_{f_k}\, ,\, u_{f_j}\,\bigr]\, u^{-1}_{f_i}\,\bigr)\> \bigl[\, u_{f_i}\, ,\, u_{f_j}\,\bigr] $$
$$ \qquad\qquad\quad \bigl(\, u_{f_j}\,\bigl[\, u_{f_i}\, ,\, u^{\pm }_{f_k}\,\bigr]\, u^{-1}_{f_j}\,\bigr)\> \bigl[\, u_{f_j}\, ,\, u^{\pm }_{f_k}\,\bigr]\> $$
\par\bigskip\noindent
that one can restrict to conjugates with ordered elements of the form as given above. In order to that one can restrict to conjugacy coefficients satisfying $\, f_{k_s} \leq f_k\, $ in the second set one notes the relations
\smallskip
$$ 0\quad =\quad \bigl[\, u_x\, ,\, u_{f_k}\,\bigr]\>
\bigl(\, u_{f_k}\,\bigl[\, u_x\, ,\, u_{f_l}\,\bigr]\, u^{-1}_{f_k}\,\bigr)\> 
\bigl[\, u_{f_k}\, ,\, u_{f_l}\,\bigr]\qquad\quad $$
$$ \quad\qquad\bigl(\, u_{f_l}\,\bigl[\, u_{f_k}\, ,\, u_x\,\bigr]\, u^{-1}_{f_l}\,\bigr)\>
\bigl[\, u_{f_l}\, ,\, u_x\,\bigr]\> \bigl(\, u_x\,\bigl[\, u_{f_l}\, ,\, u_{f_k}\,\bigr]\, u^{-1}_x\,\bigr)\> . $$
\par\medskip\noindent
Next we need to show that the generating set of Lemma 5 is independent.
Let some relation in the generators be given. We may assume the relation to be cyclically reduced and minimal (i.e. to contain no proper subrelation). Then from the basis of $\, U_B\, $ given by the union 
$\,\{\, u_c\,\}\, \cup\,\{\, u_z\,\}\, \cup\, B_3'\, $ one induces that the relation cannot contain any elements of the first set since the central bracket of these elements is in $\, B_3'\, $ and is the only appearance of such elements in all four sets so it can only be cancelled by the central bracket of another element of the first set which then must have the same coefficients $\, u^{\pm }_{x_1}\cdots u^{\pm }_{x_k}\, u_e\, $ since otherwise the elements inbetween would map to a nontrivial element in $\, J_{B / C}\, $ respectively $\, C\, $ which is impossible so the element whose central bracket cancels against the central bracket of the original element is just its inverse. Since the relation is minimal it consists only of elements of the last three types. By a similar argument the third set is independent from the second and fourth, since it is the only one containing elements $\, u_d\, , d \neq f_k\, $ and it is easy to see that the elements of the third set are mutually independent, so we have reduced to the second and fourth set.  
These two sets are easily identified with a subbasis of $\, [\, U_{C , B}\, ,\, U_B\, ]\, $. We leave the details to the reader  \qed 
\par\bigskip\noindent
{\bf Lemma 6.}\quad With assumptions as in Lemma 5 the extension
\smallskip
$$ 1\,\largerightarrow\, J_{J_{C , B}\, ,\, U_B} \,\largerightarrow\, U_{J_{C , B}\, ,\, U_B}\,
\largerightarrow\, J_{C , B}\,\largerightarrow\, 1  $$
\par\smallskip\noindent
admits a splitting which is normal for the canonical lift of $U_{B / C}$ (assuming $\, B / C \subset B\, $ by some corresponding decomposition) and the canonical lift is defined by   
\smallskip
$$ U_B\,\longrightarrow U^2_B\quad ,\quad u_x\,\longmapsto\, u_{u_x}\>  , $$
and, if the rank of $\, C\, $ is one, the lift is normal for the canonical lift of $\, U_B\, $ modulo
\smallskip
$$ \bigl[\, U_{J_{C , B}\, ,\, U_B}\, ,\, J_{U_B}\,\bigr]\>  .$$
\par\bigskip\noindent
In particular $\, K^J_2 ( J_{C , B}\, ,\, U_B )\, $ is equal to 0 if the rank of $\, C\, $ is one and in general is a free abelian group.
\par\bigskip\noindent
{\it Proof.}\quad 
Modulo the subgroup $\, [\, J_{J_{C , B}\, ,\, U_B}\, ,\, U_{U_{C , B}\, ,\, U_B}\, ]\, $ any lift of 
$\, J_{C , B}\, $ to $\, U_{J_{C , B}\, ,\, U_B}\, $ will contain a canonical copy of 
$\, [\, J_{C , B}\, ,\, J_{C , B}\, ]\, $ which does not depend on the particular lift since 
$\, J_{J_{C , B}\, ,\, U_B}\,\cap\, [\, U_{J_{C , B}\, ,\,, U_B}\, ,\, U_{J_{C , B}\, ,\, U_B}\, ]\, =\, 
[\, J_{J_{C , B}\, ,\, U_B}\, ,\, U_{J_{C , B}\, ,\, U_B}\, ]\, $ (the identity follows because $\, J_{C , B}\, $ is a free group) is divided out. 
We may then further divide by this canonical copy which amounts to dividing by 
$\, [\, U_{J_{C , B}\, ,\, U_B}\, ,\, U_{J_{C , B}\, ,\, U_B}\, ]\, $ and construct a $U_B$-normal lift of the free abelian group $\, \overline J_{C , B}\> =\> J_{C , B} / [\, J_{C , B}\, ,\, J_{C , B}\, ]\, $ into the image of 
$\, U_{J_{C , B}\, ,\, U_B}\, $. Define a lift $\, \overline J_{C , B}\longrightarrow \overline U_{J_{C , B}\, ,\, U_B}\, $ by 
\bigskip
$$ u^{\pm }_{x_1}\cdots u^{\pm }_{x_n}\, u_e\, 
\bigl( u_c\, u_x\, u^{-1}_{cx} \bigr)\, 
u^{-1}_e\, u^{\mp }_{x_n}\cdots u^{\mp }_{x_1}\quad \longmapsto  $$ 
$$ u^{\pm }_{u_{x_1}}\cdots u^{\pm }_{u_{x_n}}\, u_{u_e}\>
u_{( u_c\, u_x\, u^{-1}_{cx} )}\>
u^{-1}_{u_e}\, u^{\mp }_{u_{x_n}}\cdots u^{\mp }_{u_{x_1}}\>  , $$
\par\medskip\noindent
$$ u^{\pm }_{x_1}\cdots u^{\pm }_{x_m}\, u^{l_1}_{f_{k_1}}\cdots u^{l_s}_{f_{k_s}}\,
\bigl[ u_x\, ,\, u_{f_k} \bigr]\, u^{-k_s}_{f_{l_s}}\cdots u^{-l_1}_{f_{k_1}}
\, u^{\mp }_{x_m}\cdots u^{\mp }_{x_1}\quad\longmapsto $$
$$ u^{\pm }_{u_{x_1}}\cdots u^{\pm }_{u_{x_m}}\,  u^{l_1}_{u_{f_{k_1}}}\cdots u^{l_s}_{u_{f_{k_s}}}\>
u_{ [\, u_x\, ,\,  u_{f_k}\, ]}\> u^{-l_s}_{u_{f_{k_s}}}\cdots u^{-l_1}_{u_{f_{k_1}}}\, 
u^{\mp }_{u_{x_m}}\cdots u^{\mp }_{u_{x_1}}\>  ,  $$
\par\medskip\noindent
$$ u^{\pm }_{x_1}\cdots u^{\pm }_{x_k}\, u_e\, 
\bigl( u_c\, u^{- m_r}_{f_{j_r}}\cdots u^{- m_1}_{f_{j_1}} \bigr)
\, u^{-1}_e\, u^{\mp }_{x_k}\cdots u^{\mp }_{x_1}\quad \longmapsto $$
$$ u^{\pm }_{u_{x_1}}\cdots u^{\pm }_{u_{x_k}}\, u_{u_e}\>
u_{ (\, u_c\, u^{- m_r}_{f_{j_r}}\cdots u^{- m_1}_{f_{j_1}} )}\>
u^{-1}_{u_e}\, u^{\mp }_{u_{x_k}}\cdots u^{\mp }_{u_{x_1}}\>  ,  $$
\par\medskip\noindent
$$  u^{\pm }_{x_1}\cdots u^{\pm }_{x_k}\, u^{n_1}_{f_{i_1}}\cdots u^{n_p}_{f_{i_p}}\,  
\bigl[\, u_{f_i}\, ,\, u_{f_j}\,\bigr]\, u^{-n_p}_{f_{i_p}}\cdots u^{-n_1}_{f_{i_1}}\, 
u^{\mp }_{x_k}\cdots u^{\mp }_{x_1}\quad\longmapsto $$
$$  u^{\pm }_{u_{x_1}}\cdots u^{\pm }_{u_{x_k}}\, 
u^{n_1}_{u_{f_{i_1}}}\cdots u^{n_p}_{u_{f_{i_p}}}\, u_{[\, u_{f_i}\, ,\, u_{f_j}\, ]}\, 
u^{-n_p}_{u_{f_{i_p}}}\cdots u^{-n_1}_{u_{f_{i_1}}}\, 
u^{\mp }_{u_{x_k}}\cdots u^{\mp }_{u_{x_1}}\> . $$
\par\bigskip\noindent
One easily checks that the map is compatible with the canonical lift $\, U_B \rightarrow U^2_B\, $ modulo 
$\, \left[\, J_{ U_{C , B}\, ,\, U_B}\, ,\, U_{ J_{C , B}\, ,\, U_B}\,\right]\, $ except possibly for elements of the second and fourth type where it is not a priori clear that the lift is normal for the action of elements 
$\, u_{u_{f_l}}\, ,\, f_l > f_k , f_j\, $. However this is only a problem when the rank of $\, C\, $ is larger than one so that in any case $\, K^J_2 ( J_{C , B}\, ,\, U_B ) = 0\, $ if the rank is one. 
Now change the last set of the basis to the form
\smallskip
$$ \bigl\{\> u^{n_1}_{f_{i_1}}\cdots u^{n_p}_{f_{i_p}}\, 
u^{\pm }_{x_1}\cdots u^{\pm }_{x_k}\, \bigl[\, u^{\pm }_x\, ,\,\bigl[\, u_{f_i}\, ,\, u_{f_j}\,\bigr]\,\bigr]\, 
u^{\mp }_{x_k}\cdots u^{\mp }_{x_1}\, u^{-n_p}_{f_{i_p}}\cdots u^{-n_1}_{f_{i_1}}\>\bigr\} $$
$$ \cup\quad \bigl\{\, u^{n_1}_{f_{i_1}}\cdots u^{n_p}_{f_{i_p}}\,
\bigl[\, u^{\pm }_{f_l}\, ,\,\bigl[\, u_{f_i}\, ,\, u_{f_j}\,\bigr]\, \bigr]\, u^{-i_p}_{f_{n_p}}\cdots u^{-n_1}_{f_{i_1}} 
\>\bigr\}\quad\cup\quad \{\, \bigl[\, u_{f_i}\, ,\, u_{f_j}\,\bigr]\,\} \>  $$
\par\medskip\noindent
and the second set to
\smallskip
$$ \{\> u^{l_1}_{f_{k_1}}\cdots u^{l_s}_{f_{k_s}}\, u^{\pm }_{x_1}\cdots u^{\pm }_{x_m}\, 
\bigl[\, u^{\pm }_z\, ,\,\bigl[ u_x\, ,\, u_{f_k} \bigr]\,\bigr]\,  u^{\mp }_{x_m}\cdots u^{\mp }_{x_1}\, 
u^{-l_s}_{f_{k_s}}\cdots u^{-l_1}_{f_{k_1}}\>\} $$
$$ \cup\quad \{\> u^{l_1}_{f_{k_1}}\cdots u^{l_s}_{f_{k_s}}\,
\bigl[\, u^{\pm }_{f_l}\, ,\,\bigl[ u_x\, ,\, u_{f_k} \bigr]\,\bigr]\, u^{-l_s}_{f_{k_s}}\cdots u^{-l_1}_{f_{k_1}}\>\}
\quad\cup\quad \{\, \bigr[\, u_x\, ,\, u_{f_k}\,\bigr]\,\}\> . $$
\par\medskip\noindent
Then it is plain to see that the normal splitting extends to the first and second parts modulo the subgroups
$\, [\, U_{J_{C , B}\, ,\, U_B}\, ,\, J_{U_B}\, ]\, +\, [\, J_{J_{C , B}\, ,\, U_B}\, ,\, U^2_B\, ]\, $ 
and that modulo these subgroups all sets together lift normally for the action of $\, U_{B / C}\, $. The last set of the first part is normal modulo the image of $\, J_C \cap [\, U_C\, ,\, U_C\, ]\, $ which represents the group $\, K^J_2 ( J_C , U_C )\, $ (the map 
$\, K^J_2 ( J_C , U_C ) \rightarrowtail K^J_2 ( J_{C , B} , U_B )\, $ is injective by existence of a splitting) and it is easily checked that each relation obtained from the Jacobi identity and expressing 
an element $\, u_{f_l}\, [\, u_{f_i}\, ,\, u_{f_j}\, ]\, u^{-1}_{f_l}\, $ with $\, f_l > f_j\, $ in terms of the basis of Lemma 5 contributes a copy of $\,\mathbb Z\, $ to $\, K^J_2 ( J_C , U_C )\, $ and that they are all independent, so this subgroup is free abelian. This can be rigorously proved using an induction argument on the number of generators of $\, C\, $ noting that $\, ( J_C , U_C )\, $ is the inductive limit of the subpairs $\, ( J_{C_{\lambda }} , U_{C_{\lambda }} )\, $ for any subnet $\, c_{\lambda }\, $ converging to $\, C\, $ and it is easy to see that the functor $\, K^J_2\, $ is continuous. Similarly, each relation obtained from the Jacobi identity and expressing an element 
$\, u_{f_l}\, [\, u_x\, ,\, u_{f_k}\, ]\, u^{-1}_{f_l}\, $ in terms of the basis of Lemma 5 contributes a copy of 
$\,\mathbb Z\, $ to the quotient of $\, K^J_2 ( J_{C , B}\, ,\, U_B )\, $ by the image of 
$\, K^J_2 ( J_C , U_C )\, $, and again these are all independent of each other, noting that the group structure of $\, B / C\, $ is in no way involved, only its cardinality, so by a similar induction argument one gets the result. In conclusion one obtains a description of $\, K^J_2 ( J_{C , B}\, ,\, U_B )\, $ as a free abelian group by a set of free generators   \qed  
\par\bigskip\noindent
{\bf Lemma 7.}\quad Given a commutative diagram with exact rows and columns
\smallskip
$$ \vbox{\halign{ #&#&#&#&#\cr
$J_{N , F}$ & $\largerightarrow$ & $U_{N , F}$ & $\largerightarrow$ & $N$ \cr
\hfil $\Biggm\downarrow$\hfil && \hfil $\Biggm\downarrow$\hfil &&  $\Biggm\downarrow$ \cr
\hfil $J_F$\hfil & $\largerightarrow$ & \hfil $U_F$\hfil & $\largerightarrow$ & $F$ \cr
\hfil $\Biggm\downarrow$\hfil && \hfil $\Biggm\downarrow$\hfil &&  $\Biggm\downarrow$ \cr
\hfil $J_B$\hfil & $\largerightarrow$ & \hfil $U_B$\hfil & $\largerightarrow$ & $B$ \cr }}   $$
\par\bigskip\noindent
assume there is a map $\, U_B\buildrel\beta\over\longrightarrow F\, $ over the identity of $\, B\, $ inducing $\, J_B \longrightarrow N\, $. Then there is a natural lift 
$\, J_F\buildrel\alpha\over\longrightarrow U_{N , F}\, $ of $\, J_F\rightarrow J_B\rightarrow N\, $ such that 
$\,\alpha \, $ is the inverse on $\, J_{N , F}\, $ modulo 
$\, [\, J_{N , F}\, ,\, J_F\, ]\, +\, [\, J_{N , F}\, ,\, U_{N , F}\, ]\, $. 
\par\bigskip\noindent
{\it Proof.}\quad By Lemma 2 $\, J_F\, $ is free on the generators
\smallskip 
$$ \{\> u_x\, u_y\, u^{-1}_{xy}\> ,\> u_x\, u_{x^{-1}}\> \vert\> x , y \in F\>\}\>   . $$
\par\smallskip\noindent
Define $\,\alpha\, $ by the formulas 
\smallskip
$$ u_x\, u_y\, u^{-1}_{xy}\,\longmapsto\, u_{xy}\, u^{-1}_y\, u^{-1}_x\, 
u_{\beta ( u_{\overline x} )}\, u_{\beta ( u_{\overline y} )}\, u^{-1}_{\beta (u_{\overline {xy}} )}\>  ,  $$
$$ u_x\, u_{x^{-1}}\,\longmapsto\, 
u^{-1}_{x^{-1}}\, u^{-1}_x\, u_{\beta ( u_{\overline x} )}\, u_{\beta ( u_{\overline {x^{-1}}} )}  $$
\par\smallskip\noindent
where $\, \overline x\, $ denotes the image of $\, x\, $ in $\, B\, $. One checks that $\,\alpha\, $ is a lift of 
$\, J_F \rightarrow J_B \rightarrow N\, $ and that it induces the inverse map on the set of normal generators as given by Lemma 3 for $\, J_{N , F}\, $. Hence dividing by 
$\, [\, J_{N , F}\, ,\, J_F\, +\, U_{N , F}\, ]\, $ gives the result \qed
\par\bigskip\noindent
\bigskip\bigskip\bigskip\bigskip\bigskip\bigskip\bigskip
\par\noindent
\hfil \Large{\bf{Section II -- Regular suspension}} \hfil 
\par\bigskip\bigskip\noindent

Before coming to Theorem 2, which is the very heart of the theory, we need to slightly generalize the concept of a universal $F$-central extension. We are not yet really able to prove Theorem 2 since important preliminaries have to wait for later chapters, but we give an extensive sketch of its proof modulo the cited results especially Theorem 5 which is basic for any more sophisticated argument.  
\par\bigskip\noindent
{\bf Definition 6.}  Let $\, ( N , F ) \subseteq ( M , G )\, $ be an inclusion of normal pairs. $\, ( N , F )\, $ is said to admit a 
{\it  semiuniversal $F$-central extension relative to} $\, ( M , G )\,$ if there exists an $F$-central extension $\, {\overline N}^{ _{\sps F}}\, $ (which then embeds into the canonical $G$-central extension of $\, M\, $) such that for any other  $F$-central extension $\, \widetilde N\, $ of  
$\, N\, $ extending to a $G$-central extension $\, \widetilde M\, $ of 
$\ M\, $ there is a map $\, {\overline N}^{ _{\sps F}} \longrightarrow {\widetilde N}\, $ over $\, ( N , F )\, $, uniquely determined on the central kernel which is contained in $\, [ {\overline N}^{ _{\sps F}} , F ]\, $. In case that the kernel of $\, {\overline N}^{ _{\sps F}}\, $ is equal to $\, K^J_2 ( N , F )\, $ (equivalent to injectivity of $\, K^J_2 ( N , F ) \rightarrowtail K^J_2 ( M , G )\, $) the relative semiuniversal $F$-central extension will be called {\it faithful},  and 
{\it universal}, if it is faithful and every $F$-central extension of 
$\, {\overline N}^{ _{\sps F}}\, $ which extends to a double $G$-central extension of $\, M\, $ splits.
\par\bigskip\noindent
For simplicity the notation in this section differs from the usual notation in other chapters: the subgroups $\, U^n_{N , F}\, $ and $\, J^n_{N , F}\, $ are defined inductively with respect to the regular suspension 
$\, R_{N , F}\, $ as defined in section 1, i.e. $\, U^{n+1}_{N , F}\, =\, U_{R^n_{N , F}\, ,\, U^n_F}\, $ and 
$\, J^{n+1}_{N , F}\, =\, J_{R^n_{N , F}\, ,\, U^n_F}\, $ whereas the usual definition is by induction replacing $\, R^n_{N , F}\, $ with $\, J^n_{N , F}\, $ ($\, n\geq 0\, $). The importance of Theorem 2 is that it provides the regular theory with a sufficiently large and interesting class of examples. Once this is established the theory may have some extensions from the category of perfect groups and their normal subgroups to some larger class in one way or the other (for example the category of almost perfect groups, abstract regular pairs etc. ). This will be investigated more thoroughly in subsequent chapters.
\par\bigskip\noindent
{\bf Theorem 2.}\quad  Let $\, N\, $ be a full subgroup of a perfect group $\, F\, $ and 
$\, ( R^n_{N , F}\, ,\, U^n_F )\, $ an n-th regular suspension of $\, ( N , F )\, $. The pair 
$\, ( R^n_{N , F}\, , U^n_F )\, $ admits a relative universal $U^n_F$-central extension 
$\, {\overline {R^n_{N , F}}}^{ _{\sps U^n_F}}\, $ for the inclusion into $\, ( J^n_{N , F}\, ,\, U^n_F )\, $ for every $\, n\geq 1\, $. The kernel of this extension is contained in the subgroup  
$\, [\, [\cdots [\, {\overline {R^n_{N , F}}}^{ _{\sps U^n_F}}\, ,\, U^n_F\, ] ,\cdots ] , U^n_F\, ]\, $ of k-fold commutators modulo the canonical $U^n_F$-normal copy of 
$\, [\, J^n_{N , F}\, ,\, J^n_{N , F}\, ]\, \subseteq\, {\overline R}^{ _{\sps U^n_F}}\, $ for every $\, k \geq 1\, $.
\par\bigskip\noindent
{\it Proof.}\quad It will follow from the results in section 5 and Lemma 15 of section 6 that the natural map 
$\, {\overline {R^n_{N , F}}}^{ _{\sps U^n_F, can}} \longrightarrow 
{\overline {J^n_{N , F}}}^{ _{\sps U^n_F, can}}\, $ is injective so this assumption of Definition 5 (iii) is satisfied. The proof is by induction on $\, n\, $. We assume that we have already constructed the relative universal $U^k_F$-central extension of $\, R^k_{N , F}\, $ for all 
$\, k < n\, $. Theorem 1 gives the case $\, n = 0\, $.  One remarks the following facts : since 
$\, U^{n+1}_F / J^{n+1}_F\, =\, U^n_F\, $ and 
$\, U^{n+1}_F / U^n_{N , F}\, =\, U_{U^n_F / R^n_{N , F}}\, $ are free groups one gets
\smallskip
$$ J^{n+1}_{N , F}\,\cap\, \left[\, J^{n+1}_F\, ,\, U^{n+1}_F\,\right]\> =\> 
J^{n+1}_{N , F}\,\cap\, \left[\, U^{n+1}_{N , F}\, ,\, U^{n+1}_F\, \right]   $$
\par\smallskip\noindent
and more generally
\smallskip
$$ J^{n+1}_{N , F}\,\cap\, 
\left[\, \left[\cdots \left[\, J^{n+1}_F\, ,\, U^{n+1}_F\,\right] ,\cdots\right] , U^{n+1}_F\,\right]\> =  $$ 
$$ J^{n+1}_{N , F}\,\cap\, 
\left[\, \left[\cdots \left[\, U^{n+1}_{N , F}\, ,\, U^{n+1}_F\,\right] ,\cdots \right] , U^{n+1}_F \right]\>  .  $$
\par\smallskip\noindent
for arbitrary k-fold commutators. From the canonical lift 
$\, U^n_F \rightarrow U^{n+1}_F\, $, $\, u_x\, \mapsto\, u_{u_x}\, $, for $\, x \in U^{n-1}_F\, $ one gets a decomposition
\smallskip
$$ U^{n+1}_F\> =\> J^{n+1}_F\,\rtimes\, U^n_F\>  .  $$
\par\smallskip\noindent
Consider the following commutative diagram 
\bigskip
$$ \vbox{\halign{ #&#&#&#&#\cr
\hfil $J_{J_{J^{n-1}_{N , F}\, ,\, U^{n-1}_F}\, ,\, U^n_F}$ \hfil & $\largerightarrow$ & 
\hfil $U_{J_{J^{n-1}_{N , F}\, ,\, U^{n-1}_F}\, ,\, U^n_F}$\hfil & $\largerightarrow$ & 
$J_{J^{n-1}_{N , F}\, ,\, U^{n-1}_F}$ \cr
\hfil $\Biggm\downarrow$\hfil && \hfil $\Biggm\downarrow$\hfil && $\,\,\Biggm\downarrow$ \cr
(2) $\quad J^{n+1}_F$\hfil & \hfil $\largerightarrow$\hfil & \hfil $U^{n+1}_F$\hfil & 
\hfil $\largerightarrow$\hfil &
$U_{J^{n-1}_{N , F}\, ,\, U^{n-1}_F} \rtimes U_{U^{n-1}_F\, /\, J^{n-1}_{N , F}}$\hfil \cr
\hfil $\Biggm\downarrow$\hfil && \hfil $\Biggm\downarrow$\hfil && $\,\,\Biggm\downarrow$ \cr
$J_{U^n_F\, /\, J_{J^{n-1}_{N , F}\, ,\, U^{n-1}_F}}$ & \hfil $\largerightarrow$\hfil & 
\hfil $U_{U^n_F\, /\, J_{J^{n-1}_{N , F}\, ,\, U^{n-1}_F}}$\hfil & \hfil $\largerightarrow$\hfil &
$J^{n-1}_{N , F}\,\rtimes\, U_{U^{n-1}_F\, /\, J^{n-1}_{N , F}}$\hfil \cr}}  $$
\par\bigskip\noindent
The group 
$\, U^n_F / J_{J^{n-1}_{N , F}\, ,\, U^{n-1}_F}\, $ contains the two mutually commuting (complementary) normal subgroups 
$\, J^{n-1}_{N , F}\, $ and $\, P^{n-1}_{N , F}\, =\, J^n_F / J_{J^{n-1}_{N , F}\, ,\, U^{n-1}_F}\, $. Consider the commutative diagram
\smallskip
$$ \vbox{\halign{ #&#&#\cr
$U^n_F / J_{J^{n-1}_{N , F}\, ,\, U^{n-2}_F}\, $ & $\buildrel p_1\over\longrightarrow$ & $U^{n-1}_F$ \cr
$\> { _{p_2}} \Biggm\downarrow$\hfil &&  $\,\, \Biggm\downarrow\>$ \cr
$U^n_F / U_{J^{n-1}_{N , F}\, ,\, U^{n-1}_F}$ & $\longrightarrow$ & 
$R^{n-2}_{N , F}\,\rtimes\, U_{U^{n-2}_F / R^{n-2}_{N , F}}$\cr }}  $$
\par\smallskip\noindent
(with $\, R^{n-2}_{N , F}\,\rtimes\, U_{U^{n-2}_F /R^{n-2}_{N , F}}\> =\> F / N\, $ if $\, n = 1\, $). Fix a section 
\smallskip
$$ {\nu }_1 : R^{n-2}_{N , F}\,\rtimes\, U_{U^{n-2}_F /R^{n-2}_{N , F}}\,\nearrow\, U^{n-1}_F $$
\par\smallskip\noindent 
and a section
\smallskip 
$$ {\nu }_2 : R^{n-2}_{N , F}\,\rtimes\, U_{U^{n-2}_F /R^{n-2}_{N , F}}\,\nearrow\, 
U^n_F / U_{J^{n-1}_{N , F}\, ,\, U^{n-1}_F}\>  . $$
\par\smallskip\noindent
Then there is a section 
$\, {\mu }_1 : U^{n-1}_F\,\nearrow\, U^n_F / J_{J^{n-1}_{N , F}\, ,\, U^{n-1}_F}\, $ such that 
\smallskip
$$ {\mu }_1 ( x\, {\nu }_1 ( y ))\, =\, x\, ({\mu }_1\circ {\nu }_1 ) ( y ) $$ 
\par\smallskip\noindent
for $\, x \in J^{n-1}_{N , F}\, $ and 
$\, {\nu }_1 ( y )\, $ in the image of $\, {\nu }_1\, $. Similarly, there is a section 
$\, {\mu }_2 : U^n_F / U_{J^{n-1}_{N , F}\, ,\, U^{n-1}_F}\,\nearrow\, 
U^n_F / J_{J^{n-1}_{N , F}\, ,\, U^{n-1}_F}\, $ such that 
$$ {\mu }_2 ( x\, {\nu }_2 ( y ))\, =\, x\, ({\mu }_2\circ {\nu }_2 ) ( y )\, =\, x\, ({\mu }_1\circ{\nu }_1 ) ( y ) $$ 
\par\smallskip\noindent
for $\, x \in P^{n-1}_{N , F}\, $ and 
$\, y \in R^{n-2}_{N , F}\,\rtimes\, U_{U^{n-2}_F / R^{n-2}_{N , F}}\, $. The section 
$\, {\mu }_1\, $ induces a splitting 
\smallskip
$$\, {\sigma }_1\, :\,  U^n_F\, \longrightarrow\, 
U_{U^n_F / J_{J^{n-1}_{N , F}\, ,\, U^{n-1}_F}} $$ 
\par\bigskip\noindent
and $\, {\mu }_2\, $ induces a splitting
\smallskip 
$$\, {\sigma }_2\, :\, U_{U^n_F / U_{J^{n-1}_{N , F}\, ,\, U^{n-1}_F}}\, \longrightarrow\, 
U_{U^n_F / J_{J^{n-1}_{N , F}\, ,\, U^{n-1}_F}}\>  . $$ 
\par\bigskip\noindent
Define a section $\, s : J^{n-1}_{N , F}\, \rtimes\, U_{U^{n-1}_F\, /\, J^{n-1}_{N , F}}\, \nearrow\, U^n_F\, $
by $\, s (xy) = s (x)\, s (y)\, $ with $\, s\, $ a semicanonical lift for  $\, x \in J^{n-1}_{N , F}\, $ and $\, s\, $ equal to the lift corresponding to the decomposition in diagram (2) for 
$\, y \in U_{U^{n-1}_F\, /\, J^{n-1}_{N , F}}\, $. Let $\,\sigma\, $ be the corresponding lift 
$\, U_{U^n_F\, /\, J_{J^{n-1}_{N , F}\, ,\, U^{n-1}_F}}\,\longrightarrow\, U^{n+1}_F\, $. 
Consider the subgroups 
$\, U_{P^{n-1}_{N , F}\, ,\, U^n_F / J_{J^{n-1}_{N , F}\, ,\, U^{n-1}_F}}\, $ and  
$\, s_1 ( U_{J^{n-1}_{N , F}\, ,\, U^{n-1}_F} )\, $ of $\, U_{U^n_F / J_{J^{n-1}_{N , F}\, ,\, U^{n-1}_F}}\, $.  
The first group has a core for $\, s_1 ( U^n_F )\, $ and the second is a subgroup of 
$\, U_{J^{n-1}_{N , F}\, ,\, U^n_F / J_{J^{n-1}_{N , F}\, ,\, U^{n-1}_F}}\, $ by our choice of $\, {\mu }_1\, $ and has a core for 
\smallskip
$$ \langle \left\{ u_{{\mu }_1\circ {\nu }_1 ( y )}\right\} \rangle\, =\, 
\langle \left\{ u_{{\mu }_2\circ {\nu }_2 ( y )} \right\} \rangle \>  . $$
The commutator subgroup 
$\,C = \bigl[\, U_{P^{n-1}_{N , F}\, ,\, U^n_F / J_{J^{n-1}_{N , F}\, ,\, U^{n-1}_F}}\, ,\, 
s_1 ( U_{J^{n-1}_{N , F}\, ,\, U^{n-1}_F} )\,\bigr]\, $
is normal in $\, U_{U^n_F / J_{J^{n-1}_{N , F}\, ,\, U^{n-1}_F}}\, $ and lies in the kernel of 
\smallskip
$$ U_{J^{n-1}_{N , F}\, ,\,U^n_F / J_{J^{n-1}_{N , F}\, ,\, U^{n-1}_F}} \longrightarrow J^{n-1}_{N , F}
\>  . $$
\par\smallskip\noindent
If $\, C_P\, $ denotes the core of 
$\, U_{P^{n-1}_{N , F}\, ,\, U^n_F / J_{J^{n-1}_{N , F}\, ,\, U^{n-1}_F}}\, $ and $\, C_R\, $ the core of 
$\, s_1 ( U_{J^{n-1}_{N , F}\, ,\, U^{n-1}_F} )\, $ then a basis for the commutator subgroup is given by the set of elements
\smallskip
$$ \left\{\, u^{\pm }_{j_1}\cdots u^{\pm }_{j_n}\, v^{\pm }_{l_1}\cdots v^{\pm }_{l_m}\, 
\left[\, u^{\pm }_{j_0}\, ,\, v_{l_0}\,\right]\,
v^{\mp }_{l_m}\cdots v^{\mp }_{l_1}\, u^{\mp }_{j_n}\cdots u^{\mp }_{j_1}\,\vert\,  m\geq 1\, \right\}\> , $$
$$ \cup \left\{\, u^{\pm }_{j_1}\cdots u^{\pm }_{j_n}\,\left[\, u_{j_0}\, ,\, v^{\pm }_{l_0}\,\right]\,
u^{\mp }_{j_n}\cdots u^{\mp }_{j_1}\,\vert\, m = 0\, \right\}   $$
\par\smallskip\noindent
where the $\, \{\, u_{j_i}\,\}\, $ are taken from the set $\, \{ s_1 ( y )\, C_P\, s_1 ( y )^{-1}\, \}\, $ with $\, y\, $ in the image of the splitting 
$\, U_{{\overline R}^{n-2}_{N , F}\, ,\, U^{n-1}_F / J^{n-1}_{N , F}}\,\rtimes\, 
U^2_{U^{n-2}_F / R^{n-2}_{N , F}}\,\longrightarrow U^n_F\, $ induced by $\, {\nu }_1\, $ and the 
$\, \{\, v_{l_k}\,\}\, $ are taken from the basis of 
$\, s_1 ( U_{J^{n-1}_{N , F}\, ,\, U^{n-1}_F} )\, $. Note that the basis has a core for $\, \{\, s_1 ( y )\,\}\, $, for $\, C_P\, $ and also for $\, C_R\, $ modulo $\, [\, C\, ,\, C\, ]\, $
which together form a basis for $\, U_{U^n_F / J_{J^{n-1}_{N , F}\, ,\, U^{n-1}_F}}\, $. Let 
$\,\widetilde C\, $ denote the preimage of $\, C\, $ in $\, U^{n+1}_F\, $ and $\, \widetilde U\, $ the preimage of 
$\, J_{U^n_F / J_{J^{n-1}_{N , F}\, ,\, U^{n-1}_F}}\, +\, {\sigma }_1 ( U_{U^{n-1}_F / J^{n-1}_{N , F}} )\, $. Divide diagram (2)  by the subgroup
\smallskip
$$ \bigl[\,\bigl[\cdots \bigl[\, J_{U_{J^{n-1}_{N , F}\, ,\, U^{n-1}_F}\, ,\, 
U^n_F}\, ,\, U^{n+1}_F\,\bigr] ,\cdots\bigr] , U^{n+1}_F\,\bigr]\, +\, $$
$$ \widetilde C\> \cap\> \bigl[\,\bigl[\cdots\bigl[\, U_{U_{J^{n-1}_{N , F}\, ,\, U^{n-1}_F}\, ,\, U^n_F}\, ,\, 
\widetilde U\,\bigr]\, ,\cdots\bigr] , \widetilde U\,\bigr]\, +\, $$
$$  \bigl[\,  \sigma (\, J_{U^n_F\, /\, J_{J^{n-1}_{N , F}\, ,\, U^{n-1}_F}}\, )\, ,\, 
U_{U_{J^{n-1}_{N , F}\, ,\, U^{n-1}_F}\, ,\, U^n_F}\,\bigr]\, +\, $$
$$ \bigl[\,J^{n+1}_F\, ,\, U_{J_{J^{n-1}_{N , F}\, ,\, U^{n-1}_F}\, ,\, U^n_F}\,\bigr]\, +\,
\bigl[\, J_{J_{J^{n-1}_{N , F}\, ,\, U^{n-1}_F}\, ,\, U^n_F}\, ,\, U^{n+1}_F\,\bigr]  $$
\par\smallskip\noindent
where the the first and second expression is supposed to mean k-fold commutators for some fixed value of $\, k \geq 1\, $. Note that since 
$\, U_{U_{J^{n-1}_{N , F}\, ,\, U^{n-1}_F}\, ,\, U^n_F}\, $ is normal and since 
$\, \sigma (\, J_{U^n_F\, /\, J_{J^{n-1}_{N , F}\, ,\, U^{n-1}_F}}\, )\, $ is normal for 
$\, \sigma (\, U_{U_{U^{n-1}_F\, /\,  J^{n-1}_{N , F}\,}}\, )\, $ which is the quotient of $\, U^{n+1}_F\, $ modulo  $\, U_{U_{J^{n-1}_{N , F}\, ,\, U^{n-1}_F}\, ,\, U^n_F}\, $ the commutator subgroup of
$\, U^{n+1}_F\, $ above is normal. Its intersection with 
$\, J_{J_{J^{n-1}_{N , F}\, ,\, U^{n-1}_F}\, ,\, U^n_F}\, $ is equal to
\smallskip 
$$ J_{J_{J^{n-1}_{N , F}\, ,\, U^{n-1}_F}\, ,\, U^n_F}\> \cap\> 
\Bigl(\,\bigl[\, U_{J_{J^{n-1}_{N , F}\, ,\, U^{n-1}_F}\, ,\, U^n_F}\, ,\, 
U_{J_{J^{n-1}_{N , F}\, ,\, U^{n-1}_F}\, ,\, U^n_F}\,\bigr]\> +\qquad $$ 
$$ \qquad\qquad\qquad\bigl[\,\bigl[\cdots\bigl[\, U_{J_{J^{n-1}_{N , F}\, ,\, U^{n-1}_F}\, ,\, U^n_F}\, ,\, 
U^{n+1}_F\,\bigr] ,\cdots\bigr] , U^{n+1}_F\, \bigr]\,\Bigr)  $$
\par\bigskip\noindent 
modulo the trivial fourth and fifth relations as can be seen replacing all elements in the first three relations by the corresponding elements in $\,\sigma ( U_{U^n_F / J_{J^{n-1}_{N , F}\, \, U^{n-1}_F}} )\, $ and calculating the difference (compare the remark at the beginning of the proof). 
Consider the image of 
$\, U_{J^{n-1}_{N , F}\, ,\, U^n_F / J_{J^{n-1}_{N , F}\, ,\, U^{n-1}_F}}\, $ in the bottom line of (2). We have divided out the action of 
$\,  J_{U^n_F\, /\, J_{J^{n-1}_{N , F}\, ,\, U^{n-1}_F}}\, $ on this subgroup. 
The difficulty now consists in showing that it is possible to make this preimage into an k-fold iterated 
$U^{n-1}_F$-central extension of $\, J^{n-1}_{N , F}\, $ on dividing by a suitable free abelian subgroup in the quotient of 
$\, J^{n+1}_F\, $, i.e. reducing the action of $\, U^n_F\, $ to an action of $\, U^{n-1}_F\, $ in the bottom line.  Since $\, U^n_F\, =\, J^n_F \rtimes U^{n-1}_F\, $ (at least for $\, n\geq 2\, $) we have to worry about the action of $\, J^n_F\, $ on the preimage of $\, J^{n-1}_{N , F}\, $ which is induced by its action on the preimage of $\, U_{J^{n-1}_{N , F}\, ,\, U^{n-1}_F}\, $.
We have to check what part of $\, C\, $ remains in the quotient as above. It is clear that the conjugacy coefficients 
$\, v^{\pm}_{l_1}\cdots v^{\pm }_{l_m}\, $ are divided out since the basic commutator in the middle lies in 
$\, J_{U^n_F / J_{J^{n-1}_{N , F}\, ,\, U^{n-1}_F}}\, $. Similarly, since the element $\, v_{l_0}\, $ of the basic commutator is in $\, s_1 ( U_{J^{n-1}_{N , F}\, ,\, U^{n-1}_F} )\, $ the coefficients 
$\, u^{\pm }_{j_1}\cdots u^{\pm }_{j_n}\, $ and $\, u_{j_0}\, $ are reduced to their images in 
$\, P^{n-1}_{N , F}\, $ so that the basis of $\, C\, $ as above is reduced to the form
\smallskip
$$\bigl\{\, l^{\pm }_1\cdots l^{\pm }_m\, \bigl[\, v^{\pm }\, ,\, l_0\,\bigr]\, l^{\mp }_m\cdots l^{\mp }_1\,\bigr\}
\leqno{(*)} $$
\par\smallskip\noindent
with $\, v\, $ a basis element of $\, s_1 ( U_{J^{n-1}_{N , F}\, ,\, U^{n-1}_F} )\, $ and 
$\, l_0 , l_1,\cdots , l_n\, $ are taken from a basis of $\, P^{n-1}_{N , F}\, $ which is a free group for all 
$\, n \geq 1\, $. Of course one gets an abelian group as we have divided by 
$\, [\, J^{n+1}_F\, ,\, J_{U_{J^{n-1}_{N , F}\, ,\, U^{n-1}_F}\, ,\, U^n_F}\, ]\, $, then $\, C\, $ must be a quotient of the free abelian group obtained from this basis. One easily finds that by this replacement of the basis the relations coming from the third term of the commutator subgroup of $\, U^{n+1}_F\, $ that we have divided out are exhausted so only the first and second remain. Moreover the image of the first relation in $\, U_{U^n_f / J_{J^{n-1}_{N , F}\, ,\, U^{n-1}_F}}\, $ is reduced to 
\smallskip
$$ \bigl[\,\bigl[\cdots\bigl[\, J_{J^{n-1}_{N , F}\, ,\, U^n_F / J_{J^{n-1}_{N , F}\, ,\, U^{n-1}_F}}\, ,\, 
s_1 ( U_{U^{n-1}_F / J^{n-1}_{N , F}} )\,\bigr] ,
\cdots\bigr] , s_1( U_{U^{n-1}_F / J^{n-1}_{N , F}} )\,\bigr]\, $$
\par\bigskip\noindent
modulo the third and its intersection with the image of $\, C\, $ is contained in the image of the second relation.  We have to show two things. First, that the quotient is a free abelian group, and second, that it has a normal lift to the image of $\, J^{n+1}_F\, $. We first do the case $\, n = 1\, $, We want to show that we can reduce the basis elements 
$\, v \in s_1 ( U_{J^{n-1}_{N , F}\, ,\, U^{n-1}_F} )\, $ to basis elements in the core for $\, s_1 ( U^n_F )\, $. To see this assume that $\, v\, $ is of the form 
$\, t\, u\, t^{-1}\, $ with $\, u\, $ in the core and $\, t = s_1 ( y )\, $ with $\, y\, $ in the image of the splitting induced by $\, {\nu }_1\, $. Then 
\smallskip
$$l^{\pm }_1\cdots l^{\pm }_m\,\bigl[\, t\, u^{\pm }\, t^{-1}\, ,\, l_0\,\bigr]\, l^{\mp }_m\cdots l^{\mp }_1\, =\,$$
$$ t\, ( t^{-1}\, l^{\pm }_1\, t )\cdots ( t^{-1}\, l^{\pm }_m\, t )\, \bigl[\, u^{\pm }\, ,\, ( t^{-1}\, l_0\, t )\,\bigr]\, 
( t^{-1}\, l^{\mp }_m\, t )\cdots ( t^{-1}\, l^{\mp }_1\,  t )\, t^{-1} $$ 
\par\smallskip\noindent
and since $\ t\, $ is represented by commutators of arbitrary length modulo $\, J_{F / N}\, $ as $\, F / N\, $ is perfect
the corresponding basis element can be expanded into elements of the form 
$\, \{\, l^{\pm }_1\cdots l^{\pm }_m\, [\, u\, ,\, l_0\, ] l^{\mp }_m\cdots  l^{\mp }_1\, \}\, $ modulo 
$\, [\, J_{U^n_F\, /\, J_{J^{n-1}_{N , F}\, ,\, U^{n-1}_{N , F}}}\, ,\,
U_{J^{n-1}_{N , F}\, ,\, U^n_F\, /\, J_{J^{n-1}_{N , F}\, ,\, U^{n-1}_F}}\, ]\, $ which is contained in the image of the relations we have divided out. 
Replace the basis by the following commutator basis
\smallskip
$$\bigl\{\,\bigl[\,\bigl[\cdots\bigl[\, \bigl[\, u^{\pm }\, ,\, l_0\,\bigr] , l^{\pm }_1\,\bigr] ,
\cdots\bigr] , l^{\pm }_m\,\bigr]\,\bigr\} \leqno{(**)} $$
\par\smallskip\noindent
The next thing is that the elements $\,\{ l_k \}_{k\geq 0}\, $ can be replaced by elements from a set of representatives for the basis of the free abelian group 
$\, J_{F / N}\, /\, (\, J_{F / N}\,\cap\, [\, U_{F / N}\, ,\, U_{F / N}\, ]\, $ (which of course is isomorphic with 
$\, U_{F / N}\, /\, [\, U_{F / N}\, ,\, U_{F / N}\, ]\, $ for perfect $\, F / N\, $). Indeed, suppose that a m-fold commutator $\, \bigl[\,\bigl[\cdots\bigl[\, \bigl[\, u^{\pm }\, ,\, l_{j_0}\,\bigr] , l^{\pm }_{j_1}\,\bigr] ,
\cdots\bigr] , l^{\pm }_{j_m}\,\bigr]\, $ is such that
$\, l_{j_i}\, =\, [\, t_1\, ,\, t_2\, ]\, $ with $\, t_1 , t_2 \in U_{F / N}\, $. Then it can be written as a commutator expression of arbitrary length modulo 
$\, [\, J_{F / N}\, ,\, J_{F / N}\, ]\, $ so it can be expanded into elements of the basis above with commutator length greater than $\, m+1\, $. Since the total commutator length is bounded by the second (and first) expression of our relations one gets the result. This argument is also valid in the variable $\, l_0\, $ 
from the second relation which is thereby clearly exhausted.
Thus $(**)$ is precisely reduced to the subbasis with length 
$\, m \leq k\, $, i.e. to the free abelian group generated by such basis elements. This subgroup can be lifted normally to the image of $\, J^{n+1}_F\, $ and divided out.
To see this one notes that it suffices to show that the lift is normal for the image of $\, U_{J_F\, ,\, U_F}\, $ because modulo $\, J_{U_F}\, +\, U_{J_F\, ,\, U_F}\, $ any element can be replaced by a commutator of arbitrary length and we have divided by $\, \bigl[\, \bigl[\,\bigl[\cdots \bigl[\, J^{n+1}_F\, ,\, 
U^{n+1}_F\,\bigr] ,\cdots\bigr] , U^{n+1}_F\,\bigr] , U^{n+1}_F\,\bigr]\, $. Lifting a corresponding subbasis of $\, C\, $ by lifting the basic commutators in the middle to $\, J^{n+1}_F\, $ and completing normally with respect to some preimages of the $\,\{\, u_{l_k}\,\}_{k\geq 1}\, $ the lift will be normal for some lift of 
$\, U_{J_{F / N} }\, $, hence normal because it is normal in any case for 
$\, U_{J_{J^{n-1}_{N , F}\, ,\, U^{n-1}_F}\, ,\, U^n_F}\, $ by the fifth term of our relation subgroup. Next we do the case $\, n = 2\, $. One gets that $\, P_{N , F}\, =\, J_{U_F / J_{N , F}}\, $ with 
$\, U_F\, /\, J_{N , F}\, =\, N\,\rtimes\, U_{F / N}\, $. Since $\, N\, $ is a full subgroup of 
$\, U_F\, /\, J_{N , F}\, $ one gets by an argument as above that the elements $\, v\, $ in the basis $(*)$ are reduced to the form $\, v\, =\, s_1 ( x )\, u^{\pm }\, s_1 ( x )^{-1}\, $ where $\, u\, $ is in the core and 
$\, x \in {\nu }_1 ( U_{F / N} )\, $. Then considering the preimage in $\, U_{U_F\, /\, J_{N , F}}\, $ of the copy of $\, U_{F / N}\, $ in the decomposition $\, U_F\, /\, J_{N , F}\, =\, N\,\rtimes\, U_{F / N}\, $ denote this preimage by $\, U\, $. Then one can replace the basis $(*)$ as above by a corresponding basis where 
$\, v\, =\, s_1 ( x )\, u\, s_1 ( x )^{-1}\, $ is replaced by the basis element $\, u\, $ of the core, $\, l_0\, $ is still a basis element of $\, P_{N , F}\, $, but the outer coefficients $\, l_1 ,\cdots , l_m\, $ are replaced by elements of a basis of $\, U\, $ modulo $\, [\, U_{U_F / J_{N , F}}\, ,\, U_{U_F / J_{N , F}}\, ]\, $ (which is again equal to $\, U_{U_F\, /\, J_{N , F}}\, /\, [\, U_{U_F\, /\, J_{N , F}}\, ,\, U_{U_F\, /\, J_{N , F}}\, ]\, $). This replacement corresponds to division of the part $\, m = 0\, $ of the second relation. By an argument as above one finds that the image of $\, C\, $ modulo the second (and first) relation is equal to its image in the free abelian subgroup  with commutator basis as in $(**)$ with $\, m \leq k\, $ but where the elements 
$\, l_k\, $ are representatives of a basis for $\, U\, $ modulo 
$\, [\, U_{U_F / J_{N , F}}\, ,\, U_{U_F / J_{N , F}}\, ]\, $ (including the element $\, l_0\, $).  Again this is a free abelian group and it can be lifted normally to $\, J^{n+1}_F\, $ and divided out. The case $\, n \geq 3\, $ is much the same. By induction assumption one gets a normal decomposition of the preimage of 
$\, R^{n-1}_{N , F}\, $ in the image of $\, U_{J^{n-1}_{N , F}\, ,\, U^{n-1}_F}\, $ in the bottom line as 
\smallskip
$$ A\,\rtimes\, {\overline R}^{n-1}_{N , F} $$
\par\smallskip\noindent
where $\, {\overline R}^{n-1}_{N , F}\, $ is some minimal $U^{n-1}_F$-central extension of 
$\, R^{n-1}_{N , F}\, $ ( image of the relative universal extension) and $\, A\, $ is a (k+1)-nilpotent subgroup for the action of $\, U^{n-1}_F\, $. The decomposition 
$\, U^{n+1}_F = J^{n+1}_F \rtimes U^n_F\, $ can be adjusted so that 
$\, U^n_{N , F}\, $ maps onto the factor 
$\,  {\overline R}^{n-1}_{N , F}\, $ and $\, U_{J^{n-1}_{N , F}\, ,\, U^{n-1}_F}\, $ maps to the image of 
$\, U_{U_{J^{n-1}_{N , F}\, ,\, U^{n-1}_F}\, ,\, U^n_F}\, $. This implies by the relations we have divided out that the subgroup of $\, J_{J^{n-1}_{N , F}\, ,\, U^{n-1}_F}\, $ intersecting with the image of our relations maps to the corresponding canonical subgroup of 
$\, U_{J_{J^{n-1}_{N , F}\, ,\, U^{n-1}_F}\, ,\, U^n_F}\, $ modulo 
the fourth and fifth relation and the relation
\smallskip
$$  J_{J_{J^{n-1}_{N , F}\, \, U^{n-1}_F}\, ,\, U^n_F}\>\cap\>
\bigl(\, \bigl[\,U_{J_{J^{n-1}_{N , F}\, ,\, U^{n-1}_F}\, ,\, U^n_F}\, ,\, 
U_{J_{J^{n-1}_{N , F}\, ,\, U^{n-1}_F}\, ,\, U^n_F}\,\bigr] \qquad\qquad $$
$$ \qquad\qquad\qquad +\> \bigl[\,\bigl[\cdots\bigl[\, U_{J_{J^{n-1}_{N , F}\, ,\, U^{n-1}_F}\, ,\, U^n_F}\, ,\, 
U^{n+1}_F\,\bigr] ,\cdots\bigr] , U^{n+1}_F\,\bigr]\, \bigr)\>   $$
\par\bigskip\noindent
where the last bracket denotes $k$-fold commutators. The preimage of the image of $\, R^n_{N , F}\, $ in the top line (divided only by these relations) leads to a normal decomposition of the preimage of 
$\, R^n_{N , F}\, $ in the top line as 
\smallskip
$$ D\,\times\, R^n_{N , F}\> .  $$
\par\smallskip\noindent
If $\, n = 1\, $ this is what we want since $\, J^0_{N , F} = N\, $.  Otherwise 
consider the image of $\, J_{J_{J^{n-1}_{N , F}\, ,\, U^{n-1}_F}\, ,\, U^n_F}\, $ in the quotient modulo 
$\, U_{J^n_{N , F}\, ,\, U^n_F}\, $. It is of the form
\smallskip 
$$ J_{J_{C , B}\, ,\, \widetilde U_B}\, /\, \bigl(\,
\bigl[\,J_{J_{C , B}\, ,\, \widetilde U_B}\, ,\, U_{\widetilde U_B}\,\bigr]\, +\, 
\bigl[\, U_{J_{C , B}\, ,\, \widetilde U_B}\, ,\, J_{\widetilde U_B}\, \bigr]\,\bigr)\> . $$
\par\smallskip\noindent
with $\, C\, =\, J^{n-1}_{N , F} / R^{n-1}_{N , F}\, $ a central subgroup of 
$\, B\, =\, U^{n-1}_F / R^{n-1}_{N , F}\, $ and $\, \widetilde U_B\, $ an extension of $\, U_B\, $ by a copy of $\, R^{n-1}_{N , F}\, $ complementary to $\, J_{C , B}\, $.  ( clearly one gets a quotient of the group above but it  follows from results of the following sections  --- Theorem 5 and Lemma 15 of section 6 which states that $\, K^J_2 ( J_{C , B }\, ,\, U_B )\, =\, 0\, $ --- that the image is precisely this group ). This quotient is an extension of the free abelian group
\smallskip
$$ J_{J_{C , B}\, ,\, \widetilde U_B}\, /\, \bigl(\, J_{J_{C , B}\, ,\, \widetilde U_B}\,\cap\, 
\bigl[\,U_{J_{C , B}\, ,\, \widetilde U_B}\, ,\, U_{\widetilde U_B}\,\bigr]\,\bigr)   $$
\par\smallskip\noindent
by the quotient
\smallskip
$$ K^J_2 ( J_{C , B}\, ,\, \widetilde U_B )\> =\> {J_{J_{C , B}\, ,\, \widetilde U_B}\, \cap\, 
\bigl[\,U_{J_{C , B}\, ,\, \widetilde U_B}\, ,\, U_{\widetilde U_B}\,\bigr] \over
\bigl[\,J_{J_{C , B}\, ,\, \widetilde U_B}\, ,\, U_{\widetilde U_B}\,\bigr]\, +\, 
\bigl[\, U_{J_{C , B}\, ,\, \widetilde U_B}\, ,\, J_{\widetilde U_B}\, \bigr] }\>  . $$
\par\smallskip\noindent
By Lemma 7 $\, J_{C , B}\, /\, [\, J_{C , B}\, ,\, J_{C , B}\, ]\, $ has a core for a lift of the basis of 
$\, U_{B / C}\, $ to 
a subbasis of the basis of $\, U_B\, $ (which is again a subgroup of 
$\,\widetilde U_B = R^{n-1}_{N , F}\rtimes U_B\, $) and there is a lift of a copy of $\, U^{n-1}_F\, $ into 
$\,\widetilde U_B\, $ implementing the usual action on $\, R^{n-1}_{N , F}\, $ and such that the lift is contained in the image of $\, U_{B / C}\, $ (since both are sort of diagonal, compare with section 3). It follows from the proof of Lemma 14 that $\, J_{C , B} / [\, J_{C , B}\, ,\, J_{C , B}\, ]\, $ has a normal subgroup $\, J_0\, $ with quotient equal to $\, J_C / [\, J_C\, ,\, J_C\, ]\, $ such that $\, J_0\, $ has a core for $\, U_B / J_{C , B}\, $ whence the group 
$\, K^J_2 ( \widetilde J_0\, ,\, \widetilde U_B )\, $ where $\, \widetilde J_0\, $ is the preimage of $\, J_0\, $ in $\, J_{C , B}\, $ which is contained in $\, K^J_2 ( J_{C , B}\, ,\, \widetilde U_B )\, $ by Lemma 10 is free abelian by Lemma 11. The quotient $\, J_C / [\, J_C\, ,\, J_C\, ]\, $ can be factored again into the image of $\, [\, U_C\, ,\, U_C\, ]\, $ and a complementary free abelian group. The action of $\, U_B\, $ on these groups reduces to an action of $\, C\, $. The subgroup 
$\, A_1 =  [\, U_C\, ,\, U_C\, ] / [\, J_C\, ,\, J_C\, ]\, $ has a basis consisting of the set 
\smallskip
$$ \bigl\{\, u_e\,\bigl[ u_c\, ,\, u_d\,\bigr]\, u^{-1}_e\>\vert\> c < d , e\leq d \bigr\}  $$
\par\bigskip\noindent
and the subgroup generated by elements with $\, e = 1\, $ can be identified with 
$\, A_1 / [\, A_1\, ,\, C\, ]\, $.  The action on the complement 
$\, A_2\, $ is trivial modulo $\, A_1\, $. It follows from Lemma 10 that the map 
$\, K^J_2 ( \widetilde J_0\, ,\, \widetilde U_B ) \rightarrow K^J_2 ( J_{C , B}\, ,\, \widetilde U_B )\, $ is injective and the quotient is contained in the kernel of the map 
$\, K^J_2 ( J_C / [\, J_C\, ,\, J_C\, ]\, ,\, \widetilde U_B / \widetilde J_0 ) \rightarrow 
K^J_2 ( J_C / [\, J_C\, ,\, J_C\, ]\, ,\,  U_B / \widetilde J_0 )\, $. Now one can choose a copy of $\, U_C\, $ inside the image of $\,\widetilde U_B\, $ implementing the action of $\, C\, $ on the image of $\, J_C\, $ but acting trivially on the copy of $\, R^{n-1}_{N , F}\, $ (by the fact that the action on $\, R^{n-1}_{N , F}\, $ was implemented by the action of the copy of $\, U_{B / C}\, $ we have trivialized on $\, J_{C , B}\, $).   We may further divide the copy of $\, R^{n-1}_{N , F}\, $ by the subgroup 
$\, [\, R^{n-1}_{N , F}\, ,\, U^{n-1}_{N , F}\, ]\, $ because taking elements in the right argument of the bracket from the copy of $\, U_{B / C}\, $ implementing the usual action on 
$\, R^{n-1}_{N , F}\, $ but commuting with $\, A_1\, $ and $\, A_2\, $ its intersection with the kernel as above is clearly trivial. Consider the image of $\, [\, J^{n-1}_{N , F}\, ,\, U^{n-1}_F\, ]\, $ in the quotient. Modulo $\, A_1\, $ it is contained in the commutator subgroup of $\, U_{B / C}\, $ commuting with 
$\, A_1 \times A_2\, $. Clearly this gives a homomorphism of this image into $\, A_1\, $ and we may divide by its kernel since it consists of commutators of elements commuting with $\, A_1 + A_2\, $. The remaining quotient $\, D_1\, $ is free abelian since it is isomorphic to a subgroup of $\, A_1\, $. Also by induction we may assume that $\, D_2 = R^{n-1}_{N , F} / [\, J^{n-1}_{N , F}\, ,\, U^{n-1}_F\, ]\, $ is free abelian. The remaining part of our kernel is now contained in 
$\, K^J_2 ( D_1 \times D_2\, ,\, \overline U_B )\, $ where 
$\,\overline U_B = ( D_1 \times D_2 ) \times ( U_B / \widetilde J_0 )\, $. Since $\, D_1 \times D_2\, $ is central and free abelian Lemma 17 implies that $\, K^J_2 ( D_1 \times D_2\, ,\, \overline U_B )\, $ is a direct product of the tensor product of $\, D_1 \times D_2\, $ with the abelianization of 
$\, U_B / \widetilde J_0\, $ and $\, K^J_2 ( D_1 \times D_2 )\, $.  Since the abelianization of 
$\, U_B / \widetilde J_0\, $ is a free abelian group isomorphic to $\, A_2 \times ( U_{B / C} )^{ab}\, $ one obtains a free abelian group. 
Thus by a sequence of regular surjections each of which corresponding to a free abelian group we have reduced $\, K^J_2 ( J_{C , B}\, \widetilde U_B )\, $ (which is equal to the kernel of the regular surjection 
$\, K^J_2 ( J_{C , B}\, ,\, \widetilde U_B ) \twoheadrightarrow K^J_2 ( J_{C , B}\, ,\, U_B )\, $) to zero, implying that it is free abelian. Lifting the free abelian group
\smallskip
$$ {J_{J_{C , B}\, ,\, \widetilde U_B}\over 
\bigl[\,J_{J_{C , B}\, ,\, \widetilde U_B}\, ,\, U_{\widetilde U_B}\,\bigr]\, +\, 
\bigl[\, U_{J_{C , B}\, ,\, \widetilde U_B}\, ,\, J_{\widetilde U_B}\, \bigr]}  $$
\par\bigskip\noindent
to the image of $\, J_{J_{J^{n-1}_{N , F}\, ,\, U^{n-1}_F}\, ,\, U^n_F}\, $ and dividing out its image will project the normal copy of $\, R^n_{N , F}\, $ into the image of 
$\, U_{J^n_{N , F}\, ,\, U^n_F}\, $ modulo 
\smallskip
$$  J_{J^n_{N , F}\, ,\, U^n_F}\>\cap\>
\bigl(\, \bigl[\,U_{J^n_{N , F}\, ,\, U^n_F}\, ,\, 
U_{J^n_{N , F}\, ,\, U^n_F}\,\bigr] \qquad\qquad $$
$$ \qquad\qquad\qquad +\> \bigl[\,\bigl[\cdots\bigl[\, U_{J^n_{N , F}\, ,\, U^n_F}\, ,\, 
U^{n+1}_F\,\bigr] ,\cdots\bigr] , U^{n+1}_F\,\bigr]\, \bigr)\>   $$
\par\bigskip\noindent 
plus trivial elements and representatives of $\, K^J_2 ( J^n_{N , F}\, ,\, U^n_F )\, $ which become trivial in 
$\, K^J_2 ( J_{J^{n-1}_{N , F}\, ,\, U^{n-1}_F}\, ,\, U^n_F\, )\, $. As follows from Lemma 15 and Corollary 5.1 the map $\, K^J_2 ( J^n_{N , F}\, ,\, U^n_F ) \rightarrow 
K^J_2 ( J_{J^{n-1}_{N , F}\, ,\, U^{n--1}_F}\, ,\, U^n_F )\, $ is injective. Hence the canonical $U^n_F$-central extension of $\, J^n_{N , F}\, $ injects into the canonical $U^n_F$-central extension of 
$\, J_{J^{n-1}_{N , F}\, ,\, U^{n-1}_F}\, $ and one obtains a normal copy of $\, R^n_{N , F}\, $ in the former modulo the image of the commutator subgroup above. By construction this copy contains the image of 
$\, [\, J^n_{N , F}\, ,\, U^n_F\, ]\, $ which implies that $\, R^n_{N , F} / [\, J^n_{N , F}\, ,\, U^n_F\, ]\, $ is free abelian (compare the Remark before Lemma 14 of section 5). 
Its preimage modulo
$\, [\, J_{J^n_{N , F}\, ,\, U^n_F}\, ,\, U^{n+1}_F\, ]\, +\, [\, U_{J^n_{N , F}\, ,\, U^n_F}\, ,\, J^{n+1}_F\, ]\, $
is then a normal $U^n_F$-central extension of $\, R^n_{N , F}\, $ which will be denoted 
$\, {\overline {R^n_{N , F}}}^{ _{\sps U^n_F}}\, $. It contains a canonical copy of 
$\, [\, J^n_{N , F}\, ,\, J^n_{N , F}\, ]\, $ and modulo this subgroup the kernel is contained in the image of 
$ [\, [\cdots [\, {\overline {R^n_{N , F}}}^{ _{\sps U^n_F}}\, ,\, U^n_F\, ] ,\cdots ] , U^n_F ]\, $. It is easy to see that $\, {\overline {R^n_{N , F}}}^{ _{\sps U^n_F}}\, $ satisfies (the relative version of) property (i) of Definition 2 (c) (the uniqueness part follows from the fact that the kernel is contained in 
$\, [\, {\overline {R^n_{N , F}}}^{ _{\sps U^n_F}}\, ,\, U^n_F\, ]\, $). To complete the induction step 
$\, n-1 \rightarrow n\, $ one must show that it also satisfies (the relative version of) property (ii). Let then 
$\,{ \widetilde R}^n\, $ be a $U^n_F$-central extension of $\, {\overline {R^n_{N , F}}}^{ _{\sps U^n_F}}\, $ which extends to a double $U^n_F$-central extension of $\, J^n_{N , F}\, $. Consider the preimage of a normal copy of $\, {\overline {R^n_{N , F}}}^{ _{\sps U^n_F}}\, $ in the canonical 
$U_F$-central extension of $\, J_{N , F}\, $ in the canonical double $U_F$-central extension of 
$\, J_{N , F}\, $ which is
\smallskip
$$ U_{J^n_{N , F}\, ,\, U^n_F}\, /\, ( \bigl[\, \bigl[\, J_{J^n_{N , F}\, ,\, U^n_F}\, ,\, U^{n+1}_F\,\bigr] , 
U^{n+1}_F\,\bigr]\, +\,
  \bigl[\, J^{n+1}_F\, ,\, U_{J^n_{N , F}\, ,\, U^n_F}\,\bigr] )\>  .  $$
\par\smallskip\noindent
We claim that $\, U^n_F\, $ acts trivially on the preimage of $\, C^n\, $, the kernel of 
$\, {\overline {R^n_{N , F}}}^{ _{\sps U^n_F}} \rightarrow R^n_{N , F}\, $. Dividing this preimage by its intersection with the image of 
\smallskip
$$ \bigl[\, J^{n+1}_{N , F}\, ,\, U^{n+1}_F\,\bigr]\>\cap\> 
\bigl[\,\bigl[\, U^{n+1}_{N , F}\, ,\, U^{n+1}_F\, \bigr] , U^{n+1}_F\,\bigr] $$
\par\smallskip\noindent
one gets a splitting of the form $\, D \times {\overline {R^n_{N , F}}}^{ _{\sps U^n_F}}\, $ by the universal property (i) of the pair $\, ( {\overline {R^n_{N , F}}}^{ _{\sps U^n_F}}\, ,\, U^{n+1}_F )\, $ (and since modulo 
$\, [\, U^{n+1}_{N , F}\, ,\, J^{n+1}_F\, ]\, $ one has
\smallskip 
$$ J^{n+1}_{N , F}\,\cap\, [\, J_{J^n_{N , F}\, ,\, U^n_F}\, ,\, U^{n+1}_F\, ]\, =\, 
[\, J^{n+1}_{N , F}\, ,\, U^{n+1}_F\, ]\, $$ 
\par\smallskip\noindent
by injectivity of $\, K^J_2 ( R^n_{N , F}\, ,\, U^n_F )\rightarrow K^J_2 ( J^n_{N , F}\, ,\, U^n_F )\, $ due to Lemma 14 and 
\smallskip
$$ U^{n+1}_{N , F}\,\cap\, [\, [\, U_{J^n_{N , F}\, ,\, U^n_F}\, ,\, U^{n+1}_F\, ] , U^{n+1}_F\, ]\, =\,
[\, [\, U^{n+1}_{N , F}\, ,\, U^{n+1}_F ] , U^{n+1}_F\, ]\> ) . $$
\par\smallskip\noindent
Denote $\, \overline {{\overline R}^n}\, $ the preimage of the copy of $\, {\overline R}^n_{N , F}\, $ as above in 
$$ U_{J^n_{N , F}\, ,\, U^n_F}\, /\, ( \bigl[\, \bigl[\, J_{J^n_{N , F}\, ,\, U^n_F}\, ,\, U^{n+1}_F\,\bigr] , 
U^{n+1}_F\,\bigr]\, +\,
  \bigl[\, J^{n+1}_F\, ,\, U_{J^n_{N , F}\, ,\, U^n_F}\,\bigr] )\>  . $$
\par\smallskip\noindent
This is a minimal $U^n_F$-normal subextension. Let $\, {\overline C}^n\, $ be the kernel of 
$\, \overline {{\overline R}^n} \rightarrow {\overline {R^n_{N , F}}}^{ _{\sps U^n_F}}\, $ and $\, \overline {{\overline C}^n}\, $ the (abelian) kernel of $\, \overline {{\overline R}^n}\rightarrow R^n_{N , F}\, $. Then 
$\, {\overline C}^n\, $ is in the image of $\, [\, [\, U^{n+1}_{N , F}\, ,\, U^{n+1}_F\, ] , U^{n+1}_F\, ]\, $ (actually $\, \overline {{\overline C}^n}\, $ is). It is clear that $\,\overline {{\overline C}^n}\, $ is central for the action of $\, J^n_{N , F}\, $. Assume first that $\, N = F\, $ and $\, n = 1\, $. Then the action of 
$\, U_F\, $ on $\, \overline {\overline {C^1}}\, $ reduces to an action of $\, \overline F\, =\, U_F / R_F\, $ which is a perfect group and hence must act trivially on $\, \overline {\overline {C^1}}\, $, settling this special case.. The general case is more difficult. Put 
$\, \widetilde {U^n}\, =\, J^n_{N , F}\, +\, [\, U^n_F\, ,\, U^n_F\, ]\, $
and $\, R^{n-1}_{0 , N , F}\, :=\, R^{n-1}_{N , F}\,\cap\, [\, U^{n-1}_F\, ,\, U^{n-1}_F\, ]\, $. Since 
the pair $\, ( R^{n-1}_{0 , N , F}\, ,\, U^{n-1}_F )\, $ admits a relative semiuniversal $U^{n-1}_F$-central extension with kernel isomorphic to the kernel of $\, {\overline R}^{n-1}_{N , F}\, $ ( compare Proposition 1 of section 6 ) and since $\, R^{n-1}_{N , F}\, /\, R^{n-1}_{0 ,N , F}\, $ is free abelian, the construction above slightly modified gives that the pair $\, ( R^n_{N , F}\, ,\, \widetilde {U^n} )\, $ admits a semiuniversal 
$\,\widetilde{U^n}$-central extension relative to $\, ( J^n_{N , F}\, ,\, \widetilde {U^n} )\, $ which by the Special Deconstruction Theorem (Lemma 16 of section 6) must be isomorphic to $\, {\overline R}^n_{N , F}\, $. To see this consider the image of 
$\, U_{\widetilde {U^n}} / \sim\, $ in the decomposition 
$\, A\,\times\, ( R^{n-1}_{N , F}\rtimes U_{U^{n-1}_F\, /\, R^{n-1}_{N , F}} )\, $ of the bottom line of (2) on dividing by 
$\, [\, J^{n+1}_F\, ,\, U^{n+1}_F\, ]\, $ and note that the intersection of this commutator subgroup with $\, U_{\widetilde {U^n}}\, $ is equal to 
$\, [\, J_{\widetilde {U^n}}\, ,\, U_{\widetilde {U^n}}\, ]\, $. The quotient
$\, U^{n+1}_F / U_{\widetilde {U^n}}\, $ is of the form $\, U_C\, $ with $\, C\, $ free abelian, so that 
$\, J_C / [\, J_C\, ,\, U_C\, ]\, $ is free abelian and hence splits normally the direct factor $\, A\, $ 
( $\, U_C / [\, U_C\, ,\, U_C\, ]\, $ is free abelian and also $\, [\, U_C\, ,\, U_C\, ] / [\, J_C\, ,\, U_C\, ]\, $ since the latter is equal to $\, K^J_2\, ( C)\, $ --- compare with section 5 below --- which in case of a free abelian group is again free abelian, then $\, J_C / [\, J_C\, ,\, U_C\, ]\, $ is an abelian extension of a free abelian group by another free abelian group, hence free abelian). Now there is a projection of the subgroup 
$\,  \widetilde U\, /\, R^n_{N , F}\, $ in the copy of 
$\, R^{n-1}_{N , F}\rtimes U_{U^{n-1}_F\, /\, R^{n-1}_{N , F}}\, $
to the image of $\, U_{\widetilde {U^n}}\, $ and one easily constructs a lift of 
$\, \widetilde {U^n}\, $ to $\, U_{\widetilde {U^n}} / \sim\, $ mapping onto this normal copy in the bottom line so one gets that the pair $\, ( R^n_{N , F}\, ,\, \widetilde U^n )\, $
admits a relative semiuniversal $\widetilde U^n$-central extension (for the inclusion into 
$\, ( J^n_{N , F}\, ,\, \widetilde U^n )\, $). Then,  from the natural map 
$\, ( R^n_{N , F}\, ,\, \widetilde {U^n} )\,\subset\, ( R^n_{N , F}\, ,\, U^n_F )\, $ one gets that the kernel 
$\, C^n\, $ of $\, {\overline {R^n_{N , F}}}^{ _{\sps U^n_F}}\, $ is contained in 
$\,  [\, {\overline {R^n_{N , F}}}^{ _{\sps U^n_F}}\, ,\, \widetilde {U^n}\, ]\, $. On the other hand the subgroup $\,\widetilde {U^n}\, $ of $\, U^n_F\, $ acts trivially on the kernel 
$\, \overline {{\overline C}^n}\, $ of the double central extension $\, \overline {{\overline R}^n}\, $, so there is a splitting of the latter as 
$\, {\overline C}^n \times {\overline {R^n_{N , F}}}^{ _{\sps U^n_F}}\, $ which is normal for the action of 
$\,\widetilde {U^n}\, $. This must also be normal for the action of $\, U^n_F\, $ for otherwise there would be an element $\, c\, $ of the kernel 
$\, C^n \subseteq \overline{{\overline C}^n} = {\overline C}^n \times C^n\, $
such that $\, [\, c\, ,\, u_x\, ]\, \in {\overline C}^n\, $ is nonzero, and $\, c\, $ can be represented by a product of elements  $\, [\, r\, ,\, u\, ]\, $, where $\, r\, $ is taken from the copy of 
$\, {\overline {R^n_{N , F}}}^{ _{\sps U^n_F}}\, $ and $\, u \in \widetilde {U^n}\, $. By the Jacobi identity the element 
$\, [\, c\, ,\, u_x\, ]\, $ is again in the ($\widetilde {U^n}$-normal) copy of 
$\, {\overline {R^n_{N , F}}}^{ _{\sps U^n_F}}\, $ giving a contradiction. This argument settles the general case modulo the Special Deconstruction Theorem which will have to wait until section 6. It remains to construct the almost canonical map 
$\, {\alpha }_n :\, U^{n+1}_{N , F} \rightarrow {\overline {R^n_{N , F}}}^{ _{\sps U^n_F}}\, $. One first divides by the subgroup $\, [\, J^{n+1}_{N , F}\, ,\, U^{n+1}_F\, ]\, +\, [\, U^{n+1}_{N , F}\, ,\, J^{n+1}_F\, ]\, $. Then the quotient contains the kernel of $\, {\overline {R^n_{N , F}}}^{ _{\sps U^n_F}}\, $ as a subgroup, moreover the commutators 
$\, [\, {\overline {R^n_{N , F}}}^{ _{\sps U^n_F}}\, ,\, U^n_F\, ]\, $ are uniquely determined as the image of the commutator subgroup $\, [\, U^{n+1}_{N , F}\, ,\, U^{n+1}_F\, ]\,$ and contain the kernel of 
$\, {\overline {R^n_{N , F}}}^{ _{\sps U^n_F}}\, $. By the universal property of the pair 
$\, ( U^{n+1}_{N , F}\, ,\, U^{n+1}_F )\, $ one gets an equivariant lift of the natural evaluation map 
$\, U^{n+1}_{N , F}\,\longrightarrow\, R^n_{N , F}\, $ to $\, {\overline {R^n_{N , F}}}^{ _{\sps U^n_F}}\, $ which must factor over the canonical $U^n_F$-central extension and is canonical on commutators. In particular it gives the identity map on the kernel of $\, {\overline {R^n_{N , F}}}^{ _{\sps U^n_F}}\, $ \qed
\bigskip\bigskip\bigskip
\par\noindent
\hfil\hfil \Large{\bf{Section III -- Strict splitting}} \hfil 
\bigskip\bigskip
\par\noindent

Before coming to the main result we will fix some notations. Let $\, F\, $ be  a (not necessarily perfect) group and $\, R\, $ a normal subgroup of $\, U_F\, $ contained in $\, J_F\, $. One has the following commutative diagram with exact rows and columns
\bigskip
$$ \vbox{\halign{ #&#&#&#&#\cr
$J_{R\, ,\, U_F}$ & $\largerightarrow$ & $U_{R\, ,\, U_F}$ & $\largerightarrow$ & $\, R$ \cr
\hfil $\Biggm\downarrow$\hfil && \hfil $\Biggm\downarrow$\hfil && $\,\,\Biggm\downarrow$ \cr
\hfil $J_{U_F}$\hfil & $\largerightarrow$ & \hfil $U^2_F$\hfil & $\largerightarrow$ & $\,\, U_F$ \cr
\hfil $\Biggm\downarrow$\hfil && \hfil $\Biggm\downarrow$\hfil && $\,\,\Biggm\downarrow$ \cr
$J_{U_F\, /\, R}$ & $\largerightarrow$ & $U_{U_F\, /\, R}$ & $\largerightarrow$ & $U_F\, /\, R$ \cr}}\> .
\leqno\raise48pt\hbox{\, (3)} $$
\par\bigskip\noindent
To save notation we will write 
$\, U = U_{R\, ,\, U_F}\, ,\, J = J_{R\, ,\, U_F}\, ,\, V = J_{U_F}\, ,\, P = J_{U_F\, /\, R}\, $. Fix a section 
$\, \{\, \overline x\, \}\, $ to the quotient map $\, U_F\,\rightarrow\, U_F\, /\, R\, $ such that each basis element $\, u_x\, $ is contained in the image set $\, \{\, \overline x\, \}\, $ for $\, x \in F\, $. Such a section will be referred to as a {\fndef good section} in the following. Consider the natural map 
$\, U_F\,\rightarrow\, U_{U_F\, /\, R}\, $ given by composing the canonical lift 
$\, U_F\,\rightarrow U^2_F\, $ with the quotient map $\, U^2_F\,\rightarrow\, U_{U_F\, /\, R}\, $. This map has a left inverse over $\, U_F / R\, $ given by composing the lift 
$\, U_{U_F\, /\, R}\,\rightarrow\, U^2_F\, $ induced by the good section $\, \{\, \overline x\, \}\, $ with the canonical projection onto $\, U_F\, $. Let $\, P_1\, $ be the intersection of $\, P\, $ with the image of 
$\, U_F\, $ by the map above and $\, P_0\, $ the kernel of the left inverse to this map, so that 
$\, P\, =\, P_0\,\rtimes\, P_1\, $ and $\, P_1\,\simeq\, R\, $.Then $\, P_0\, $ lifts naturally to $\, V\, =\, J_F\, $ where it is normal for the image of the canonical lift $\, U_F\,\rightarrow\, U^2_F\, $ and in fact has a core for the $\,\{\, u_{u_x}\,\}\, $ given by
$\,\{\, u^{\pm }_{u_{1 , \overline x}}\cdots u^{\pm }_{u_{m , \overline x}}\, u^{-1}_{\overline x}\,\}\, $ where 
$\, {\overline x}\, =\, u^{\pm }_{1 , \overline x}\cdots u^{\pm }_{m , \overline x}\, $ and 
$\, {\overline x}\,\notin\, \{\, u_x\,\}\, $ is understood. Let $\, C_U\, $ be the standard core of $\, U\, $ for the 
$\, u_{\overline x}\, $ and $\, C_V\, $ the core of $\, V\, $ for the $\,\{\, u_{u_x}\,\}\, $ given by
\smallskip
$$ \left\{\, u^{\pm}_{u_{x_1}}\cdots u^{\pm }_{u_{x_n}}\, u^{-1}_{u^{\pm }_{x_1}\cdots u^{\pm }_{x_n}}\,
\right\}\>  .  $$
\par\smallskip\noindent
As we have seen $\, P_0 \subset V\, $ has a core $\, C_{P_0}\subset C_V\, $ for the 
$\, \{\, u_{u_x}\,\}\, $. Project $\, C_V\, $ to $\, P\, $, then each basis element $\, v \in C_V\, $ is mapped to a product $\, c_v\, b_v\, $ with $\, b_v \in C_{P_0}\, $ and $\, c_v \in P_1\, $. Then 
$\, \{\, v\, b^{-1}_v\, \}\, $ is a core for the $\,\{\, u_{u_x}\,\}\, $ of some subgroup 
$\,\widetilde L_1\subset L_1\, $ where $\, L_1\, $ denotes the preimage of $\, P_1\, $ in $\, V\, $. Take the full basis $\,\widetilde Y\, $ of $\,\widetilde L_1\, $ determined by this core. Let 
$\, \overline L_1\, =\, V\, \widetilde L_1\, V^{-1}\, $ be the normal closure of $\,\widetilde L_1\, $ in 
$\, V\, $. A basis $\,\overline Y\, $ of $\,\overline L_1\, $ is obtained by taking adjoints of $\,\widetilde Y\, $ with elements of $\, P_0\, $. Since the basis $\, B\, $ of $\, P_0\, $ and $\,\widetilde Y\, $ both have a core for the $\,\{\, u_{u_x}\,\}\, $, so has $\,\overline Y\, $. Let us take a look at some reflection maps of the diagram (3). Define a map from the standard basis of $\, U\, $ to the basis of $\,\overline L_1\, $ by the following rule: for $\, u_x\in U^2_F\, $ a basis element, let 
$\, x\, =\, u^{\pm }_{1 , x}\cdots u^{\pm }_{n , x}\, $ and 
$\,\overline x\, =\, u^{\pm }_{1 , \overline x}\cdots u^{\pm }_{m , \overline x }\, $. Then define
\smallskip 
$$ \rho ( u_x\, u^{-1}_{\overline x} )\, =\, u^{\pm }_{u_{1 , x}}\cdots u^{\pm }_{u_{n , x}}\, u^{-1}_x\,
u_{\overline x}\, u^{\mp }_{u_{m , \overline x}}\cdots u^{\mp }_{u_{1 , \overline x}}  $$
\par\smallskip\noindent
and extend this map equivariantly for the action of the $\, \{\, u_{\overline x}\,\}\, $, i.e. if 
$\, u\, =\, u^{\pm }_{\overline x_1}\cdots u^{\pm }_{\overline x_l}\, $ then 
$\, \rho ( u\, u_x\, u^{-1}_{\overline x}\, u^{-1})\, =\, u\, \rho ( u_x\, u^{-1}_{\overline x} )\, u^{-1}\, $. One checks that by this rule the basis of $\, U\, $ is mapped to $\,\overline Y\, $ in a one-to-one fashion, so that one gets an isomorphism $\, U\buildrel\rho\over\longrightarrow \overline L_1\, $ with inverse 
$\, \overline L_1\buildrel\lambda\over\longrightarrow U\, $ given by reflection of the basis. By construction $\,\rho\, $ and $\,\lambda\, $ extend to the whole of $\, U^2_F\, $, and the restriction of 
$\,\rho\, $ to the subgroup $\,\widetilde U\, $ of $\, U\, $ which is generated by the standard core 
$\, C_U\, $ and conjugates by arbitrary elements in the image of the canonical lift of $\, U_F\, $ to 
$\, U^2_F\, $ maps $\,\widetilde U\, $ onto $\,\widetilde L_1\, $. Note that $\,\rho\, $ is a lift of the map 
$\, R\rightarrow P_1\, $ as above in the sense that the composition 
$\, U\buildrel\rho\over\longrightarrow {\overline L}_1\longrightarrow P\longrightarrow P_1\, $ is equal to 
$\, U\longrightarrow R\longrightarrow P_1\, $, and that $\,\lambda\, $ is a lift of (the image of 
$\,\overline L_1\, $ in) $\, P \rightarrow R\, $. However, $\, \rho ( J ) \subsetneq J\, $ and 
$\, J \subsetneq \overline L_1\, $. To see this note that the basis of $\,\widetilde L_1\, $ can be completed to a full basis $\, X\, $ of $\, V\, $ by adding the basis $\, B\, $ of $\, P_0\, $ and since $\, \overline L_1\, $ is the normal closure of $\,\widetilde L_1\, $ in $\, V\, $ (it is in fact normal in $\, U^2_F\, $), there is a natural projection $\, {\pi }_1: V = \overline L_1 P_0\longrightarrow P_0\, $. If $\, J\, $ would be contained in $\, \overline L_1\, $ then the map $\, {\pi }_1\, $ would factor over the projection 
$\, \pi : V\longrightarrow P\, $. But the image of $\,\overline L_1\, $ under this map is the normal closure of $\, P_1\, $ in $\, P\, $ which has nontrivial intersection with $\, P_0\, $. However on 
$\, P_0 \subset V\, $ both maps agree so that some part of $\, P_0\, $ would be mapped into the image of $\,\overline L_1\, $ under $\,\pi \, $ and hence must lie in the kernel of $\, {\pi }_1\, $, a contradiction! We will also consider $\, L_1\, =\, J\, +\, \widetilde L_1\, $.
\par\smallskip\noindent
A good section $\, \{ \overline x \} : U_F / R \nearrow U_F\, $ will be called a 
{\fndef very good section } iff the following holds: if 
$\, \overline x = u^{\pm }_{x_1}\cdots u^{\pm }_{x_n}\, $ and if $\, y\, $ is the image of 
$\, u^{\pm }_{x_k}\cdots u^{\pm }_{x_{k+l}}\, $ with $\, k+l \leq n\, $ in $\, U_F / R\, $ then 
$\, \overline y = u^{\pm }_{x_k}\cdots u^{\pm }_{x_{k+l}}\, $. 
To construct a very good section choose a well order on the elements of $\, F\, $ inducing a well order on the elements of $\, U_F\, $ by lexicografic ordering (and defining the inverse of a basis element 
$\, u^{-1}_x\, $ to be larger than any basis element $\, u_y\, $). Then choosing the minimal element among all possible lifts of some element in the quotient $\, U_F / R\, $ defines a very good section which in addition comes with a natural well order on its elements.
In the following, if $\, A\, $ is a subset of some group we let $\,\langle A\rangle\, $ denote the subgroup generated by $\, A\, $. In particular if $\, C_U\, $ is the core of $\, U\, $, let 
$\,\mathcal{C}_U = \langle C_U\rangle\, $. Then $\, {\mathcal C}_U\, $ is naturally defined, independent of the chosen section $\, U_F / R \nearrow U_F\, $. We make the following
\par\bigskip\noindent
{\bf Definition 7.} Let $\,\widetilde C_U\, $ be a basis of $\, {\mathcal C}_U\, $. A splitting $\, R\buildrel s\over\longrightarrow U\, $ will be called {\fndef semicanonical} (with respect to $\,\widetilde C_U\, $) if there exists a basis $\,\{ f_k \}\, $ of $\, R\, $ such that $\,\{ s( f_k ) \} \subset \{\, u \widetilde C_U u^{-1}\,\}\, $ with $\, u \in \langle\,\{\, u_{u_x}\,\}\,\rangle\, $. It will be called a {\fndef strict splitting} if there exists a subset 
$\,\{ t_j \} \subset \widetilde C_U\, $ and a fixed section $\,\{ \overline x \} \subseteq U_F\, $ containing the identity such that 
\smallskip
$$ s( f_k ) =  \overline u\, t_j\, \overline u^{-1}  $$
\par\medskip\noindent 
where $\, \overline u\, $ is the image of $\, \overline x\, $ by the canonical lift 
$\, U_F \rightarrow U^2_F\, $, and if secondly the section $\, \{ \overline x \}\, $ has the following convergence property: if $\, \overline x_1\, ,\, \overline x_2\, $ are two given elements of the section consider the difference $\, \overline x_1^{-1}\, \overline x_2\, ( \overline x_1^{-1}\, \overline x_2 )^{-1}\, $, where 
$\, ( x )\, $ denotes the element of $\, \{ \overline x \}\, $ corresponding to some fixed $\, x \in U_F\, $. This difference may be written as a product
\smallskip
$$ \bigl( \overline x_{j_1}\, e_{j_1}\, \overline x_{j_1}^{-1} \bigr)\cdots 
\bigl( \overline x_{j_m}\, e_{j_m}\, \overline x_{j_m}^{-1} \bigr) $$
\par\medskip\noindent
where $\, e_j \subseteq \{ f_k \}\, $ is the basis element corresponding to $\, t_j\, $ and each bracket 
$\, \overline x_{j_l}\, e_{j_l}\, \overline x _{j_l}^{-1}\, $ corresponds to a basis element of $\, R\, $. Also we put $\,\overline x_{j_{m+1}} := ( \overline x_1^{-1}\, \overline x_2 )\, $. Then iterating this process with all possible pairs $\, ( \overline x_{j_k}\, ,\, \overline x_{j_l} )\, $ for $\, k \neq l\, $ one arrives at the identity in each path after a finite number of steps (independent of the path). 
\par\bigskip\noindent
We remark that a splitting satisfying only the first condition of a strict splitting is semicanonical and normal for the action of $\, U_F\, $ modulo $\, [\, V\, ,\, U\, ]\, $. In particular $\, K^J_2 ( R\, ,\, U_F ) = 0\, $. The second condition is certainly most difficult to verify in practice. One way to do this is by introducing a level function
\smallskip
$$ l : \{ \overline x \} \longrightarrow \mathbb N $$
\par\medskip\noindent
on the section (which could be for example the length of an element by a given basis) such that the only element of level zero is the identity element and such that the level of the elements 
$\, ( \overline x_1^{-1}\, \overline x_2 )\, $ decreases in each path by the process as above. To see that strict splittings are not just theoretical objects we provide the following example which is basic in many technical applications.
\par\medskip\noindent
{\it Example.}\quad Let $\, F = D \times E\, $ be a direct product and consider the intersection 
$\, R = U_{D , F} \cap U_{E , F} \subseteq J_F\, $. We may change the standard basis of $\, U_F\, $ to the union of the standard bases of $\, U_D\, $ and $\, U_E\, $ and a core 
$\, \{ u_{c\, x}\, u^{-1}_x\, u^{-1}_c\,\vert\, c \in D\, ,\, x \in E \}\, $ contained in $\, R\, $. Then $\, R\, $ has a basis of the form
\smallskip
$$ \bigl\{\> u^{\pm }_{d_1}\cdots u^{\pm }_{d_m}\, u^{\pm }_{y_1}\cdots u^{\pm }_{y_n}\, 
\bigl( u_{c x}\, u^{-1}_x\, u^{-1}_c \bigr)\,
u^{\mp }_{y_n}\cdots u^{\mp }_{y_1}\, u^{\mp }_{d_m}\cdots u^{\mp }_{d_1}\> 
\bigr\} $$
$$\cup\quad 
\bigl\{\> u^{\pm }_{d_1}\cdots u^{\pm }_{d_m}\, u^{\pm }_{y_1}\cdots u^{\pm }_{y_n}\, 
\bigl[\, u_x\, ,\, u_c\,\bigr]\, u^{\mp }_{y_n}\cdots u^{\mp }_{y_1}\, u^{\mp }_{d_m}\cdots u^{\mp }_{d_1}\> 
\bigr\}\quad $$
\par\medskip\noindent
so that the elements $\, e_j\, $ correspond to brackets $\, [\, u_x\, ,\, u_c\, ]\, $ or to brackets 
$\, ( u_{c x}\, u^{-1}_x\, u^{-1}_c )\, $ with $\, c \in D\, ,\, x \in E\, $ and the section is given by partially ordererd expressions of the form
$\, u^{\pm }_{d_1}\cdots u^{\pm }_{d_m}\, u^{\pm }_{y_1}\cdots u^{\pm }_{y_n}\, $ where 
$\, d_1, \cdots , d_m \in D\, $ and $\, y_1 , \cdots , y_n \in E\, $ is understood. It is clear that the first condition of a strict splitting is satisfied. The level function is given by the length of an element with respect to the basis and one has to check that this function decreases strictly monotonous by the process indicated in Definition 7. This is left as an exercise to the reader.
\par\medskip\noindent
For $\, R\, $ as in diagram (3) define $\, R^n\, $ inductively by $\, R^1 = R\, $ and 
$\, R^n\, =\, J_{R^{n-1}\, ,\, U^{n-1}_F}\, $, let $\, U^n\, =\, U_{R^n\, ,\, U^n_F}\, $ and similarly for the other groups appearing in (3). 
\par\bigskip\noindent
{\bf Theorem 3.}\quad Assume that there exists a strict splitting $\, R\buildrel s\over\longrightarrow U\, $. Then for every $\, n\geq 1\, $ there is a strict splitting $\, R^n\buildrel s_n\over\longrightarrow U^n\, $.
\par\bigskip\noindent
{\it Proof.}\quad The proof is by induction on $\, n\, $, so we can assume that $\, n = 1\, $ and must show that given a strict splitting $\, R\buildrel s\over\longrightarrow U\, $ one gets a strict splitting 
$\, J \longrightarrow U^2\, $. With notation as in Definition 7 there is, for every basis element $\, f_k\, $ of $\, R\, $, a unique basis element 
$\, t_j \in \widetilde C_U\, $ such that $\, s( f_k ) = u\, t_j\, u^{-1}\, $ for some $\, u\, $ (in the following 
$\, u\, $ always denotes a typical element of 
$\,\langle\{\, u_{u_x}\,\}\rangle\, $) and $\, t_j = s ( e_j )\, $. Let $\,\widetilde C_U\, =\,\{ z_i \} \cup \{ t_j \}\, $. Then a core of 
$\,\widetilde U\, $ for the $\,\{ u_{u_x} \}\, $ is given by the set 
$\,\{\, z_i\, s( \overline z_i)^{-1}\, ,\, t_j\, \}\, =\, \{ x_k \}\, $ with $\,\overline z_i\, $ the image of $\, z_i\, $ in 
$\, R\, $ (not to be mistaken with the images $\,\overline x\, $ of the section 
$\, U_F / R \nearrow U_F\, $). One easily checks that a basis of the normal subgroup 
$\, [\, P_0\, ,\, U\, ]\, =\, [\, P_0\, ,\,\widetilde U\, ]\, $ is given by the set
\smallskip
$$ \bigl\{\, y^{\pm }_{m_1}\cdots y^{\pm }_{m_n}\, b^{\pm }_{l_1}\cdots b^{\pm }_{l_p}\,
\bigl[\, y_{m_0}\, ,\, b_{l_0}\,\bigr]\, 
b^{\mp }_{l_p}\cdots b^{\mp }_{l_1}\, y^{\mp }_{m_n}\cdots y^{\mp }_{m_1}\,\bigr\}  $$
\par\smallskip\noindent
with $\, y_{m_j} \in \{ u\, x_k\, u^{-1} \}\, $ and $\, b_{l_j}\in B\, $, the basis of $\, P_0\, $. Write 
$\, [\, P_0\, ,\, U\, ]\, $ as a semidirect product of two parts 
$\, [\, P_0\, ,\, U\, ]_0 \rtimes [\, P_0\, ,\, U\, ]_1\, $ where $\, [\, P_0\, ,\, U\, ]_0\, $ is generated normally by that part of the basis where at least one of the $\, y_{m_j}\, $ is of the form 
$\, u\, z_i\, s(\overline z_i )^{-1}\, u^{-1}\, $ and $\, [\, P_0\, ,\, U\, ]_1\, $ is generated by the remaining basis elements. It is clear that $\, [\, P_0\, ,\, U\, ]_0\, $ has a core for the $\,\{ u_{u_x} \}\, $, for 
$\,\widetilde C_U\, $ and for $\, C_{P_0}\, $ modulo $\, [\, J\, ,\, J\, ]\, $. A lift of $\, [\, P_0\, ,\, U\, ]_0\, $ which is normal modulo $\, [\, V^2\, ,\, U^2\, ]\, $ is then determined by sending the commutator 
$\, [\, y_{m_0}\, ,\, b_{l_0}\, ]\, $ to 
$\, [\, \iota ( u )\, u_{x_k}\, \iota ( u )^{-1}\, ,\, u_{b_{l_0}}\, ]\, $ if $\, b_{l_0}\in C_{P_0}\, $ and 
$\, y_{m_0}\, =\, u\, z_i\, s( \overline z_i )^{-1}\, u^{-1}\, $, and to 
$\, u_{[\, y_{m_0}\, ,\, b_{l_0}\, ]}\, $ if $\, b_{l_0} \in C_{P_0}\, $ and 
$\, y_{m_0} \in \{ u\, t_j\, u^{-1} \}\, $ and completing normally for the lift 
\smallskip
$$ \iota : U^2_F \rightarrow U^3_F\, ,\quad \iota ( u_{u_x} )\, =\, u_{u_{u_x}}\, ,\quad 
\iota ( x_k )\, =\, u_{x_k}\, ,\quad \iota ( b )\, =\, u_b\, ,\, b \in C_{P_0}\>  . $$
\par\smallskip\noindent
Since $\, J\, =\, [\, P_0\, ,\, U\, ]\,\rtimes\,\widetilde J\, $ with $\,\widetilde J\, =\, J\cap \widetilde U\, $ one also has to construct a normal lift of $\,\widetilde J\, $ to $\, U^2\, $. We use a similar decomposition 
$\, \widetilde J\, =\, \widetilde J_0\,\rtimes\,\widetilde J_1\, $ into a semidirect product. Let 
$\,\widetilde J_0\, $ be the normal subgroup of $\,\widetilde U\, $ generated by the elements
\smallskip
$$  \left\{\, u\> \widetilde C_{U , 0}\> u^{-1}\,\right\}\quad ,\quad 
\widetilde C_{U , 0}\, =\, \left\{\, z_i\, s(\overline z_i )^{-1}\,\right\}\>  . $$
\par\smallskip\noindent
This subset defines a core of $\,\widetilde J_0\, $ for the $\,\{ s( f_j ) \}\, $, so that a lift which is normal for 
$\,\widetilde U\, +\, \langle\{ u_{u_x} \}\rangle\, $ (modulo $\, [\, V^2\, ,\, U^2\, ]\, $) is given by restriction of the map $\,\iota\, $ defined above. Clearly, it is compatible with the lift of $\, [\, P_0\, ,\, U\, ]_0\, $, hence the combined lift will be normal also for $\, C_{P_0}\, $. If $\,\widetilde J_1\, $ is the quotient of 
$\,\widetilde J\, $ by $\,\widetilde J_0\, $, let $\, J_1\, =\, [\, P_0\, ,\, U\, ]_1\> \widetilde J_1\, $, so that 
$\, J_1\, $ is the quotient of $\, J\, $ by $\, J_0\, =\, [\, P_0\, ,\, U\, ]_0\> \widetilde J_0\, $. Then 
$\, \widetilde J_1\, $ is naturally isomorphic to the subgroup 
$\,\langle\{ u\, s( f_j )\, u^{-1} \}\rangle \cap \widetilde J\, $ and is completely contained in 
$\, [\, V\, ,\, U\, ]\, $ by assumption on $\, s\, $. Let $\,\{ g_k\}\,$ be the basis of $\, P_1\, $ corresponding to 
$\,\{ f_k \}\, $ and $\,\widetilde s( g_k )\, =\, u_k\, s( f_k )^{-1}\in L_1\, $ where $\, u_k\, $ is the image of 
$\, f_k\, $ by the canonical lift $\, U_F \rightarrow U^2_F\, $, so that $\, \{ \widetilde s( g_k ) \}\, $ defines a splitting $\, \widetilde s : P_1 \rightarrow L_1\, $. By assumption on $\, s\, $ every element 
$\, u\, s( f_k )\, u^{-1}\, $ can be written as 
$\,  v\, s( f_l )\, v^{-1}\, $ with $\, v \in \langle \{ u_k \}\rangle\, $. Since 
$\, u_k = \widetilde s ( g_k )\, s ( f_k )\, $ one gets that $\,\widetilde J_1\, $ is generated (freely) by elements of the form
\smallskip
$$ \bigl\{\> w\, y
\, \bigl[\, \widetilde s ( g_{l_0} )\, ,\, s ( f_{k_0} )\,\bigr]\, y^{-1} \, w^{-1} \> \bigr\} $$
\par\medskip\noindent
with $\, y = s ( f_{j_1} )^{\pm }\cdots s ( f_{j_n} )^{\pm }\, ,\, 
w = \widetilde s ( g_{i_1} )^{\pm }\cdots \widetilde s ( g_{i_m} )^{\pm }\, $ or equivalently by
\smallskip
$$ \bigl\{\> v\, y\, \bigl[\, \widetilde s ( g_{l_0} )\, ,\, s ( f_{k_0} )\, \bigr]\, y^{-1}\, v^{-1}\>\bigr\} $$
\par\medskip\noindent
with $\, v \in \langle \{ u_k \} \rangle\, $ and $\, y\, $ as above. We want our basis to have a core for the elements $\, u \in \langle \{ u_{u_x} \} \rangle\, $ and not just for $\, v \in \langle \{ u_k \} \rangle\, $. To this end we propose that the following set  
\smallskip
$$ A\> =\> \bigl\{\> u\, s ( f_{k_1} )^{\pm }\cdots s ( f_{k_m} )^{\pm }\, 
\bigl[\, t_j\, ,\, \widetilde s ( g_l )\,\bigr]\, 
s ( f_{k_m} )^{\mp }\cdots s ( f_{k_1} )^{\mp }\, u^{-1} \>\bigr\} $$
\par\medskip\noindent
with arbitrary coefficients $\, s ( f_{j_1} )^{\pm }\cdots s ( f_{j_m} )^{\pm }\, u\, $, arbitrary 
$\, t_j\, ,\, \widetilde s ( g_l )\, $ is linear independent and generates 
$\, \widetilde J_1\, $. We first show that it generates $\,\widetilde J_1\, $. Consider basis elements of the form $\, v\, [\, s( f_k )\, ,\, \widetilde s( g_l )\, ]\, v^{-1}\, $ with $\, v \in \langle \{ u_k \}\rangle\, $ which may be written in the form
\smallskip
$$ v\, \overline u_1\, \bigl[\, t_j\, ,\, u^{\pm }_{k_1}\cdots u^{\pm }_{k_m}\, \widetilde s ( g_{l'} )\, 
u^{\mp }_{k_m}\cdots u^{\mp }_{k_1}\, \bigr]\, \overline u_1^{-1}\, v^{-1} $$
\par\medskip\noindent
with $\, s( f_k ) = \overline u_1\, t_j\, \overline u^{-1}_1\, ,\, 
\widetilde s( g_l ) = \overline u_2\, \widetilde t_m\, \overline u^{-1}_2\, ,\, 
\widetilde s( g_{l'} ) = ( \overline u^{-1}\, \overline u_2 )\, \widetilde t_m\, 
( \overline u^{-1}_1\, \overline u_2 )^{-1}\, $ and 
$\, u^{\pm }_{k_1}\cdots u^{\pm }_{k_m} = 
\overline u^{-1}_1\, \overline u_2\, ( \overline u^{-1}_1\, \overline u_2 )^{-1}\, $. It is easily checked that the bracket may be extended into certain conjugates by the elements $\, s ( f_{k_j} )^{\pm }\, $ and 
$\, u_{k_j}\, $  of basic brackets 
$\, [\, \overline u_{j_p}\, t_{j_p}\, \overline u^{-1}_{j_p}\, ,\, 
\overline u_{j_q}\,\widetilde t_{j_q}\, \overline u^{-1}_{j_q}\, ]\, $ such that iterating this process (i.e. extracting the coefficient $\, \overline u_{j_p}\, $ from the inner bracket each time one arrives at an expression in the elements of $\, A\, $ after finitely many steps. The argument for an arbitrary basis element with nontrivial $\, s( R )$-coefficient follows by induction on the length of the coefficient in the basis $\,\{ s( f_k ) \}\, $. Assume the result is true for all coefficients of length less or equal to $\, n\, $. Then to prove the result for length $\, n+1\, $ it suffices that any expression 
$\, s( f_k )\, u\, $ is equivalent to an expression $\, u'\, y\, $ with $\, y \in s( R )\, $ modulo the range of 
$\, A\, $. Clearly one gets $\, s( f_k )\, u = u'\, s( f_{k'} )\, v\, $ so it suffices to consider $\, u = v\, $. Then one may again use induction on the length of $\, v\, $ in the basis $\,\{ u_k \}\, $, i.e. one can assume that 
$\, v = u^{\pm }_l\, $. One gets 
\smallskip
$$ u_l\, s( f_k )\, u^{-1}_l\> =\>
\bigl[\, \widetilde s( g_l )\, ,\, s( f_l )\, \bigr] \cdot 
\bigl( s( f_l )\, \bigl[\, \widetilde s( g_l )\, ,\, s( f_k )\, \bigr]\, s( f_l )^{-1}\bigr) \cdot $$ 
$$ \bigl( s( f_l )\, s( f_k )\, s( f_l )^{-1}\, \bigl[\, s( f_l )\, ,\, \widetilde s( g_l )\, \bigr]\, 
s( f_l )\, s( f_k )^{-1}\, s( f_l )^{-1} \bigr) $$
$$ s( f_l )\, s( f_k )\, s( f_l )^{-1} \> . $$  
\par\medskip\noindent
The first bracket is in $\, A\, $ and the second has coefficient length equal to 1 so in case that 
$\, n \geq 1\, $ it is in the range of $\, A\, $ by induction assumption, otherwise the level of $\, s( f_l )\, $ (i.e. the level of the corresponding coefficient $\, \overline u\, $) is smaller than the level of 
$\, s( f_k )\, $ if a level function with strictly monotonous decreasing behaviour by the process of Definition 7 is defined, in any case the process converges so that the elements 
$\, s( f_l )\, [\, \widetilde s( g_l )\, ,\, s( f_k )\, ]\, s( f_l )^{-1}\, $ can be assumed to lie in the range of $\, A\, $. As for the third bracket which has coefficient length 3, one notes that only the coefficient $\, s( f_k )\, $ counts and then the same argument as above applies. This proves that $\, A\, $ generates 
$\,\widetilde J_1\, $. Next we show that $\, A\, $ is linear independent. For this purpose we may consider the equivalent set which is the union of the two sets
\smallskip
$$ \bigl\{\> u\, s( f_{k_1} )^{\pm }\cdots s( f_{k_m} )^{\pm }\, 
\bigl[\, s( f_k )^{\pm }\, ,\, \bigl[\, t_j\, ,\, \widetilde s( g_l )\,\bigr]\,\bigr]\, 
s( f_{k_m} )^{\mp }\cdots s( f_{k_1} )^{\pm }\, u^{-1}\>\bigr\} $$
$$ \quad\cup\quad \bigl\{\> u\, \bigl[\, t_j\, ,\, \widetilde s( g_l )\,\bigr]\, u^{-1}\>\bigr\}\> . $$
\par\medskip\noindent 
Considering the expansion in the subbasis $\, \{ u \} \cup \{ t_j \}\, $ define the meridian of an element of the first kind to be the basis element $\, t^{\mp }_q\, $ corresponding to $\, s( f_k )^{\mp }\, $ appearing in the inner bracket $\, [\, s( f_k )^{\pm }\, ,\, [\, t_j\, ,\, \widetilde s( g_l )\, ]\, ]\, $. It is clear that the meridian of an element of the first kind cannot be cancelled by any letter of an adjacent element of either first or second kind unless this element is precisely the inverse of the given one. Also one defines the meridian for elements of the second kind to be the letter $\, t_j^{-1}\, $ in the inner bracket 
$\, [\, t_j\, ,\, \widetilde s( g_l )\, ]\, $. One readily checks that if the meridian of an element of the second kind is cancelled by any letter of an adjacent element which is also of the second kind then this element must be its inverse. In particular the first and second set taken separately are both linear independent.
To show that the two sets are also mutually independent we will change the second set once more to the union of the two sets
\smallskip
$$ \bigl\{\> u\, \bigl[\, u_k\, ,\, \bigl[\, t_j\, ,\, \widetilde s( g_l )\,\bigr]\,\bigr]\, u^{-1}\>\bigr\} 
\quad\cup\quad \bigl\{\> \overline u\, \bigl[\, t_j\, ,\, \widetilde s( g_l )\,\bigr]\, \overline u^{-1}\>\bigr\} $$
\par\medskip\noindent
where in the first set the coefficient $\, u\, $ is arbitrary whereas the coefficient $\, \overline u\, $ in the second set is a canonical lift of an element of the section $\,\{ \overline x \}\, $ and $\, u_k\, $ is the canonical lift of the basis element $\, f_k\, $. Considering the union $\, B\, $ of all three sets we first show that the first two sets are independent. Define the left (resp. right) side of an element of the first or second kind to be the corresponding $\, u\cdot s( R )$-conjugate of the expression 
$\, s( f_k )\, [\, t_j\, ,\, \widetilde s( g_l )\, ]\, s( f_k )^{-1}\, $ for elements of the first kind and 
$\, u_k\, [\, t_j\, ,\, \widetilde s( g_l )\, ]\, u^{-1}_k\, $ for elements of second kind (resp. of
$\, [\, \widetilde s ( g_l )\, ,\, t_j\, ]\, $) and
assume given a minimal nontrivial relation in the generators of the first two sets which must contain elements of either kind since both sets taken separately are independent. Consider the relation as a word in the letters $\,\{ u\, t_j\, u^{-1} \}\, $. Then the letter cancelling against the first letter of the first element must terminate some (left or right) part of some other element occurring in the relation, since otherwise the remaining part of the relation to the right of the cancelling letter wouldn't be contained in $\, [\, V\, ,\, U\, ]\, $ which it must being trivial. Call the part up to the letter cancelling the first letter of the first element a compartment. By rotating the elements of the relation if necessary we may assume this compartment to be minimal, i.e. all the elements involved to be mutually connected by cancelling letters. 
Define the left (resp. right) meridian of an element of second kind to be the conjugate of $\, t^{-1}_j\, $ (resp. $\, t_j\, $) appearing in the left part 
$\, u\, u_k\, [\, t_j\, ,\, \widetilde s( g_l )\, ]\, u^{-1}_k\, u^{-1}\, $ (resp. the right part 
$\, u\, [\, \widetilde s( g_l )\, ,\, t_j\, ]\, u^{-1}\, $) of the bracket. In case that $\, t_j = s( f_l )\, $ the inner bracket reduces to $\, t_j\cdot ( u_j\, t^{-1}_j\, u^{-1}_j )\, $ of two letters instead of an expression with four letters so that in this case one defines the second letter to be the left meridian. Then it is clear that the right meridian of such an element cannot be cancelled by a letter in an adjacent element to the left (of first or second kind) unless it is the inverse of the given element, and the same holds for the left meridian with respect to an adjacent element on the right side 
(left and right have to be exchanged considering the inverse of an element). Also if the left meridian is cancelled by a letter in an adjacent element on the left then this element has the same inner bracket 
$\, [\, t_j\, ,\, \widetilde s( g_l )\, ]\, $ and the coefficient $\, u\, $ differs at most by 
$\, u^{\pm }_k\, u^{\mp }_{k'}\, $ and a corresponding statement holds for the right meridian. For each element in the compartment whose (left or right) meridian is not cancelled by a letter from an adjacent element there is at least one element of second kind inbetween whose left or right meridian cancelles by an adjacent element so that by the argument above the whole left or right part cancelles against the right or left part of the adjacent element. We make the following induction assumption: for given $\, n\, $ if a left or right meridian cancelles against a letter in an element with less than $\, n\, $ elements inbetween, then this element has the same inner bracket $\, [\, t_j\, ,\, \widetilde s( g_l )\, ]\, $ and the cancelling extends to the whole left or right part of the corresponding element. Obviously this is true for 
$\, n = 1\, $. Suppose it holds for $\, n\, $. Then if the meridian cancelles by a letter of an element with 
$\, n\, $  elements inbetween all the meridians of the intermediate elements must cancel from a shorter distance so by induction assumption they all have the same inner brackets 
$\, [\, t_j\, ,\, \widetilde s( g_l )\, ]\, $ since they must all be mutually connected so that the result also holds for $\, n+1\, $ and what is more one finds that the number of left parts must equal the number of right parts so the compartment gives a subrelation in the elements of the first two sets and from minimality this is the whole relation so each right part of any element, except the last one which cancelles against the left part of the first element, cancelles against the left part of the adjacent element on the right. A simple argument then shows that the relation must be elementary, giving a contradiction.  
Now include the third set into consideration. Each element of the first two sets lies in the normal closure of the range of the subbasis
\smallskip
$$ \bigl\{\>  v\, y\, \bigl[\, s( f_m )^{\pm }\, ,\, \bigl[\, s( f_k )\, ,\, \widetilde s( g_l )\,\bigr]\,\bigr]\, y^{-1}\, v^{-1}
\> \bigr\} $$
$$\quad\cup\quad \bigl\{\> v\, \bigl[\, u^{\pm }_m\, ,\, \bigl[\, s( f_k )\, ,\, \widetilde  s( g_l )\, \bigr]\,\bigr]\, 
v^{-1}\> \bigr\}\qquad $$
\par\medskip\noindent
and modulo this subbasis each element of the third set is congruent to a unique basis element 
\smallskip
$$ \bigl\{\> \bigl[\, s ( f_k )\, ,\, \widetilde s( g_l )\,\bigr]\>\bigr\}\> $$
\par\medskip\noindent
so that projecting a relation in the elements of $\, B\, $ to the subspace generated by the 
elementary brackets as above one gets an elementary relation in the images of the third set.
Let $\, C\, $ be the set of all conjugates of the first two sets of $\, B\, $ by arbitrary elements in the range of the third set. Then the relation may be considered as a relation in the elements of $\, C\, $ and triviality of the latter implies triviality of the original relation, thus independence of $\, B\, $. For elements of $\, C\, $ one extends the notion of right and left part of an element to be the corresponding conjugates of the central element in $\, B\, $. Also the left and right meridian of an element of $\, C\, $ is defined to be the left or right meridian of the central element in $\, B\, $ which if it cancels by some conjugacy coefficient of the element is to be replaced by the corresponding letter of the inverse of this coefficient appearing on the other side of the central element. Then the proof given above for the first two sets of $\, B\, $ extends with minor changes which are left to the reader to the set $\, C\, $ proving independence of $\, B\, $ and thus of $\, A\, $. Let 
$\,\overline A\, $ the set obtained from $\, A\, $ by allowing arbitrary 
$\, \overline u \in \langle\{ u_{\overline x} \}\rangle\, $ instead of $\, u\, $ in the formula for the generators. Then it is plain to see that together with the subbasis of $\, [\, P_0\, ,\, U\, ]\, $ consisting of
\smallskip
$$ \bigl\{ u\, b^{\pm }_{l_1}\negthinspace\negthinspace\cdots b^{\pm }_{l_p}\, 
s( f_{j_1} )^{\pm }\negthinspace\negthinspace\cdots s( f_{j_n} )^{\pm }\, 
\bigl[\, s ( f_{j_0} )\, ,\, b_{l_0}\, \bigr]\, 
s( f_{j_n} )^{\mp }\negthinspace\negthinspace\cdots s( f_{j_1} )^{\mp }\, 
b^{\mp }_{l_p}\negthinspace\negthinspace\cdots b^{\mp }_{l_1}\, u^{-1} \bigr\}   $$
\par\medskip\noindent
where $\, b_{l_j} \in B\, $ for $\, j = 1 ,\cdots , p\, $ and $\, b_{l_0} \in C_{P_0}\, ,\,\overline A\, $ freely generates $\, J_1\, =\, [\, P_0\, ,\, U\, ]_1\,\widetilde J_1\, $. Indeed, consider the set 
$\,\{\, \overline u\, t_j\, {\overline u}^{-1}\,\}\, $ with 
$\, t_j \in \widetilde C_U \cap \{ s( f_j ) \}\, ,\, \overline u \in \langle\{ u_{\overline x} \}\rangle\, $ arbitrary. Then $\, J_1\, =\, \langle\{ \overline u\, t_j\, {\overline u}^{-1} \}\rangle \cap J\, $. The set 
$\, \{ u\, t_j\, u^{-1} \}\, $ freely generates $\,\widetilde J_1\cdot s( R )\, $ (semidirect product) so that
\smallskip 
$$ J_1\> =\> \bigl[\, P_0\, ,\, {{\widetilde J}_1}\cdot s( R )\,\bigr] \cdot {\widetilde J}_1\>  $$
$$ \> =\> \bigl(\, \bigl(\, \bigl[\, P_0\, ,\, {\widetilde J}_1\,\bigr]\, {\widetilde J}_1\,\bigr)\,
\bigl[\, P_0\, ,\, s( R )\,\bigr]\, \bigl(\, \bigl[\, P_0\, ,\, {\widetilde J}_1\,\bigr]\, {\widetilde J}_1\,\bigr)\,\bigr)^{-1} \cdot \bigl[\, P_0\, ,\, {\widetilde J}_1\,\bigr] \cdot {\widetilde J}_1  $$
\par\smallskip\noindent
(semidirect product). $\,\overline A\, $ is a basis for 
$\, [\, P_0\, ,\, {\widetilde J}_1\, ]\, {\widetilde J}_1\> =\> P_0\,  {\widetilde J}_1\, P^{-1}_0\, $ whereas the subbasis of $\, [\, P_0\, ,\, U\, ]_1\, $ as above is elementwise congruent to a basis of 
$\, [\, P_0\, ,\, s( R )\, ]\, $ modulo $\,\langle\{ \overline A \}\rangle\, $ since for every 
$\, u \in \langle\{ u_{u_x} \}\rangle\, $ the set $\, \{ s( u\, f_j\, u^{-1} ) \}\, $ is a basis of $\, s( R )\, $ and 
$\, u\, s( f_j )\, u^{-1}\,\equiv\, s( u\, f_j\, u^{-1} )\, $ modulo $\,\langle A \rangle\, $.  Then one only has to check that the subbasis as above becomes a core for $\,\langle\{ \overline A \}\rangle\, $ on deleting the $\, u\, $ and letting $\, b_{l_0} \in B\, $ instead of $\, b_{l_0} \in C_{P_0}\, $, which is obvious from the general basis of $\, [\, P_0\, ,\, U\, ]_1\, $. The subbasis of $\, [\, P_0\, ,\, U\, ]_1\, $ has a core for the 
$\,\{ u_{\overline x} \}\, $ and for the $\, \{ s( f_j ) \}\, $ modulo $\, [\, J_1\, ,\, J_1\, ]\, $ so adding this set to the existing core of $\, J_0\, $ for $\, J_1\, $ given by the set 
\smallskip
$$  \bigl\{\> \overline u\, s( f_{j_1} )^{\pm }\cdots s( f_{j_n} )^{\pm }\, 
\,{\widetilde C}_{U , 0}\, s( f_{j_n} )^{\mp }\cdots s( f_{j_1} )^{\mp }\, {\overline u}^{-1}\>\bigr\} $$
\par\medskip\noindent
one gets a basis of $\, J\, $ which has a core for the fixed section 
\smallskip
$$ U^2_F / J \nearrow U^2_F :\quad u\, b\, r\, \mapsto \iota ( u b )\, s( r ) $$
where $\, \iota : U_{U_F / R} \rightarrow U^2_F\, $ is the lift defined by some good section as described above and $\, s :  R \rightarrow U\, $ is the given strict splitting consisting of the union of the three sets 
\smallskip
$$ \widetilde C_{U , 0}\quad\cup\quad \bigl\{\> \bigl[\, s( f_k )\, ,\, b_l\, \bigr]\,\vert\, b_l \in C_{P_0}\>\bigr\}
\quad\cup\quad \bigl\{\> \bigl[\, t_j\, ,\, \widetilde s( g_l )\,\bigr]\>\bigr\} $$
\par\medskip\noindent
One checks that the convergence property needed to construct a strict splitting 
$\, s^2:  J \rightarrow U^2\, $ is inherited from the corresponding convergence property one dimension below. For example given a level function for the section $\, \{ \overline x \}\, $ used in the presentation of the basis of $\, R\, $ which also serves to define a level on each basis element $\, f_k\, $ by considering the conjugacy coefficient of the section, and thus also on the elements $\, \{ s( f_k )^{\pm } \}\, $ as well as 
the elements $\, \{ u^{\pm }_k \}\, $ one may define a level function on the section 
$\, u\, b\, s( r )\, $ by adding up the levels of its pure components 
$\, u = \overline u\, v \in \langle\{ u_{u_x} \}\rangle\, ,\, b_l \in C_{P_0}\, $ and $\, s( r )\, $ where the level of an element $\, \overline u\, $ equals the level of the corresponding element of the section 
$\, \{ \overline x \}\, $, the level of $\, v \in \langle \{ u_k \}\rangle \, $ is the sum of the levels of the basis elements appearing in its expression, the same for the level of an element of $\, s( R )\, $ and each basis element $\, b^{\pm }_l \in C_{P_0}\, $ raises the level by 1. Then one checks that this level function has the corresponding strictly monotonous decreasing property by the process of Definition 7 \qed
\par\bigskip\noindent
{\bf Remark.} It would be very interesting to have a wider and more flexible notion of strict splitting since the conditions of Definition 7 are very restrictive. For example one may try to replace the property of 
$\, s\, $ being semicanonical in the definition of a strict splitting by the following weaker condition: there exists a subset $\, \{ t_l \} \subset \widetilde C_U\, $ such that $\, \langle\{ u\, s( f_j )\, u^{-1} \}\rangle\, =\, \langle\{ u\, t_l\, u^{-1} \}\rangle\, $ and for every $\, u\, $ and $\, s( f_k )\, $ there exist elements 
$\, v_{j_1} ,\cdots , v_{j_l}\, $ and $\, s( f_{j_1} ),\cdots , s( f_{j_l} )\, $ such that 
\smallskip
$$  u\, s( f_k )\, u^{-1}\, =\, ( v_{j_1}\, s( f_{j_1} )^{\pm }\, v^{-1}_{j_1} )\cdots 
( v_{j_l}\, s( f_{j_l} )^{\pm }\, v^{-1}_{j_l} )\> .   $$
\par\smallskip\noindent
One can show that such a decomposition is necessarily unique. In case of a semicanonical splitting one of course gets 
$\, l = 1 \, $ and the condition is then equivalent to saying that there is a basis $\,\{ f_j \}\, $ of $\, R\, $ such that for every basis element $\, f_k\, $ and each $\, x \in U_F\, $ there is  $\, v \in R\, $ such that 
$\, (v x)\, f_k\, (vx)^{-1}\, =\, f_l\, $ is again a basis element of $\, R\, $ (but not necessarily by a fixed section as in Definition 7). In the general case (without assuming that $\, s\, $ is semicanonical one only gets that there exists a basis $\, \{ g_k \}\, $ of $\, R\, $ such that for each basis element $\, g_k\, $ and each $\, x \in U_F\, $ there is another basis element $\, g_l\, $ congruent to $\, x\, g_k\, x^{-1}\, $ modulo 
$\, [\, R\, ,\, R\, ]\, $. Indeed, from the condition above one gets that 
$\, J \cap \langle \{ u\, t_j\, u^{-1} \}\rangle\, $ is contained in $\, [\, V\, ,\, U\, ]\, $ so dividing by the normalization of the elements $\, v \in \langle \{ u_k \}\rangle\, $ the quotient will be isomorphic to 
$\, R\, $ on one hand and has a basis consisting of the images of the $\, \{ u\, t_j\, u^{-1} \}\, $ on the other hand, so by this identification one constructs a basis $\,\{ g_k \}\, $ having a core for the quotient modulo 
$\, [\, R\, ,\, R\, ]\, $. However, it seems extremely difficult to find an appropriate notion of convergence in this case which is hereditary with respect to the suspension.
An even more general concept which is useful at times is given by the following Definition.
\par\bigskip\noindent
{\bf Definition 8.} Assume that $\, N \subset F\, $ is a normal subgroup and choose a section 
$\, s : F / N \nearrow F\, $. We will say that {\fndef $\,\bf N\, $ has a weak core for}  $\,\bf F\, $ if 
$\, N / [ N , N ]\, $ is free abelian and has a basis of the form 
$\,\{\, f_{i , x}\, =\, s ( x )\, e_i\, s ( x )^{-1}\,\}\, $ such that $\, f_{i , x} \neq f_{j , y}\, $ for $\, i \neq j\, $ or 
$\, x \neq y\, $.
\par\bigskip\noindent
{\bf Lemma 8.} Suppose that $\, N\, $ is free and has a weak core for  $\, F\, $. Then there exists an $F$-normal splitting $\, N \longrightarrow U_{N , F} / [\, U_{N , F}\, ,\, J_F\, ]\, $ of the quotient map. In particular 
$\, K^J_2 ( N , F )\, $ is trivial.
\par\bigskip\noindent
{\it Proof.}\quad Since $\, N\, $ is free there exists a splitting in any case. Moreover the subgroup 
corresponding to $\, [ N , N ]\, $ in the image is $F$-normal as an elementary computation shows. Dividing by this normal subgroup the quotient is abelian and since $\, N\, $ has a weak core for $\, F\, $ there an $F$-normal lift of the free abelian group $\, N / [ N , N ]\, $ to this abelian quotient of 
$\, U_{N , F}\, $. The preimage of this subgroup with respect to the normal copy of $\, [ N , N ]\, $ then gives a $F$-normal copy of $\, N\, $ in $\, U_{N , F} / [\, U_{N , F}\, ,\, J_F\, ]\, $ \qed   
\bigskip\bigskip
\par\noindent
\hfil\hfil \Large{\bf{Section IV -- Excision.}} \hfil
\bigskip\bigskip
\par\noindent

\noindent
{\bf Definition 9.}\quad For any normal pair define inductively
\smallskip
$$ ( J^{n+1}_{N , F}, U^{n+1}_F ) := ( J_{J^n_{N , F}, U^n_F} , U_{U^n_F} ). $$
\smallskip
Then define the {\fndef  $\bf K^J_n$-functors} of $\, ( N , F )\, $ by 
\smallskip
$$ K^J_n ( N , F ) := K^J_2 ( J^{n-2}_{N , F} , U^{n-2}_F ) ,\quad n>2 .$$
\par\medskip\noindent
Define the {\fndef extended $\bf K^J_2$-functor} of $\, ( N , F )\, $ by the formula 
\smallskip
$$ \widetilde K^J_2 ( N , F )\, =\, {J_{N , F} \cap \bigl[\, U_{N , F}\, ,\, U_F\, \bigr]  \over 
\bigl[\, J_{N , F}\, ,\, U_F\, \bigr]} \>  .  $$
\par\bigskip\noindent
Clearly, $\, K^J_n\, $ (and $\, \widetilde K^J_2\, $) give a (covariant) functors from the category of normal pairs to the category of abelian groups for all $\, n\geq 2\, $. Although we are primarily interested in the functors $\, K^J_n\, $ it is convenient for some applications to consider the extended $\, K^J_2$-functor.
In section 5 we will see that any surjective map 
$\, ( N , F )\rightarrow ( P , E )\, $ can be factored into a map $\, ( N , F )\rightarrow ( N/K , F/K )\, $ which will be called a {\fndef normal projection }, and a (regular) surjection of the form 
$\, ( P , F/K )\rightarrow ( P , E )\, $ inducing a surjective map on $K^J_n$-groups. Consider an arbitrary embedding $\, ( N , F )\subseteq ( M , G )\, $. Put 
$\, M_F = M\cap F\, $ and $\, M_F^G = G M_F G^{-1}\, $, then our map factors as
$\, ( N , F )\subseteq ( M_F , F )\subseteq ( M_F^G , G )\subseteq ( M , G )\, $. The first and the last map will be called {\fndef normal inclusions}. Then the intermediate map is a natural candidate for a 
(possibly regular) injection, i.e. an inclusion inducing an injective map on ${\rm K}^J_2$-groups. Unfortunately, things are not so easy. Let $\, ( N , F )\subseteq ( M , G )\, $ be an inclusion satisfying 
$\, N = M_F = M\cap F\, $. It can again be factored into the composition
$\, ( N , F )\subseteq ( M ,\, F+M )\subseteq ( M , G )\, $ where both parts are of the "same type", i.e.
$\, N = M_F\, $ and $\, M = M_{M+F}\, $. The Excision Theorem below gives a (partial) answer to the question under what conditions one of the two types of inclusions above or both give rise to an injection in $K^J$-theory.  In low dimensions ($\, n = 2,3\, $) one has to make an additional assumption concerning the inclusion $\, N \subseteq M\, $ (resp. $\, F / N \subseteq G / M\, $). Any inclusion 
$\, ( N , F ) \subseteq ( M , G )\, $ with $\, N = M \cap F\, $ and such that 
$\, K^J_2 ( N , F ) \rightarrowtail K^J_2 ( M , G )\, $ is injective will be called {\fndef preexcisive }. An example is an inclusion $\, N \subseteq M\, $ where $\, M\, $ is a product $\, N \cdot M_0\, $  It means  that there exists a "complementary" subgroup $\, M_0 \subseteq M\, $ for $\, N\, $ such that every left coset of $\, N\, $ contains one and only one element of $\, M_0\, $. Check that this implies the same for right cosets of $\, N\, $. If $\, N\, $ or $\, M_0\, $ is normal in $\, M\, $ one gets a semidirect product.
Up to now this is the only known condition ensuring injectivity of 
$\, K^J_2 ( N ) \rightarrowtail K^J_2 ( M )\, $ except trivial ones but this property does not seem to pass on to $\, K^J_3\, $.
The name "excision" is in some accordance with the use of this notion in ordinary algebraic $K$-theory (in the case when $\, N = M\, $) but at the same time much more general and also weaker (since it only gives an injectivity result and not an isomorphism). Indeed, in the classical case one looks at an ideal 
$\,\mathcal J\, $ of a unital ring $\,\mathcal R\, $ and asks whether the relative $K$-groups of the ideal are the same as the "absolute" $K$-groups. Translating into our setting the ideal is to be replaced by a normal subgroup $\, N\, $ of the group $\, E (\mathcal R )\, $ of elementary matrices over 
$\,\mathcal R\, $, which is the kernel corresponding to division by $\,\mathcal J\, $. Since $\, \mathcal J\, $ is not (necessarily) unital its "absolute" $K$-groups are given by the kernel of the evaluation map from the $K$-groups of its unitization to the $K$-groups of $\,\mathbb Z\, $. In certain cases there will be a splitting of the map 
$\,\mathbb Z \rightarrow {\mathcal R} /  {\mathcal J}\, $ corresponding to the condition 
$\, G / N = G_0 \rtimes F / N\, $ in the Lemma below (but in fact without loss of generality one can replace $\,\mathcal R\, $ by its unitization adjoining another unit showing that in this particular case the condition is not really essential).
\par\bigskip\noindent
{\bf Theorem 4. } ( Excision )\quad Let $\, ( N , F )\subseteq ( M , G )\, $ be an inclusion satisfying 
$\, N = M\, \cap\, F\, $. Then $\, K^J_n ( N , F ) \rightarrowtail K^J_n ( M , G )\, $ is injective for all 
$\, n \geq 4\, $.  If the image of $\, K^J_3 ( N )\, $ in $\, K^J_3 ( N , F )\, $ is trivial one gets an injective map $\, K^J_3 ( N , F ) \rightarrowtail K^J_3 ( M , G )\, $. If there exists a retraction 
$\, r : M \searrow N\, $ which is a homomorphism and $F$-equivariant modulo a given subgroup 
$\,N_0 \subseteq N\, $ (meaning that $\, r ( c\, d ) r ( d )^{-1} r ( c )^{-1} \in N_0\, $ and  
$\, r ( \xi\, c\, {\xi }^{-1} )\, \xi\, r ( c )^{-1}\, {\xi }^{-1} \in N_0\, $ for $\, c , d \in M\, ,\, \xi \in F\, $), the kernel of 
$\, K^J_2 ( N , F ) \rightarrow K^J_2 ( M , M + F )\, $ is contained in the image of 
$\, K^J_2 ( N_0 , F )\, $. Considering the case $\, N = M\, $ suppose there exists a retraction 
$\, r : G / N\, \searrow\, F / N \, $ which is a homomorphism modulo a subgroup $\, F_0 / N\, $ (the case 
$\, F_0 = N\, $ corresponds to a decomposition  $\, G / N\, =\, G_0 / N \rtimes F / N\, $) and let 
$\, G_0\, $ denote the subgroup generated by $\, F_0\, $ and all elements $\, \{ s ( x \cdot r ( x )^{-1} ) \}\, $  where $\, s : G / N \nearrow G\, $ is some chosen section. Suppose that $\, F_0\, $ is normal in $\, G\, $ and that the retraction $\, r\, $ has a left kernel, which is to say that 
$\, r ( G_{0 , 0} ) = 1\, $ where $\, G_{0 , 0}\, $ is the subgroup generated by the elements 
$\,\{ x \cdot r ( x )^{-1} \}\, $. By normality of 
$\, F_0\, ,\, G_0\, $ is normal in $\, G\, $ and equal to the preimage of $\, F_0\, $ for $\, r\, $.  Each element of $\, G\, $ can be expressed as a product of (the image by a given section of) an element of 
$\, G_{0 , 0}\, $ and an element of $\, F\, $. Assume that the following condition is fulfilled : 
$\, G_0\, $ contains the union of the $d$-commutants 
$\, \bigcup\, G_d\, $ where $\, G_d\, $ is the subgroup of elements commuting with a fixed element 
$\, d \in N\, $, and the union is taken over all  
$\, d \neq 1\, $. Also assume there exists a section $\, F / F_0 \nearrow F / N\, $ such that, denoting an element in its image by $\, \overline\xi\, $ the subgroup $\, G_{0 , 0}\, $ is normal for the adjoint action 
of $\, \overline\xi\, $, i.e. $\, \overline\xi\, x_0\, {\overline\xi }^{-1} \in G_{0 , 0}\, $ for every 
$\, x_0 \in G_{0 , 0}\, $ and denoting an element in $\, F_0 / N\, $ by $\, {\xi }_0\, $ there exists a section 
$\, s : G / N \nearrow G\, $ such that $\, s\, $ satisfies the relations $\, s ( x_0\, \overline\xi\, {\xi }_0 ) = s ( x_0 )\, s ( \overline\xi )\, s ( {\xi }_0 ) = 
s ( \overline\xi )\, s ( {\overline\xi }^{-1}\, x_0\, \overline\xi )\, s ( {\xi }_0 )\, $ and 
$\, s ( \overline\zeta\,\overline\xi ) = s ( \overline\zeta )\, s ( \overline\xi )\, $. Then 
$\, K^J_2 ( N , F ) \rightarrowtail K^J_2 ( N , G )\, $ is injective modulo the image of 
$\, K^J_2 ( N , F_0 )\, $. If $\, N\, $ is abelian and if the retraction 
$\, r : G / N\, \searrow\, F / N\, $ preserves the action on $\, N\, $ then the map 
$\, K^J_2 ( N , F ) \rightarrow K^J_2 ( N , G )\, $ is injective modulo the image of 
$\, K^J_2 ( N , F_0 )\, $ (even without the extra conditions, only assuming that $\, r\, $ is a homomorphism modulo $\, F_0 / N\, $). 
The map $\, \widetilde K^J_2 ( N , F ) \rightarrow \widetilde K^J_2 ( N , G )\, $ is always injective up to the image of $\,  K^J_2 ( N )\, $ and injective if $\, N\, $ is abelian and there exist a retraction 
$\, r : G / N \searrow F / N\, $ which preserves the action on $\, N\, $.
If $\, M = M_0\cdot N\, $ the map $\, K^J_2 ( N ) \rightarrowtail K^J_2 ( M )\, $ is injective.
\par\bigskip\noindent
{\it Proof.}\quad Throughout we use the terminology of Lemma 3 in section 1. Consider the subset
$\, B'_3\, =\, \{ u_c\, u_y\, u^{-1}_{cy} \}\, $ of $\, B_3\, $ (contained in the core of $\, J_{N , F}\, $ for 
$\, \lambda ( J_{F / N} )\, $) and note that a general element of $\, B_3\, $ can be written as
\smallskip
$$ u_c\, u_{dy}\, u^{-1}_{cdy}\, =\, \bigl[ u_c\, ( u_{dy}\, u^{-1}_y\, u^{-1}_d )\, u^{-1}_c \bigr]\,
( u_c\, u_d\, u^{-1}_{cd} )\, ( u_{cd}\, u_y\, u^{-1}_{cdy} )   , $$
\par\smallskip\noindent
so that dividing by the normalization of $\, B'_3\, $, the set $\, B_3\, $ becomes congruent with $\, B_1\, $.
Let now $\, ( N , F ) \subseteq ( M , G )\, $ be an inclusion such that $\, N = M\,\cap\, F\, $. One can choose a section $\,s : G / M \nearrow G\, $ restricting to $\, F / N \nearrow F\, $ giving compatible lifts
$\,\lambda : J_{G / M}\rightarrow J_G\, $ and $\, {\lambda }_F : J_{F / N}\rightarrow J_F\, $ together with compatible lifts $\, \sigma : U_{G / M}\rightarrow U_G\, $ and $\, {\sigma }_F : U_{F / N}\rightarrow U_F\, $, all of which are semicanonical with respect to the standard $U$- and $J$-bases (but note that
$\, \lambda ( J_{G / M} )\,\neq\,\sigma ( J_{G / M} )\, $ since the latter is not necessarily contained in 
$\, J_G $). This leads to natural inclusions of the cores of $\, U_{N , F}\, $ (resp. $\, J_{N , F} $) in the corresponding cores of $\, U_{M , G}\, $ (resp. $\, J_{M , G} $) for $\,\sigma\, $ (resp. $\, \lambda $). We may alter the standard core $\, \{ u_x\, u^{-1}_{\overline x} \}\, $ of $\, U_{M , G}\, $ replacing the elements 
$\, \{ u_x\, u^{-1}_{\overline x}\,\vert\, {\overline x} \neq 1 \}\, $ by the set 
$\, B'_3\, =\,\{ u_{x{\overline x}^{-1}}\, u_{\overline x}\, u^{-1}_x \}\, $ where $\, {\overline x}\, $ is in the image of the section $\, s\, $. Then $\, B'_3\, $ is a common subset for the (altered) core of 
$\, U_{M , G}\, $ as well as the $\, \{ u_z\, B\, u^{-1}_z \} $-core of $\, J_{M , G}\, $. Dividing by its normalization in 
$\, U_G\, $ gives a projection of $\, U_{M , G}\, $ onto the normal closure of $\, U_M\, $ in the subgroup of $\, U_G\, $ generated by $\, U_M\, $ and $\,\sigma ( U_{G / M} )\, $. It will send an element of 
$\, J_{M , G} \cap [\, U_{M , G}\, ,\, U_G\, ]\, $ representing an element of $\, K^J_2 ( M , G )\, $ to another representative in the same equivalence class and of the same degree of regularity because $\, B'_3\, $ is part of a basis of $\, U_G\, $. In particular trivial elements in 
$\, [\, J_{M , G}\, ,\, U_G\, ]\, +\, [\, U_{M , G}\, ,\, J_G\, ]\, $ are sent to trivial elements and k-fold commutators in 
$\, [\, [\cdots [\, U_{M , G}\, ,\, U_G\, ]\cdots ], U_G\, ]\, $ are sent to k-fold commutators. After this division the core of $\, U_{M , G}\, $ is reduced to the standard basis $\,\{ u_c \}\, $ of $\, U_M\, $ and the 
$\, \{ u_z\, B\, u^{-1}_z \}$-core of $\, J_{M , G}\, $ is reduced to 
$\, \{ u_z\, B_1\, u^{-1}_z \}\, +\,\{ u_z\, B_2\, u^{-1}_z \}\, $. All of this is effected in a compatible way for the inclusion
$\, ( J_{N , F} , U_{N , F} , U_F ) \subseteq ( J_{M , G} , U_{M , G} , U_G )\, $, since one has 
\smallskip
$$ J_{M , G}\,\cap\, [\, U_{M , G}\, ,\, U_G\, ]\,\cap U_F\, =\, J_{N , F}\,\cap\, [\, U_{N , F}\, ,\, U_F\, ]  $$
\par\smallskip\noindent
and
\smallskip
$$ \bigl[\, \bigl[\cdots\bigl[\, U_{M , G}\, ,\, U_G\, \bigr] ,\cdots\bigr] , U_G\,\bigr]\,\cap U_F\, =\,
\bigl[\, \bigl[\cdots\bigl[\, U_{N , F}\, ,\, U_F\, \bigr] ,\cdots\bigr] , U_F\, \bigr]   $$
\par\smallskip\noindent
for k-fold commutators. In the following $\, U_G\, ,\, U_F\, $ etc. will refer to this reduced form. We use letters $\, a , b , c , d ,\ldots\, $ to denote elements of $\, M\, $ and $\, x , y , z\, $ for elements in 
$\, s\, ( G / M )\, $. Correspondingly, elements of $\, N\, $ are denoted by greek letters $\,\alpha ,\beta ,\gamma ,\delta ,\ldots\, $ and we use $\,\xi ,\eta ,\zeta\, $ for $\, s\, ( F / N )\, $. Consider the subgroups of the (reduced) groups $\, U_G\, $ and $\, U_F\, $ generated by $\, J_G\, $ (resp. $\, J_F $) and 
$\, U_{M , G}\, $ (resp. $\, U_{N , F} $). One checks that a basis for these subgroups is given by the union of the sets $\, B_2\, ,\, \{\, u_c\,\}\, ,\,\{\, u_x\, u_y\, u^{-1}_{xy}\,\}\, $ 
(resp. $\, B^F_2\, ,\, \{\, u_\gamma\,\}\, ,\,\{ u_\xi\, u_\eta\, u^{-1}_{\xi\eta }\,\} $). Since $\, J_{M , G}\, $ has a core for $\, J_{G / M}\, $ one gets that its image dividing by $\, [\, J_{M , G}\, ,\, J_{G / M}\, ]\, $ is still a free group (resp. free abelian dividing by $\, [\, J_{M , G}\, ,\, J_G\, ] $) with basis given by the subcore
$\,\{\, u_z\, B_1\, u^{-1}_z\,\}\, +\,\{\, u_z\, B_2\, u^{-1}_z\,\}\, $. A representative for a k-regular element 
(k = 0 stands for a trivial element) can be expressed  modulo $\, [\, J_{M , G}\, ,\, J_G\, ]\, $ by a combination of elements (written symbolically)
\smallskip
$$ \bigl[\, u_z\, B_1\, u^{-1}_z\, ,\, u_c\, u_x\,\bigr]\, +\, \bigl[\, u_z\, B_2\, u^{-1}_z\, ,\, u_c\, u_x\,\bigr]\,
+\,\bigl[\, u_c\, ,\, u_x\, u_y\, u^{-1}_{xy}\,\bigr]\,    $$
$$  +\,\bigl[\,\bigl[\cdots\bigl[\, u_a\, ,\, u_{r_1}\,\bigr] ,\cdots \bigr] , u_{r_k}\,\bigr]   $$
\par\smallskip\noindent
where $\, r_1 ,\ldots , r_k \in M \cup s\, ( G / M )\, $, on a first glance, but taking a second look one reduces to the form
\smallskip
$$ \bigl[\, B_1\, ,\, u_c\, \bigr]\, +\, \bigl[\, B_1\, ,\, u_z\,\bigr]\, +\, \bigl[\, B_2\, ,\, u_c\,\bigr]\, +
\,\bigl[\, B_2\, ,\, u_z\,\bigr]\, +\,\bigl[\, u_c\, ,\, u_x\, u_y\, u^{-1}_{xy}\,\bigr]   $$
$$ +\, \bigl[\,\bigl[\cdots\bigl[\, u_a\, ,\, u_{r_1}\,\bigr] ,\cdots\bigr] , u_{r_k}\,\bigr]\,  .    $$
\par\smallskip\noindent 
$\, [\, J_{M , G}\, ,\, J_G\, ]\, $ adds commutators of elements in the basis $\, B_1\, \cup\, B_2\,\cup\,
[\, B_1\, ,\, u_z\, ]\,\cup\, [\, B_2\, ,\, u_z\, ]\, $ with elements conbined from this core and the standard basis 
$\,\{\, u_x\, u_y\, u^{-1}_{xy}\,\}\, $ of $\,\lambda (J_{G / M})\, $. An element of $\, B_2\, $ may be written in the form $\, (\, [\, u_x , u_d\, ]\, u_d\, u^{-1}_{d_x}\, )\, $ where $\, d_x\, =\, x\, d\, x^{-1}\, $. Also put
$\, xy\, =\, c_{x,y}\, (xy)\, $ with $\, c_{x,y} \in M\, $ and $\, (xy) \in s\, ( G / M )\, $ (note that this differs from a previous definition !), then one gets $\, u_{xy}\,\equiv\, u_{c_{x,y}}\, u_{(xy)}\, $ modulo $\, B'_3\, $, and an element in $\, [\, u_z\, ,\, B_2\, ]\, $ can be written modulo $\, [\, J_{M , G}\, ,\, J_{M , G}\, ]\, $ as
\smallskip 
$$\bigl[ u_z , \bigl( \bigl[\, u_x , u_d\,\bigr]\, u_d\, u^{-1}_{d_x}\,\bigr)\bigr]\,\equiv\, 
\bigl[\, u_z\, u_x\, u^{-1}_{zx} , u_{c_{x,y}}\, \bigl(\bigl[\, u_{(zx)} , u_d\,\bigr]\, 
u_d\, u^{-1}_{d_{(zx)}}\,\bigr)\, u^{-1}_{c_{x,y}}\,\bigr]  $$
$$ \bigl[\, u_z\, u_x\, u^{-1}_{zx}\, ,\, u_{c_{z,x}}\, u_{d_{(zx)}}\, u^{-1}_{c_{z,x}}\, u^{-1}_{d_{zx}}\,\bigr]\quad\bigl[\, u_{c_{z,x}}\, ,\, \bigl(\,\bigl[\, u_{(zx)}\, ,\, u_d\,\bigr]\, u_d\, u^{-1}_{d_{(zx)}}\,\bigr)\,\bigr]\! $$
$$ \bigl(\,\bigl[\, u_{(zx)}\, ,\, u_d\,\bigr]\, u_d\, u^{-1}_{d_{(zx)}}\,\bigr)\quad  
\bigl(\, u_{d_{zx}}\, u^{-1}_{d_x}\,\bigl[\, u_{d_x}\, ,\, u_z\,\bigr]\,\bigr)\quad
\bigl(\, u_{d_x}\, u^{-1}_d\,\bigl[\, u_d\, ,\, u_x\,\bigr]\,\bigr)\!  $$
$$ \bigl(\, u_{c_{z,x}}\, u_{d_{(zx)}}\, u^{-1}_{c_{z,x}}\, u^{-1}_{d_{zx}}\,\bigr)\quad
\bigl[\, u_z\, u_x\, u^{-1}_{zx}\, ,\, u_{d_{zx}}\,\bigr]                               $$
\par\smallskip\noindent
showing that $\, [\, u_z\, ,\, B_2\, ]\, $ and $\, [\, u_c\, ,\, u_x\, u_y\, u^{-1}_{xy}\, ]\, $ (as sets) are congruent modulo $\, B_1\, +\, B_2\, +\, [\, B_2\, ,\, u_c\, ]\, +\,\{\, u_x\, u_y\, u^{-1}_{xy}\,\}\, $. Similarly, an element of 
$\, [\, u_z\, ,\, B_1\, ]\, $ can be represented modulo $\, [\, J_{M , G}\, ,\, J_{M , G}\, ]\, $ by the combination $\, \bigl[\, u_z\, ,\, u_c\, u_d\, u^{-1}_{cd}\,\bigr]\,\equiv\, $
\smallskip
$$ \bigl(\,\bigl[\, u_z\, ,\, u_c\,\bigr]\, u_c\, u^{-1}_{c_z}\,\bigr)\quad
\bigl[\, u_{c_z}\, ,\, \bigl(\,\bigl[\, u_z\, ,\, u_d\,\bigr]\, u_d\, u^{-1}_{d_z}\,\bigr)\,\bigr]\quad  
\bigl(\,\bigl[\, u_z\, ,\, u_d\,\bigr]\, u_d\, u^{-1}_{d_z}\,\bigr)  $$
$$ \bigl(\, u_{c_z d_z}\, u^{-1}_{cd}\,\bigl[\, u_{cd}\, ,\, u_z\,\bigr]\,\bigr)\quad
\bigl(\, u_{c_z}\, u_{d_z}\, u^{-1}_{c_z d_z}\, \bigr)\quad \bigl(\, u_{cd}\, u^{-1}_d\, u^{-1}_c\,\bigr) $$
\par\smallskip\noindent
so that modulo $\, B_1 + B_2\, $ the sets $\, [\, B_1\, ,\, u_z\, ]\, $ and $\, [\, B_2\, ,\, u_c\, ]\, $ are congruent. We may then replace the basis elements $\, [\, B_2\, ,\, u_z\, ]\, $ by the corresponding elements in $\, [\, u_c\, ,\, u_x\, u_y\, u^{-1}_{xy}\, ]\, $ and $\, [\, B_1\, ,\, u_z\, ]\, $ by the elements of 
$\, [\, B_2\, ,\, u_c\, ]\, $ if we wish to do so, which then together with $\, B_1\, ,\, B_2\, $ and 
$\,\{\, u_x\, u_y\, u^{-1}_{xy}\,\}\, $ give another basis of $\, J_G\, $. 
\par\medskip\noindent
We begin by showing that an inclusion of type $\, ( N , F ) \subseteq ( N , G )\, $ induces an injection 
modulo the image of $\, K^J_2 ( N )\, $ on 
$\, \widetilde K^J_2\, $. 
A trivial element for $\, \widetilde K^J_2 ( N , G )\, $ can be represented modulo 
$\, [\, J_{N , G}\, ,\, J_G\, ]\, $ by a combination of elements (written symbolically) 
\smallskip
$$ \bigl[\, u_z\, ,\, B_2\, \bigr]\> +\> \bigl[\, u_z\, ,\, B_1\, \bigr]\> +\> \bigl[\, u_c\, ,\, B_2\,\bigr]\> +\> 
\bigl[\, u_c\, ,\, B_1\,\bigr]  $$
\par\medskip\noindent
We may project such a representative to $\, J_{N , F}\, $ by considering the basis of $\, J_{N , G}\, $ consisting of (iterated) commutators of elements of the core $\, B_1 , B_2 , [\, u_z\, ,\, B_1\, ] , 
[\, u_z\, ,\, B_2\, ]\, $ by basis elements $\, u_x\, u_y\, u^{-1}_{xy}\, $ of $\,\lambda ( J_{G / N} )\, $ and elements of the core and projecting the first (commutator) part of the basis to 
$\, [\, J_{N , F}\, ,\, \lambda ( J_{F / N} )\, ]\, $ in the obvious way eliminating a basis element which is not a basis element of the latter group and also projecting the elements of the core in the same manner. The formula for the interdependence of elements $\, [\, u_c\, ,\, B_2\, ]\, $ and $\, [\, u_z\, ,\, B_1\, ]\, $ given above shows that this projection maps a trivial element to a trivial element modulo an element of
$\, J_N\, $ which must be contained in $\, J_N \cap [\, U_N , U_N\, ]\, $ since the projection leaves invariant our representative which comes from $\, J_{N , F} \cap [\, U_{N , F}\, ,\, U_F\, ]\, $. Now assume that $\, N\, $ is abelian and one is given an action preserving retraction $\, r : G / N \searrow F / N\, $. We may the use the same projection as above on $\, [\, J_{N , G}\, ,\, \lambda ( J_{G / N} )\, ]\, $ but modify the projection on the core by sending $\, u_z\, $ to $\, u_{r ( z )}\, $ and $\, u_c\, $ to $\, u_c\, $. One checks that this projection maps a trivial element to a trivial one giving the result. 
\par\medskip\noindent
Next we consider the case $\, N = M\, $ for $\, K^J_2\, $. 
Let $\, r : G / N\,\searrow\, F / N\, $ be the given retraction and $\, x\, $ (resp. $\, z\, $) an arbitrary element of $\, G / N\, $. For simplicity of notation we will write $\, \xi\, $ for the element $\, r ( x )\, $ (resp. $\,\zeta\, $ for $\, r ( z )\, $) and $\, x_0 = x\,{\xi }^{-1}\, $ (resp. $\, z_0 = z\, {\zeta }^{-1}\, $) for the difference. 
Let a trivial element in $\, [\, J_{N , G}\, ,\, U_G\, ]\, +\, 
[\, U_{N , G}\, ,\, J_G\, ]\, $ be given represented modulo $\, [\, J_{N , G}\, ,\, J_G\, ]\, $ by a symbolic expression as above and assume that it comes from $\, J_{N , F}\,\cap\, [\, U_{N , F}\, ,\, U_F\, ]\, $. We want to show that the terms $\, [\, u_z\, ,\, B_2\, ]\, $ can be reduced to the four cases
\smallskip 
$$ \bigl[\, u_{\zeta }\, ,\, B^{\xi }_2\, \bigr]\, ,\, \bigl[\, u_{\zeta }\, ,\, B^{x_0}_2\, \bigr]\, ,\, 
\bigl[\, u_{z_0}\, ,\, B^{\xi }_2\, \bigr]\, ,\, \bigl[\, u_{z_0}\, ,\, B^{x_0}_2\, \bigr] $$
\par\smallskip\noindent
modulo trivial elements of other types. A staightforward calculation replacing 
$\,u_x\, $ by $\, u_{x_0}\, u_{\xi }\, $ and 
$\, u_z\, $ by 
$\, u_{z_0}\,u_{\zeta }\, $ respectively by $\, u^{-1}_{c_{\eta\, ,\, y_0}}\, u_{\eta }\, u_{y_0}\, $ modulo 
$\, J_G \, $  (note that the section 
$\, s : G / N\,\nearrow\, G\, $ satisfies $\, s ( x ) = s ( x_0 )\, s ( \xi )\, $ so the terms $\, c_{x_0\, ,\, \xi}\, $ vanish) shows that one can reduce the $\, [\, u_z\, ,\, B_2\, ]\, $ terms to the four cases above plus the cases 
$\, u_{z_0 }\, [\, u_{\zeta }\, ,\, B^{\xi }_2\, ]\, u^{-1}_{z_0 }\, $ and 
$\, u_{\zeta }\, [\, u_{z_0 }\, ,\, B_2^{x_0}\, ]\, u^{-1}_{\zeta }\, $ modulo $\, [\, J_{N , G}\, ,\, J_G\, ]\, $ plus expressions of type $\, [\, u_c\, ,\, u_x\, u_y\, u^{-1}_{xy}\, ]\, ,\, [\, u_c\, ,\, B^{x_0}_2\, ]\, $ and trivial expressions in $\, J_{N , F} \cap [\, U_{N , F}\, ,\, U_F\, ]\, $.  Then replacing the inner bracket of the latter modulo 
$\, [\, J_{N , G}\, ,\, J_{N , G}\, ]\, $ by the corresponding expression in the elements 
$\, [\, u_c\, ,\, u_x\, u_y\, u^{-1}_{xy}\, ]\, $ and the other terms given above and noting that the space of expressions $\,\langle\{ [\, u_c\, ,\, u_x\, u_y\, u^{-1}_{xy}\, ]\}\rangle\, $ is normal modulo 
$\, [\, J_{N , G}\, ,\, J_G\, ]\, +\,  [\, J_N\, ,\, U_N\, ]\, $ and elements of type 
$\, [\, u_d\, ,\, B^{c_{x , y}}_2\, ]\, $, the two last cases are reduced again to the four other cases (modulo expressions of different type). Considering the conjugacy classes of $\, N\, $ for the adjoint action of 
$\, G_0\, $, one gets by normality of $\, G_0\, $ that $\, F / F_0\, $ acts on the set of these classes and 
by the assumption that $\, G_0\, $ contains all $d$-commutants this action is free, so that $\, F / F_0\, $ permutes these classes. Choose a representative class for this action in each $G$-conjugacy class and denote the elements in this class by $\, \{ d^0 \}\, $.  By assumption on the section 
$\, s : G / N \nearrow G\, $ the corresponding curvature coefficients $\, c_{x , y}\, $ satisfy the relations 
$\, c_{ z_0\, ,\, \overline\xi\, {\xi }_0 } = 1\, , c_{{\overline\zeta }\, ,\, x_0\, {\xi }_0} = 1\, $ and 
$\, c_{\overline\zeta\, ,\, \overline\xi } = 1\, $  where 
$\, {\zeta }_0 \in s ( F_0 / N )\, ,\, 
z_0\, ,\, x_0  \in s ( G_{0 , 0} )\, $ and $\, \{ \overline \zeta \}\, $ denotes the image  of the section  
$\, F / F_0 \nearrow F / N\, $ given in Theorem 4. One finds that one can reduce the elements of type 
$\, [\, u_z\, ,\, B_2 \, ]\, $ further modulo trivial elements of other types to the types 
\smallskip
$$  \bigl[\, u_{\zeta }\, ,\, B_2^{\xi }\, \bigr]\> ,\> 
\bigl[\, u_{\overline \zeta}\, ,\, B_2^{ x_0}\,\bigr] \> ,\> 
 \bigl[\, u_{{\zeta }_0}\, ,\, B_2^{x_0}\,\bigr] \> ,\> 
 \bigl[\, u_{z_0}\, ,\, B_2^{\xi }\,\bigr] \> ,\> 
 \bigl[\, u_{z_0}\, ,\, B_2^{x_0}\,\bigr]\> .  $$
\par\medskip\noindent
A similar reduction can be carried through for the elements of type 
$\, [\, u_z\, ,\, B_1\, ]\, $ and $\, [\, u_c\, ,\, B_2\, ]\, $  modulo elements 
$\, [\, u_c\, ,\, u_x\, u_y\, u^{-1}_{xy}\, ]\, $ and certain $U_N$-conjugates of the types above splitting them into the types
\smallskip
$$ \bigl[\, u_{\zeta }\, ,\, B_1\,\bigr]\> ,\> \bigl[\, u_{z_0}\, ,\, B_1\,\bigr]\> ,\> 
\bigl[\, u_c\, ,\, B_2^{\xi }\,\bigr]\> ,\> \bigl[\, u_c\, ,\, B_2^{x_0}\,\bigr] $$
\par\medskip\noindent
with $\,\zeta \in F\, ,\, x_0  \in G_{0 , 0}\, $ modulo $\, [\, J_{N , G}\, ,\, J_G\, ]\, $, the 
$\, [\, u_z\, ,\, B^x_2\, ]\, $ types as above, elements 
$\, [\, u_c\, ,\, u_x\, u_y\, u^{-1}_{xy}\, ]\, $ (we allow arbitrary $U_N$-coefficients for the 
$[\, u_z\, ,\, B_2\, ]$-types in our list except those of type $\, [\, u_{\overline\zeta }\, ,\, B^{x_0}_2\, ]\, $ as well as for the elements $\, [\, u_c\, ,\, u_x\, u_y\, u^{-1}_{xy}\, ]\, $) and trivial elements in $\, J_{N , F}\, $. 
It is possible to further reduce the types 
$\, [\, u_{z_0}\, ,\, B_2^{\overline\xi }\, ]\, $ and $\, [\, u_{z_0}\, ,\, B_2^{x_0}\, ]\, $ to types
\smallskip
$$   u_e\, \bigl[\, u_{z_0}\, ,\, B_2^{\overline\xi   , d^0}\,\bigr]\, u^{-1}_e \quad ,\quad
u_e\, \bigl[ u_{z_0}\, ,\, B_2^{x_0 , d^0}\,\bigr]\, u^{-1}_e   $$
\par\medskip\noindent
modulo $\, [\, J_{N , G}\, ,\, J_G\, ]\, $, trivial elements in $\, J_{N , F}\, $  plus elements of type 
$\, [\, u_c\, ,\, u_x\, u_y\, u^{-1}_{xy}\, ]\, $ (by assumption on the section $\, s\, $ this reduction does not 
result in additional elements of type $\, [\, u_{z_0}\, ,\, B_1\, ]\, $ and $\, [\, u_c\, ,\, B^{x_0}_2\, ]\, $).
The elements of the second type can be expressed by a combination of elements 
$\, u_e\, [\, u_{x_0}\, ,\, u_{d^0}\, ]\, u^{-1}_e\, $ and an expression in $\, [\, U_N\, ,\, U_N\, ]\, $ modulo trivial elements of other types (compare the formula above connecting the $[ u_z\, ,\, B_2 ]$-elements and corresponding elements of type $\, [\, u_c\, ,\, u_x\, u_y\, u^{-1}_{xy}\, ]\, $).   
Modulo 
$\ [\, J_{N , G}\, ,\, J_G\, ]\, $ and $[\, u_z\, ,\, B_2\, ]$-elements contained in our list we may reduce the elements of type 
$\, [\, u_{z_0}\, ,\, B_1\, ]\, $ and 
$\, [\, u_c\, ,\, B_2^{x_0}\, ]\, $ to
\smallskip
$$  \bigl[\, u_{z_0}\, ,\, u_c\, u_{d_1}\, u^{-1}_{c d_1}\,\bigr]\quad ,\quad 
\bigl[\, u_c\, ,\, B_2^{x_0 , d^0}\,\bigr] $$
\par\medskip\noindent
where $\, d^0\, $ is taken from the set of representatives for the conjugacy classes of 
$\, N \, $ modulo the action of $\, s ( F / F_0 )\, $ as explained and 
$\, d_1 \notin \{ d^0 \}\, $ can be assumed. Additional coefficients $\, u_e\, $ can be suppressed in case of $\, [\, u_{z_0}\, ,\, u_c\, u_{d_1}\, u^{-1}_{c d_1}\, ]\, $  by the assumption that the section satisfies 
$\, s ( x_0 )\, s ( \overline\zeta ) = 
s ( \overline\zeta )\, s ( {\overline\zeta }^{-1} x_0 {\overline\zeta } )\, $. One now reduces the subgroup 
$\, [\, J_{N , G}\, ,\, J_G\, ]\, $  using the core 
\smallskip
$$ [\, u_c\, ,\, u_x\, u_y\, u^{-1}_{xy}\, ]\> ,\> u_c\, B_2\, u^{-1}_c\> ,\> B_2\> ,\> B_1  $$ 
\par\smallskip\noindent
for $\, \lambda (J_{G / N})\, $ by dividing out as much as possible of the commutators 
$\, [\, u_x\, u_y\, u^{-1}_{xy}\, ,\, J_{N , G}\, ]\, $ without altering expressions of type
$\, u_e\, [\, u_z\, ,\, B_2\, ]\, u^{-1}_e\, $ contained in our list. The congruence of $\, [\, u_z\, ,\, B_2\, ]\, $ with 
$\, [\, u_c\, ,\, u_x\, u_y\, u^{-1}_{xy}\, ]\, $ above shows that we can divide by all ($\, J_G$-conjugates of) double commutators of the form 
$\, [\, u_x\, u_y\, u^{-1}_{xy}\, ,\, [\, u_z\, u_w\, u^{-1}_{zw}\, ,\, J_{N , G}\, ]\, ]\, $ with either 
$\, x\, ,\, y\, $,$\,  z\, $ or $\, w\, $ not in $\, F\, $. Also the basis elements 
$\, [\, u_w\, u_z\, u^{-1}_{wz}\, ,\, [\, u_c\, ,\, u_x\, u_y\, u^{-1}_{xy}\, ]\, ] $, and all elements 
$\, [\, u_x\, u_y\, u^{-1}_{xy}\, ,\, u_c\, B^z_2\, u^{-1}_c\, ]\, $ with 
$\, z\,\neq\, (xy)\, $ or $\, ( x , y )\, $ is not equal to one of the types as above may be divided out. Modulo commutators by elements from the core of $\, J_{N , G}\, $ this leaves us with the following types (and certain 
$\,\{\, u_c\,\}$-conjugates thereof) which are
\smallskip
$$ \bigl[\, u_c\, ,\, u_x\, u_y\, u^{-1}_{xy}\, \bigr]\> ,\> \bigl[\, B_1\, ,\, u_z\,\bigr]\> ,\> \bigl[\, B_2\, ,\, u_c\, \bigr]    $$
\par\smallskip\noindent
the $\, [\, u_z\, ,\, B^x_2\, ]\, $ types as above and
\smallskip
$$ \bigl[\, u_z\, u_x\, u^{-1}_{zx}\, ,\, u_c\, B^{(zx)}_2\, u^{-1}_c\, \bigr]\>\equiv\>
 u_c\, \bigl[\, u_z\, u_x\, u^{-1}_{zx}\, ,\, \bigl[\, u_{(zx)}\, ,\, u_d\,\bigr]\,\bigr]\, u^{-1}_c  $$
\par\smallskip\noindent
modulo the former types (note that also the element $\, [\, u_x\, u_y\, u^{-1}_{xy}\, ,\, B_1\, ]\, $ is represented by the former types). 
Consider the following basis for $\, [\, U_{N , G}\, ,\, U_G\, ]\, $ consisting of the union of the two sets
\smallskip 
$$ \bigl\{\, u^{\pm}_{c_1}\cdots u^{\pm}_{c_m}\, \bigl[\, u_c\, ,\, u_d\,\bigr]\, u^{\mp}_{c_m}\cdots 
u^{\mp}_{c_1}\,\bigr\}\quad\cup $$ 
$$ \bigl\{\, u^{\pm }_{c_1}\cdots u^{\pm }_{c_m}\, u^{\pm }_{x_1}\cdots u^{\pm }_{x_n}\,
\bigl[\, u_x\, ,\, u_c\, \bigr]\, u^{\mp}_{x_n}\cdots u^{\mp }_{x_1}\, u^{\mp }_{c_m}\cdots u^{\mp }_{c_1}\,
\bigr\}   $$
\par\smallskip\noindent
where $\, c_1\leq\cdots\leq c_m\leq d\, ,\, c < d\, $ in the first set and $\, c_1\leq\cdots\leq c_m\, $ in the second (with respect to a chosen order on the elements of $\, N\, $ by the ordering axiom ! ) It may be divided into two parts. Let $\, S_0\, $ be the subspace generated by conjugates of basic commutators 
$\, [\, u_x\, ,\, u_d\, ]\, $ with $\, u^{\pm }_{x_1}\cdots u^{\pm }_{x_n} \in \sigma ( J_{G / N} )\, $ and arbitrary $\, u^{\pm }_{c_1}\cdots u^{\pm }_{c_m}\, $ together with the first set of the basis above. Then the complementary space $\, S_1\, $ consists of those conjugates of basic commutators 
$\, u_z\, [\, u_x\, ,\, u_d\, ]\, u^{-1}_z\, $ ($\, z\, =\, x^{\pm }_1\cdots x^{\pm }_n\, $) with elements from 
$\, \sigma ( J_{G / N} )\, $ and arbitrary $\, u^{\pm }_{c_1}\cdots u^{\pm }_{c_m}\, $. We will make a change of basis sending $\, u_z\, [\, u_x\, ,\, u_d\, ]\, u^{-1}_z\, $ to
\smallskip 
$$ \bigl[\, u_z\, ,\, \bigl(\, \bigl[\, u_x\, ,\, u_d\, \bigr]\, u_{[d,x]}\, \bigr)\, \bigr]\,
=\, \bigl(\, u_z\, \bigl[\, u_x\, ,\, u_d\, \bigr]\, u^{-1}_z\,\bigr)\, \bigl[\, u_z\, ,\, u_{[d,x]}\,\bigr]\, 
\bigl[\, u_d\, ,\, u_x\,\bigr]    $$
\par\smallskip\noindent
if $\, ( z , x )\, $ corresponds to one of the remaining $\, [\, u_z\, ,\, B_2^x\, ]\, $ types and 
$\, x \neq (z^{-1})\, $ in order that the elements $\, [\, u_z\, ,\, B_2\, ]\, $ (we have changed 
$\, [\, u_z\, ,\, B_2\, ]\, $ to this form modulo 
$\, [\, B_2\, ,\, [\,B_1\, ,\, u_z\, ]\, ]\, +\, [\, B_1\, ,\, u_z\, ]\, $) are contained in  $\, S_1\, $. Let us write down the expression for 
$ u_{c_{x,y}}\, [\, u_x\, u_y\, u^{-1}_{xy} , u_c\, ]\, u^{-1}_{c_{x,y}}\, $ in this basis. It reads 
$\, u_{c_{x,y}}\, [\, u_x\, u_y\, u^{-1}_{xy}\, ,\, u_c\, ]\, u^{-1}_{c_{x,y}}\> =\, $
$$ \bigl[\, u_{c_{x,y}}\, , u_x\,\bigr]\, \bigl( u_x\, \bigl[\, u_{c_{x,y}}\, , u_y\,\bigr]\, u^{-1}_x \bigr)\,
\bigl(\, u_x\, u_y\, u^{-1}_{(xy)}\,\bigl[\, u_{(xy)}\, ,\, u_{c_{x,y}}\,\bigr]\, u_{(xy)}\, u^{-1}_y\, u^{-1}_x\,\bigr)\, $$
$$ \bigl(\, u_x\> u_y\> u^{-1}_{(xy)}\,\bigl[\, u_c\, ,\, u_{(xy)}\, \bigr]\, u_{(xy)}\, u^{-1}_y\, u^{-1}_x\,\bigr)\quad \bigl(\, u_x\> \bigl[\, u_y\, ,\, u_c\,\bigr]\, u^{-1}_x\,\bigr)\, $$
$$ \bigl[\, u_x\, ,\, u_c\,\bigr] \quad \bigl(\, u_c\> u_x\> u_y\> u^{-1}_{(xy)}\,\bigl[\, u_{c_{x,y}}\, ,\, u_{(xy)}\,
\bigr]\, u_{(xy)}\> u^{-1}_y\, u^{-1}_x\, u^{-1}_c\,\bigr) $$
$$ \bigl(\, u_c\> u_x\>\bigl[\, u_y\, ,\, u_{c_{x,y}}\,\bigr]\, u^{-1}_x\, u^{-1}_c\,\bigr)\quad
\bigl(\, u_c\>\bigl[\, u_x\, ,\, u_{c_{x,y}}\,\bigr]\, u^{-1}_c\,\bigr)\quad\bigl[\, u_c\, ,\, u_{c_{x,y}}\,\bigr]  . $$    
\par\smallskip\noindent
If $\, (xy)\,\neq\, 1\, $ and $\, c\,\neq\, c_{x,y}\, $ we can change the basis elements 
\smallskip
$$ \bigl(\, u_x\> u_y\> u^{-1}_{(xy)}\, \bigl[\, u_c\, ,\, u_{(xy)}\, \bigr]\, u_{(xy)}\> u^{-1}_y\, u^{-1}_x\,\bigr) $$
\par\smallskip\noindent 
to $\, u_{c_{x,y}}\, [\, u_c\, ,\, u_x\, u_y\, u^{-1}_{xy}\, ]\, u^{-1}_{c_{x,y}}\, $ and for $\, (xy)\,\neq 1\, $ and 
$\, c\, =\, c_{x,y}\, $ we can change 
\smallskip
$$ \bigl(\, u_x\> u_y\> u^{-1}_{(xy)}\,\bigl[\, u_{c_{x,y}}\, ,\, u_{(xy)}\,\bigr]\, 
     u_{(xy)}\> u^{-1}_y\, u^{-1}_x\, \bigr)         $$ 
\par\smallskip\noindent
to $\, [\, u_{c_{x,y}}\, ,\, u_x\, u_y\, u^{-1}_{xy}\, ]\, $.  If $\, (xy) = 1\, $ and $\, c \neq c_{x , x^{-1}}\, $ we may change 
\smallskip
$$ u_x\, \bigl[\, u_{(x^{-1})}\, ,\, u_c\,\bigr]\, u^{-1}_x  $$
\par\smallskip\noindent
to $\, u_{c_{x , x^{-1}}}\, [\, u_c\, ,\, u_x\, u_{x^{-1}}\, ]\, u^{-1}_{c_{x , x^{-1}}}\, $ and if 
$\, c = c_{x , x^{-1}}\, $ change 
\smallskip
$$ u_x\, \bigl[\, u_{(x^{-1})}\, ,\, u_c\,\bigr]\, u^{-1}_x  $$
\par\smallskip\noindent
to $\, [\, u_c\, ,\, u_x\, u_{x^{-1}}\, ]\, $. We extend these changes to (ordered) 
$U_N$-conjugates of these basis elements. By this change of base one finds that the 
conjugates of (new) basis elements 
$\, \{\, [\, u_z\, ,\, (\, [\, u_x\, ,\, u_d\, ]\, u_d\, u^{-1}_{d_x}\, )\, ]\, \vert\, x \neq (z^{-1})\, \}\, $,
$\, \{\, [\, u_c\, ,\, u_z\, u_x\, u^{-1}_{zx}\, ]\, \}\, $ and 
$\, \{\, [\, u_x\, ,\, u_d\, ]\, \}\, $ and are all independent of each other. 
We will change the ($U_N$-conjugates of) basis elements 
$\, [\, u_{x_0\, \overline\zeta}\, ,\, u_d\, ]\, $ to the corresponding conjugates of
\smallskip 
$$ \bigl(\, \bigl[\, u_{x_0\, \overline\zeta}\, ,\, u_d\, \bigr]\, u_d\, u^{-1}_{d_{x_0 \overline\zeta }}\, \bigr)\> 
\bigl(\, u_{d_{{\overline\zeta } \cdot ({\overline\zeta }^{-1} x_0 \overline\zeta )}}\, 
u^{-1}_{d_{( {\overline\zeta }^{-1} x_0 {\overline\zeta } )}}\, 
\bigl[\, u_{d_{( {\overline\zeta }^{-1} x_0 {\overline\zeta } )}}\, ,\, u_{\overline\zeta }\, \bigr]\, \bigr) $$ 
$$ \bigl(\, u_{d_{( {\overline\zeta }^{-1} x_0 {\overline\zeta } )}}\, u_d^{-1}\, 
\bigl[\, u_d\, ,\, u_{( {\overline\zeta }^{-1} x_0 {\overline\zeta } )}\, \bigr]\, \bigr)\>
\bigl(\, u_{c_{{\overline\zeta }\, ,\, ({\overline\zeta }^{-1} x_0 {\overline\zeta } )}}\, 
u_{d_{{\overline\zeta } \cdot ({\overline\zeta }^{-1} x_0 {\overline\zeta })}}\, 
u^{-1}_{c_{{\overline\zeta }\, ,\, ({\overline\zeta }^{-1} x_0 {\overline\zeta } )}}\, 
u^{-1}_{d_{x_0 {\overline\zeta }}}\, \bigr)\> . $$
\par\medskip\noindent
Correspondingly, we change the $U_N$-conjugates of $\, [\, u_{x_0 {\zeta }_0}\, ,\, u_d\, ]\, $ for 
$\, d \notin \{ d^0 \}\, $ to the form 
\smallskip
$$ \bigl(\, \bigl[\, u_{x_0 {\zeta }_0}\, u_d\,\bigr]\, u_d\, u^{-1}_{d_{x_0 {\zeta }_0}}\,\bigr)\quad 
\bigl(\, u_{d_{( x_0 {\zeta }_0 x_0^{-1} ) \cdot x_0 }}\, 
u^{-1}_{d_{x_0}}\, 
\bigl[\, u_{d_{x_0}}\, ,\, u_{( x_0 {\zeta }_0 x_0^{-1} )}\, \bigr]\, \bigr) $$
$$ \bigl(\, u_{d_{x_0}}\, u^{-1}_{d^0}\, 
\bigl[\, u_d\, ,\, u_{x_0}\, \bigr]\, \bigr)\quad 
\bigl(\, u_{c_{( x_0 {\zeta }_0 x_0^{-1} )\, ,\, x_0}}\, 
u_{d_{( x_0 {\zeta }_0 x_0^{-1} ) \cdot x_0}}\, 
u^{-1}_{c_{( x_0 {\zeta }_0 x_0^{-1} )\, ,\, x_0}}\, 
u^{-1}_{d_{x_0 {\zeta }_0}}\, \bigr) $$
\par\medskip\noindent
We change the conjugates of basis elements 
$\, [\, u_{x_0\, \zeta }\, ,\, u_{d^0}\, ]\, $ for $\, \zeta  \notin \{ \overline\zeta \}\, ,\, \zeta \neq 1\, $ to the form
\smallskip
$$ \bigl(\, \bigl[\, u_{x_0\, \zeta }\, ,\, u_{d^0}\, \bigr]\, u_{d^0}\, u^{-1}_{d^0_{x_0\,\zeta }}\, \bigr)\> 
\bigl(\, u_{d^0_{\zeta }}\, u^{-1}_{d^0}\, \bigl[\, u_{d^0}\, ,\, u_{\zeta }\, \bigr]\, \bigr)\> 
\bigl(\, u_{d^0_{x_0\, \zeta }}\, u^{-1}_{d^0_{\zeta }}\, \bigl[\, u_{d^0_{\zeta }}\, ,\, u_{x_0}\, \bigr]\, \bigr)\> .  $$
\par\medskip\noindent 
Assume that $\, d = (d^0)_{\overline\zeta}\, $. Then we may change the basis elements 
$\, [\, u_{x_0}\, ,\, u_d\, ]\, $ (no $U_N$-conjugates) to the form
\smallskip
$$ \bigl(\, \bigl[\, u_{x_0}\, ,\, u_d\, \bigr]\, u_d\, u^{-1}_{d_{x_0}}\,\bigr)\quad 
\bigl(\, \bigl[\, u_{\overline\zeta }\, ,\, u_{d^0}\,\bigr]\, u_{d^0}\, u^{-1}_d\, \bigr)  $$
$$ \bigl(\> u_{(d^0)_{{\overline\zeta } \cdot ( {\overline\zeta }^{-1} x_0 {\overline\zeta } )}}\, 
u^{-1}_{(d^0)_{( {\overline\zeta }^{-1} x_0 {\overline\zeta } )}}\, 
\bigl[\, u_{(d^0)_{( {\overline\zeta }^{-1} x_0 {\overline\zeta } )}}\, ,\, u_{\overline\zeta }\,\bigr]\> \bigr) $$
$$ \bigl(\> u_{(d^0)_{( {\overline\zeta }^{-1} x_0 {\overline\zeta } )}}\, 
u^{-1}_{d^0}\, \bigl[\, u_{d^0}\, ,\, u_{( {\overline\zeta }^{-1} x_0 {\overline\zeta } )}\, \bigr]\> \bigr)\> 
\bigl(\> u_{c_{{\overline\zeta } \cdot ( {\overline\zeta }^{-1} x_0 {\overline\zeta } )}}\> 
u_{(d^0)_{{\overline\zeta } ( {\overline\zeta }^{-1} x_0 {\overline\zeta } )}}\> 
u^{-1}_{c_{{\overline\zeta } \cdot ( {\overline\zeta }^{-1} x_0 {\overline\zeta } )}}\> 
u^{-1}_{d_{x_0}}\> \bigr)\> .  $$
\par\medskip\noindent 
We also change the basis elements $\, u_c\, [\, u_{x_0}\, ,\, u_{d_1}\, ]\, u^{-1}_c\, $ to 
\smallskip
$$ \bigl[\, u_{x_0}\, ,\, u_c\, u_{d_1}\, u^{-1}_{c d_1}\,\bigr] $$
\par\medskip\noindent
if $\, d_1 \notin \{ d^0 \}\, ,\, c \neq 1\, $ leaving unchanged the other basis elements of 
$\, [\, U_{N , G}\, ,\, U_G\, ]\, $. Then consider the projection
\smallskip
$$ \bigl[\, U_{N , G}\, ,\, U_G\, \bigr]\,\longrightarrow\,\bigl[\, U_{N , F}\, ,\, U_F\,\bigr]  $$
\par\smallskip\noindent
sending the ($U_N$-conjugates of) basis elements $\, [\, u_c\, ,\, u_z\, u_x\, u^{-1}_{zx}\, ]\, $ to (conjugates of)
$\, [\, u_c\, ,\, u_{r (z )}\, u_{r ( x )}\, u^{-1}_{r ( z ) r ( x )}\, ]\, $, the corresponding conjugates of
$\, [\, u_z\, ,\, B^x_2\, ]\, $ with $\, ( z , x )\, $ one of the remaining types to (conjugates of) 
$\, [\, u_{r ( z ) }\, ,\, B^{r ( x )}_2\, ]\, $, the other new basis elements as above and the $U_N$-conjugates of  $\, [\, u_{x_0}\, ,\, u_{d^0}\, ]\, $ to zero and extending this projection on the complementary subbasis by means of the retraction $\, r\, $ (i.e. if some $U_G$-conjugate of a basic commutator 
$\, [\, u_x\, ,\, u_d\, ]\, $ is contained in the complementary subbasis then the basic commutator is projected as above and the conjugacy coefficients are projected by means of $\, r\, $). Using the formulas for the interdependence of the trivial elements of types 
$\, [\, u_c\, ,\, u_x\, u_y\, u^{-1}_{xy}\, ]\, $ and $\, [\, u_z\, ,\, B_2\, ]\, $ respectively  
$\, [\, u_z\, ,\, B_1\, ]\, $ and $\, [\, u_c\, ,\, B_2\, ]\, $ written out above one checks the following cases:
the terms giving the difference (modulo $\, [\, J_{N , G}\, ,\, J_G\, ]\, $) of elements of types 
$\, [\, u_{z_0}\, ,\, B_2^{\xi}\, ]\, $ ( $\, \xi \neq {\overline\xi }\, $, resp. 
$\, [\, u_{\overline\zeta}\, ,\, B_2^{x_0}\, ]\, $) and $\, u_e\, [\, u_{z_0}\, ,\, B_2^{x_0 , d^0}\, ]\, u^{-1}_e\, $ which appear in our list with corresponding elements $\, [\, u_c\, ,\, u_x\, u_y\, u^{-1}_{xy}\, ]\, $ are mapped to trivial elements modulo $\, [\, U_N\, ,\, U_N\, ]\, $ (this is true also in the case 
$\, [\, u_{z_0}\, ,\, B_2^{( z_0^{-1} ) , d^0}\, ]\, $). The expressions  
$\, [\, u_{z_0}\, ,\, u_c\, u_{d_1}\, u^{-1}_{c d_1}\, ]\, $ and 
$\, [\, u_c\, ,\, B_2^{x_0 , d^0}\, ]\, $ are sent to trivial elements modulo $\, [\, U_N\, ,\, U_N\, ]\, $. The difference terms modulo $\, [\, J_{N , G}\, ,\, J_G\, ]\, $ involving the elements of type 
$\, u_e\, [\, u_{{\zeta }_0 }\, ,\, B_2^{x_0 , d^0}\, ]\, u^{-1}_e\, $ appearing in our reduction give rise to the expression
\smallskip 
$$ \bigl(\, u_{d^0_{( x_0 {\zeta }_0 x_0^{-1} ) \cdot x_0 }}\, 
u^{-1}_{d^0_{x_0}}\, 
\bigl[\, u_{d^0_{x_0}}\, ,\, u_{( x_0 {\zeta }_0 x_0^{-1} )}\, \bigr]\, \bigr) $$ 
$$ \bigl(\, u_{d^0_{x_0}}\, u^{-1}_{d^0}\, 
\bigl[\, u_{d^0}\, ,\, u_{x_0}\, \bigr]\, \bigr)\>
\bigl(\, u_{c_{( x_0 {\zeta }_0 x_0^{-1} )\, ,\, x_0}}\, 
u_{d^0_{( x_0 {\zeta }_0 x_0^{-1} ) \cdot x_0}}\, 
u^{-1}_{c_{( x_0 {\zeta }_0 x_0^{-1} )\, ,\, x_0}}\, 
u^{-1}_{d_{x_0 {\zeta }_0}}\, \bigr) $$
$$ \bigl(\, \bigl[\, u_{x_0}\, ,\, u_{d^0_{{\zeta }_0}}\, \bigr]\, 
u_{d^0_{{\zeta }_0}}\, u^{-1}_{d^0_{x_0 {\zeta }_0}}\,\bigr)\> 
\bigl(\, \bigl[\, u_{{\zeta }_0}\, ,\, u_{d^0}\,\bigr]\, u_{d^0}\, u^{-1}_{d^0_{{\zeta }_0}}\,\bigr)\>   $$
\par\medskip\noindent
which maps to $\, [\, U_{N , F_0}\, ,\, U_{F_0}\, ]\, $. 
One checks that the projection sends the commutator subgroup of the core of $\, J_{N , G}\, $ for 
$\, \lambda ( J_{G / N} )\, $ to trivial elements plus 
$\, [\, U_N\, ,\, U_N\, ]\, $ (recall that we have used the core consisting of 
$\, [\, u_c\, ,\, u_x\, u_y\, u^{-1}_{xy}\, ]\, ,\, [\, u_c\, ,\, B_2\, ]\, ,\, B_2\, ,\, B_1\, $ for our first projection). Since our representative from $\, J_{N , F} \cap [\, U_{N , F}\, ,\, U_F\, ]\, $ must remain invariant, one induces that it lies in $\, J_{N , F_0} \cap [\, U_{N , F_0}\, ,\, U_{F_0}\, ]\, $. The argument in case that the retraction is action preserving is even much simpler. It then suffices to consider a single projection 
$\, J_{N , G} \cap [\, U_{N , G}\, ,\, U_G\, ] \rightarrow J_{N , F} \cap [\, U_{N , F}\, ,\, U_F\, ]\, $ induced by the obvious projection $\, u_z \mapsto u_{r ( z )}\, $ on the core elements and 
$\, u_x\, u_y\, u^{-1}_{xy} \mapsto u_{r ( x )}\, u_{r ( y )}\, u^{-1}_{r ( x ) r ( y )}\, $ (the main thing is that elements of type $\, B_2\, $ are again mapped to elements of type $\, B_2\, $).  
\par\medskip\noindent
We now consider the cases of higher $\, n\, $ beginning with the case $\, n = 3\, $. The condition 
$\, N = M \cap F\, $ implies $\, J_{N , F} = J_{M , G} \cap U_F\, $. We may then replace $\, U_F\, $ and 
$\, U_G\, $ by their reduced forms as above (for the inclusion $\, ( N , F ) \subseteq ( M , G )\, $) since the normalization of the $B_3'$-elements is a cone in both cases. Let $\,\{\, c_0\,\}\, $ be a collection of representatives of left cosets $\, N \backslash M\, $ in $\, M\, $ and $\,\{\, x_0\,\}\, $ a collection of representatives of left cosets $\, (F / N) \backslash (G / M)\, $ in $\, s (G / M)\, $. Consider the following basis of the (reduced) subgroup $\, J_G + U_F\, $  consisting of the union of the sets 
\smallskip
$$ \bigl\{\, u_{\gamma }\,\bigr\}\,\cup\,
\bigl\{\, u_{\xi }\,\bigr\}\,\cup\,\bigl\{\, u_x\, u_d\, u^{-1}_x\, u^{-1}_{d_x}\,\vert\, ( d , x ) \notin ( N , F )\,\bigr\}\,\cup\, $$ 
$$ \bigl\{\, u_{\xi }\, u_{x_0}\, u^{-1}_{\xi x_0}\, \vert\, x_0 \in (F / N) \backslash (G / M)\,\bigr\}\,\cup\,
\bigl\{\, u_{x_0}\, u_y\, u^{-1}_{x_0 y}\, \vert\, x_0 \in (F / N) \backslash (G / M)\,\bigr\}\, $$
$$ \cup\,\bigl\{\, [\, u_{z_0}\, ,\, u_{c_0}\, u_d\, u^{-1}_{c_0 d}\,]\,\vert\, z_0 \in (F / N) \backslash (G / M )\, ,\, c_0 \in N \backslash M\,\bigr\}\, $$
$$ \cup\,\bigl\{\, u_{\gamma }\, u_{c_0}\, u^{-1}_{\gamma c_0}\,\vert\, c_0 \in N \backslash M\,\bigr\}\,\cup\, 
\bigl\{\, u_{c_0}\, u_d\, u^{-1}_{c_0 d}\,\vert\, c_0 \in N \backslash M\,\bigr\}\, $$ 
$$ \cup\,\bigl\{\, [\, u_{c_0}\, ,\, u_x\, u_y\, u^{-1}_{xy}\, ]\,\vert\, c_0 \in N \backslash M\,\bigr\}\> . $$
\par\smallskip\noindent
Then there is a projection $\, J_G + U_F \longrightarrow U_F\, $ sending $\, J_{M , G}\, $ to 
$\, J_{N , F}\, $ and its restriction to $\, J_{M , G} + U_F\, $ shows that the map 
$\, K^J_2 ( J_{N , F}\, ,\, U_F\, ) \rightarrowtail K^J_2 ( J_{M , G}\, ,\, J_{M , G} + U_F )\, $ is injective in the first place by existence of a splitting. Then consider the inclusion 
$\, ( J_{M , G}\, ,\, J_{M , G} + U_F ) \subseteq ( J_{M , G}\, ,\, U_G )\, $. Dividing by $\, J_{M , G}\, $ gives  
$\, N \rtimes U_{F / N} \subseteq M \rtimes U_{G / M}\, $. Combining a projection
$\, U_{G / N} \rightarrow U_{F / N}\, $ with a retraction $\, M \searrow N\, $ gives a retraction 
$\, M \rtimes U_{G / M} \searrow N \rtimes U_{F / N}\, $ which is a homomorphism modulo 
$\, N\, $. The argument above applies to give that the image of an element of 
$\, K^J_2 ( J_{M , G}\, ,\, J_{M , G} + U_F )\, $ which becomes trivial in $\, K^J_2 ( J_{M , G}\, ,\, U_G )\, $ must be trivial modulo the image of 
$\, K^J_2 ( J_{M , G}\, ,\, J_{M , G} + U_N )\, $. Now there is a natural map 
$\, U_{M , G} \longrightarrow U_M\, $ sending an element of the standard basis of $\, U_{M , G}\, $ to the corresponding basis element in $\, U_M\, $ (under evaluation). Clearly the kernel of this map is a cone. This implies that
$\, K^J_2 ( J_{M , G}\, ,\, J_{M , G} + U_N )\, =\, K^J_2 ( J_M\, ,\, J_M + U_N )\, $. Then consider the basis of $\, J_M + U_N\, $ given by the union of sets
\smallskip
$$ \bigl\{\, u_{\gamma }\,\bigr\}\,\cup\,
\bigl\{\, u_{\gamma }\, u_{c_0}\, u^{-1}_{\gamma c_0}\,\vert\, c_0 \in N \backslash M\,\bigr\}\,\cup\,
\bigl\{\, u_{c_0}\, u_d\, u^{-1}_{c_0 d}\,\vert\, c_0 \in N \backslash M\,\bigr\}  $$
\par\smallskip\noindent
Then there is a projection 
$\, J_M + U_N \longrightarrow U_N\, $ and the kernel of this map is a cone, so that 
$\, K^J_2 ( J_M\, ,\, J_M + U_N )\, =\, K^J_2 ( J_N\, ,\, U_N )\, $ showing that our map is injective by assumption that the image of $\, K^J_3 ( N )\, $ in $\, K^J_3 ( N , F )\, $ is trivial. 
\par\medskip\noindent
The case $\, n \geq 4\, $ follows from the case $\, n = 3\, $ by putting $\, ( N , F ) = ( J_{N , F}\, ,\, U_F )\, $ and noting that $\, K^J_3 ( J_{N , F} ) = 0\, $ by the Cone Lemma (see below) because $\, J_{N , F}\, $ is a free group.
\par\medskip\noindent
Next, consider the case $\, K^J_2 ( N ) \rightarrow K^J_2 ( M )\, $ where $\, M = M_0\cdot N\, $.
The problem is to find a projection
\smallskip
$$ J_M\>\cap\>\bigl[\, U_M\, ,\, U_M\, \bigr]\>\longrightarrow\> J_N\>\cap\> \bigl[\, U_N\, ,\, U_N\,\bigr]  $$
\par\medskip\noindent
sending $\, [\, B_1\, ,\, u_c\, ]\, $ to $\, [\, B^N_1\, ,\, u_\gamma \, ]\, $ (or at least a substantial part). One can assume that $\, N\, $ contains an element $\,\mu\, $ of infinite order which is central for $\, M\, $, by injectivity of $\, K^J_2 ( N )\rightarrow K^J_2 ( N \times\mathbb Z )\, $. Each element $\, c\, $ can be written uniquely as $\, c\, =\, c_0\, \gamma\, $ with $\, c_0 \in M_0\, $ and 
$\,\gamma \in N\, $, and each element $\, b\, $ can be written uniquely as $\, b\, =\, \beta\, b_1\, $ with 
$\,\beta \in N\, ,\, b_1 \in M_0\, $. Then modulo $\, [\, J_M\, ,\, J_M\, ]\, $ each element of 
$\, [\, B_1\, ,\, u_c\, ]\, $ factors into expressions of the types
\smallskip
$$ \bigl[\, u_a\, u_{b_1}\, u^{-1}_{ab_1}\, ,\, u_\gamma\,\bigr]\> ,\> 
\bigl[\, u_\alpha\, u_b\, u^{-1}_{\alpha b}\, ,\, u_{c_0}\, \bigr]\> ,\>
\bigl[\, u_a\, u_\beta\, u^{-1}_{a\beta }\, ,\, u_{c_0}\, \bigr]\> ,\> 
\bigl[\, u_a\, u_\beta\, u^{-1}_{a\beta }\, ,\, u_\gamma\,\bigr]    $$
\par\smallskip\noindent
with $\, \alpha , \beta , \gamma \in N\, $. Write 
\smallskip
$$ u_a\> u_\beta\> u^{-1}_{a\beta }\> =\> \bigl(\, u_a\, u_{a^{-1}}\,\bigr)\, u^{-1}_{a^{-1}}\,
\bigl(\, u_\beta\, u_{{\beta }^{-1} a^{-1}}\, u^{-1}_{a^{-1}}\,\bigr)\, u_{a^{-1}}\,
\bigl(\, u^{-1}_{{\beta }^{-1} a^{-1}}\, u^{-1}_{a\beta }\,\bigr)   $$
$$ \equiv\> \bigl(\, u_a\, u_{a^{-1}}\,\bigr)\,\bigl[\, u_a\,
\bigl(\, u_\beta\, u_{{\beta }^{-1} a^{-1}}\, u^{-1}_{a^{-1}}\,\bigr)\, u^{-1}_a\,\bigr]\, 
\bigl(\, u^{-1}_{{\beta }^{-1} a^{-1}}\, u^{-1}_{a\beta }\,\bigr)     $$
\par\medskip\noindent
modulo $\, [\, J_M\, ,\, J_M\, ]\, $. Then modulo $\, [\, J_M\, ,\, J_M\, ]\, $ and terms of the form
$\, [\, u_a\, u_{a^{-1}}\, ,\, u_{c_\gamma }\, ]\, $ with $\, c_\gamma \in N\, $ or $\, c_\gamma \in M_0\, $ one gets expressions of the form $\, [\, u_\alpha\, u_b\, u^{-1}_{\alpha b}\, ,\, u_{c_\gamma }\, ]\, $ which can be factored modulo $\, [\, J_M\, ,\, J_M\, ]\, $ into expressions
\smallskip
$$ \bigl[\, u_\alpha\, u_{b_1}\, u^{-1}_{\alpha b_1 }\, ,\, u_c\,\bigr]\> ,\>
\bigl[\, u_\alpha\, u_\beta\, u^{-1}_{\alpha\beta }\, ,\, u_{c_0}\,\bigr]\> ,\>
\bigl[\, u_\alpha\, u_\beta\, u^{-1}_{\alpha\beta }\, ,\, u_\gamma\,\bigr]\>    $$
\par\smallskip\noindent
with $\, b_1\, ,\, c_0\, \in\, M_0\, $. We may assume that $\, a^2\,\neq\, 1\, $ and $\, ( a\beta )^2\,\neq\, 1\, $ in this decomposition, otherwise write $\, u_a\> u_\beta\> u^{-1}_{a\beta }\> =\> $ 
\smallskip
$$ \bigl(\, u_a\, u^{-1}_{a\mu }\, u^{-1}_\mu\,\bigr)\>
u_{{\mu }^{-1}}\>\bigl(\, u_{a\mu }\, u_\beta\, u^{-1}_{a\beta\mu }\,\bigr)\>
\bigl(\, u_{a\beta\mu }\, u^{-1}_{a\beta }\, u^{-1}_\mu\, \bigr)\> u^{-1}_{{\mu }^{-1}}\,
\bigl(\, u_{{\mu }^{-1}}\, u_\mu\, \bigr)\> ,    $$
\par\smallskip\noindent
so that modulo $\, [\, J_M\, ,\, J_M\, ]\, $ the element 
$\, [\, u_a\> u_\beta\> u^{-1}_{a\beta }\, ,\, u_{c_\gamma }\, ]\, $ can be replaced by elements 
$\, [\, u_{{\alpha }'}\, u_{b'}\, u^{-1}_{{\alpha }' b'}\, ,\, u_{c_{\gamma }'}\, ]\, $ and
$\, [\, u_{a'}\, u_{{\beta }'}\, u^{-1}_{a'{\beta }'}\, ,\, u_{c_{\gamma }'}\, ]\, $ with $\, ( a' )^2\,\neq\, 1\, $ and 
$\, ( a' {\beta }' )\,\neq\, 1\, $. Choose an order on the set of elements of $\, M\, $ inducing an order of the standard basis of $\, U_M\, $, Changing the basis to 
\smallskip
$$ \bigl\{\, u_c\, u_{c^{-1}}\,\vert\, c\, <\, c^{-1}\, ,\, c\,\not\in\, N\, ,\, c\,\not\in\, M_0\,\bigr\}\>\cup\>
\bigl\{\, u_{d_0}\,\vert\, d_0\,\in\, M_0\,\bigr\}\>\cup\> $$
$$ \bigl\{\, u_d\,\vert\, d\,\geq\, d^{-1}\, ,\, d\,\not\in\, N\, ,\, d\,\not\in\, M_0\,\bigr\}\>\cup\>
\bigl\{\, u_\gamma\, \vert\, \gamma\, \in\, N\,\bigr\}\>     $$
\par\smallskip\noindent
define an order on this new basis by the condition that any element of the first set is larger than any element of the other three, any element of the second set is larger than any element of the third and the fourth, and any element of the third set is larger than any of the last, combined with the induced orderings of the $\, {u_c }'s\, ,\, {u_{d_0}}'s\, ,\, {u_d}'s\, $ and $\, {u_\gamma }'s\, $. Then a basis for 
$\, [\, U_M\, ,\, U_M\, ]\, $ is given by the set
\smallskip
$$ \bigl\{\, v^{\pm }_1\cdots v^{\pm }_n\,\bigl[\, u\, ,\, v\,\bigr]\, v^{\mp }_n\cdots v^{\mp }_1\,\bigr\}   $$ 
\par\smallskip\noindent
where $\, u<v\> ,\> v_1\leq\cdots v_n\leq v\, $ and $\, u\, ,\, v\, ,\, v_1\, , \cdots ,\, v_n\, $ are taken from the altered basis of $\, U_M\, $ above. In particular, each element $\, [\, u_\gamma\, ,\, u_{a^{-1}}\, u_a\, ]\, $ is represented by the combination
\smallskip
$$ \bigl(\, u_\gamma\,\bigl[\, u_{a^{-1}}\, ,\, u_a\, u_{a^{-1}}\, ]\, u^{-1}_\gamma\,\bigr)\>
\bigl[\, u_\gamma\, ,\, u_a\, u_{a^{-1}}\,\bigr]\> \bigl[\, u_a\, u_{a^{-1}}\, ,\, u_{a^{-1}}\,\bigr]\>  . $$
\par\smallskip\noindent
Let $\, b\, =\, b_1 \in M_0\, $ and assume that either $\, a\, $ or $\, ab\, $ corresponds to a basis element of the first set in the basis of $\, U_M\, $. Then replacing 
$\, [\, u_\gamma\, ,\, u_a\, u_{b_1}\, u^{-1}_{ab_1}\, ]\, $ modulo $\, [\, J_M\, ,\, J_M\, ]\, $ and elements of the form $\, [\, u_\gamma\, ,\, u_a\, u_{a^{-1}}\, ]\, $ by 
$\, [\, u_\gamma\, ,\, u^{\pm }_{a^{\pm }}\, u_{b_1}\, u^{\mp }_{(ab_1)^{\pm }}\, ]\, $ and writing this element in terms of the basis of $\, [\, U_M\, ,\, U_M\, ]\, $ one finds that the latter may be changed according to the rule
\smallskip
$$  u^{\pm }_{a^{\pm }}\, u^{\mp }_{(ab_1)^{\pm }}\,\bigl[\, u_\gamma\, ,\, u_{b_1}\,\bigr]\,
u^{\pm }_{(ab_1)^{\pm }}\, u^{\mp }_{a^{\pm }}\>\longrightarrow\>
\bigl[\, u_\gamma\, ,\, u_a\, u_{b_1}\, u^{-1}_{ab_1}\,\bigr]   $$
\par\smallskip\noindent
where the signs on the left hand side are chosen so that the corresponding element $\, u_{a^{\pm }}\, $ resp. $\, u_{(ab_1)^{\pm }}\, $ is contained in one of the last three subsets of the basis of $\, U_M\, $ (and at least one of them in the third). Note that modulo $\, [\, J_M\, ,\, J_M\, ]\, $ one can assume 
$\, a\, <\, ab_1\, $ by changing $\, u_{b_1}\, $ to $\, u^{-1}_{{b_1}^{-1}}\, $ modulo expressions 
$\, [\, u_c\, ,\, u_{b_1}\, u_{{b_1}^{-1}}\, ]\, $ if necessary. Also the order may be chosen such that 
$\, a\, <\, a b_1\, $ implies $\, a^{\pm} < (a b_1)^{\pm}\, $. Making this change will ensure that the (new) basis elements $\, [\, u_\gamma\, ,\, u_a\, u_{b_1}\, u^{-1}_{ab_1}\, ]\, $  are contained in a complementary subspace to the space generated by conjugates of the 
$\, [\, u_\gamma\, ,\, u_a\, u_{a^{-1}}\, ]\, ,\, a\,\not\in\, N\, ,\, M_0\, $. If on the other hand $\, a\,\in\, M_0\, $ or $\, c\,\not\in\, N\, ,\, M_0\, ,\, c\,\geq\, c^{-1}\, $ one easily checks that the elements 
$\, [\, u_c\, ,\, u_a\, u_{b_1}\, u^{-1}_{ab_1}\, ]\, $ are in the complement of the space of conjugates of 
$\, [\, u_\gamma\, ,\, u_a\, u_{a^{-1}}\, ]\, ,\, a\,\not\in\, N\, ,\, M_0\, $ (for the new basis). If 
$\, c\,\not\in\, N\, ,\, M_0\, ,\, c\,<\, c^{-1}\, $ replacing $\, [\, u_c\, ,\, u_a\, u_{b_1}\, u^{-1}_{ab_1}\, ]\, $ modulo $\, [\, J_M\, ,\, J_M\, ]\, $ by 
$\, u^{-1}_{c^{-1}}\, [\, u_a\, u_{b_1}\, u^{-1}_{ab_1}\, ,\, u_{c^{-1}}\, ]\, u_{c^{-1}}\, $ (all the while we are assuming $\ a\, <\, ab_1\, $ so that $\, ab_1\,\not\in\, N\, $) this element also will be orthogonal to
$\, \langle\, [\, u_\gamma\, ,\, u_a\, u_{a^{-1}}\, ]\,\rangle\, $. The same is true for all other types of expressions given by our reduction modulo $\, [\, J_M\, ,\, J_M\, ]\, $ above which are basically
$\, [\, u_\alpha\, u_\beta\, u^{-1}_{\alpha\beta }\, ,\, u_{c_0}\, ]\, $ and 
$\, [\, u_\alpha\, u_\beta\, u^{-1}_{\alpha\beta }\, ,\, u_\gamma\, ]\, $. One now divides 
$\, [\, U_M\, ,\, U_M\, ]\, $ by all conjugates of the basis elements 
$\, [\, u_\gamma\, ,\, u_a\, u_{a^{-1}}\, ]\, ,\, a\, <\, a^{-1}\, ,\, a\,\not\in\, N\, ,\, M_0\, $ and
$\, [\, u_{a^{-1}}\, ,\, u_a\, u_{a^{-1}}\, ]\, =\, [\, u_{a^{-1}}\, ,\, u_a\, ]\, $ with 
$\, a\, <\, a^{-1}\, ,\, a\,\not\in\, N\, ,\, M_0\, $. Note that both elements can be seen as part of a basis for 
$\, J_M\, $ by changing the basis elements $\, (\, u_{\gamma a}\, u_{a^{-1}}\, u^{-1}_\gamma\, )\, $ of 
$\, B_1\, $ to 
\smallskip
$$    \bigl[\, u_\gamma\, ,\, u_a\, u_{a^{-1}}\, \bigr]\> =\> \bigl(\, u_\gamma\, u_a\, u^{-1}_{\gamma a}\,\bigr)
\>\bigl(\, u_{\gamma a}\, u_{a^{-1}}\, u^{-1}_\gamma\,\bigr)\>\bigl(\, u^{-1}_{a^{-1}}\, u^{-1}_a\,\bigr)   $$
\par\smallskip\noindent
and 
$\, (\, u_{a^{-1}}\, u_a\, )\, $ to $\, [\, u_{a^{-1}}\, u_a\, ]\, $ for $\, a\, <\, a^{-1}\, ,\, a\,\not\in\, N\, ,\, M_0\, $. 
We claim that our projection maps $\, [\, J_M\, ,\, J_M\, ]\,\subset\, [\, U_M\, ,\, U_M\, ]\, $ to itself. To see this one notes that modulo $\, [\, J_M\, ,\, J_M\, ]\, $ arbitrary $\, U_M$-conjugates of the elements 
$\, [\, u_\gamma\, ,\, u_a\, u_{a^{-1}}\, ]\, $ and $\, [\, u_{a^{-1}}\, ,\, u_a\, ]\, $ are mapped to zero, since up to commutators with elements $\, u_d\, u_{d^{-1}}\, ,\, d\, <\, d^{-1}\, ,\, d\,\not\in\, N\, ,\, M_0\, $ and 
$\, d\, >\, a\, $ all such conjugates appear as $\, [\, U_M\, ,\, U_M\, ]$-conjugates of some basis element 
$\, v^{\pm}_1\cdots v^{\pm}_n\, [\, u\, ,\, v\, ]\, v^{\mp}_n\cdots v^{\mp}_1\, $ where $\, [\, u\, ,\, v\, ]\, $ is of the form $\, [\, u_\gamma\, ,\, u_a\, u_{a^{-1}}\, ]\, $ or $\, [\, u_{a^{-1}}\, ,\, u_a\, u_{a^{-1}}\, ]\, $. The commutators $\, [\, u_a\, u_{a^{-1}}\, ,\, u_d\, u_{d^{-1}}\, ]\, $ and their conjugates however are orthogonal to our kernel and hence invariant under projection. Now suppose given an element of 
$\, [\, J_M\, ,\, J_M\, ]\, $ expressed in terms of the basis of $\, [\, U_M\, ,\, U_M\, ]\, $, say 
$\, x\, =\, x^{\pm}_1\cdots x^{\pm}_k\, $. If some $\, x_j\, $ is in the kernel it is some $\, U_M$-conjugate of 
$\, [\, u_\gamma\, ,\, u_a\, u_{a^{-1}}\, ]\, $ for example. Then also 
$\, (\, x^{\mp}_k\cdots x^{\mp}_{j+1}\, x^{\pm}_j\, x^{\pm}_{j+1}\cdots x^{\pm}_k\, )\,$ is a conjugate of the basis element $\, [\, u_\gamma\, ,\, u_a\, u_{a^{-1}}\, ]\, $ of $\, J_M\, $ by some $\, u\,\in\, U_M\, $. Changing the basis of $\, J_M\, $ by conjugation with $\, u\, $ then tells us that since 
$\, x\,\in\, [\, J_M\, ,\, J_M\, ]\, $ the expression 
$\, x^{\pm }_1\cdots x^{\pm}_{j-1}\, x^{\pm}_{j+1}\cdots x^{\pm}_k\, $ must contain an inverse of 
$\, x^{\pm}_j\, $ which therefore since the expression inbetween is in $\, J_M\, $ corresponds to some other $\, U_M$-conjugate of $\, [\, u_\gamma\, ,\, u_a\, u_{a^{-1}}\, ]\, $ which is shifted by an element from $\, J_M\, $, but since all such conjugates are divided out, the whole commutator in 
$\, [\, J_M\, , J_M\, ]\, $ corresponding to the intermediate terms is divided out and the resulting expression is again in $\, [\, J_M\, ,\, J_M\, ]\, $. Going on with this procedure until there are only basis elements in the complement of the kernel left shows that $\, [\, J_M\, , J_M\, ]\, $ maps to itself. Modulo 
$\, [\, J_M\, , J_M\, ]\, $ we have reduced our trivial element to expressions of the form 
$\, [\, u_a\, u_b\, u^{-1}_{ab}\, ,\, u_c\, ]\, $ where either $\, b\,\in\, M_0\, $ or $\, c\,\in\, M_0\, $ or 
$\, a\, ,\, b\, ,\, c\,\in\, N\, $. Consider the projection $\, J_M\rightarrow J_N\, $ given by 
\smallskip
$$  u_{a_0\alpha }\> u_{\beta b_1}\> u^{-1}_{a_0\alpha\beta b_1}\quad\longmapsto\quad
  u_\alpha\> u_\beta\> u^{-1}_{\alpha\beta }   $$
\par\smallskip\noindent
with $\, a_0\, $ and $\, b_1\, $ in $\, M_0\, $. It maps the  elements 
$\, [\, u_a\> u_{b_1}\> u^{-1}_{ab_1}\> ,\, u_c\, ]\, $ and $\, [\, u_a\> u_b\> u^{-1}_{ab}\> ,\, u_{c_0}\, ]\, $ to zero, so that modulo $\, [\, J_N\, , J_N\, ]\,$ our element is contained in $\, [\, B^N_1\,,\, u_\gamma\, ]\, $ and the map $\, K^J_2 ( N )\rightarrow K^J_2 ( M )\, $ is injective  \qed 
\par\bigskip\noindent
{\bf Remark.}\quad  The product of two product type inclusions need not be a product type inclusion itself. Nonetheless $\,  N  \subseteq  M \, $ is preexcisive whenever there exists 
$\,  M  \subseteq  M_1\, $ such that the inclusion $\,  N  \subseteq  M_1\, $ factors into a series of product type inclusions. 
\par\medskip\noindent
The estimate for the kernel of $\, K^J_2 ( N , F ) \rightarrow K^J_2 ( N , G )\, $ can be improved under mild restrictions. If the union of the $d$-commutants is contained in $\, N\, $ and if the section $\, s\, $ satisfies that its curvature coefficients $\, c_{z_0 \overline\zeta\, ,\, x_0 \overline\xi}\, $ are contained in the ($G$-normal) commutator subgroup $\, N_0\, =\, F_0\, [ N , G_{0 , 0} ]\, F^{-1}_0\, $, a similar proof as above gives injectivity modulo the image of $\, K^J_2 ( N )\, +\, K^J_2 ( N_0 , F )\, $, so both results taken together render injectivity modulo the intersection of the images of 
$\, K^J_2 ( N , F_0 )\, $ and $\, K^J_2 ( N )\, +\, K^J_2 ( N_0 , F )\, $. This sharper estimate may prove useful in certain contexts, for example when $\, F_0\, $ is large. The proof is a bit more extensive than the one given here, but since the main line of argument is very similar, the details are omitted. 
The assumption about $\, G_0\, $ containing the union of the $d$-commutants is necessary, there are counterexamples when this is not the case.

\end{document}